\documentclass{rwth}
\usepackage{stix}			
\usepackage{algorithm}		
\usepackage[noend]{algpseudocode}
\usepackage{tikz}
\usepackage{pgfplots}
\usepackage{pgfplotstable}
\usepackage{subfig}
\pgfplotsset{compat=1.8}		
\usepackage{diagbox}
\usepackage{hhline}
\def\aposteriori{\emph{a posteriori\/}}
\newcommand{\K}{\mathbb{K}}

\newcommand{\norm}[1]{\|{#1}\|}
\newcommand{\uN}{u_N^\ast} 
\newcommand{\uast}{u^\ast}
\DeclareMathOperator*{\argmax}{arg\,max}
\newcommand{\taucon}{\tau_{\textrm{con}}}
\newcommand{\cB}{\mathcal{B} } 
\newcommand{\cD}{\mathcal{D} } 
\newcommand{\cL}{\mathcal{L} } 
\newcommand{\cN}{\mathcal{N} } 
\newcommand{\cO}{\mathcal{O} } 
\newcommand*\diff{\mathop{}\!\mathrm{d}}	
\DeclareMathOperator{\spn}{span}
\newcommand{\alphaLB}{\ensuremath{\alpha_\mathrm{LB}(\mu)}}
\newcommand{\LB}{\ensuremath{\mathrm{LB}}}

\usepackage{color}
\usepackage{xcolor}
\newcommand{\cby}[1]{{\color{blue}#1}}

\begin{document}
\title{Reduced Order Model Predictive Control for Parametrized \\ Parabolic Partial Differential Equations}
%
\author{Saskia Dietze}\address{Numerical Mathematics, RWTH Aachen University, Templergraben 55, 52056 Aachen, Germany; \email{dietze@igpm.rwth-aachen.de\ \&\ grepl@igpm.rwth-aachen.de }}
\author{Martin A. Grepl}\sameaddress{1}
%
%
\begin{abstract} 
Model Predictive Control (MPC) is a well-established approach to solve infinite horizon optimal control problems. Since optimization over an infinite time horizon is generally infeasible, MPC determines a suboptimal feedback control by repeatedly solving finite time optimal control problems. Although MPC has been successfully used in many applications, applying MPC to large-scale systems -- arising, e.g., through discretization of partial differential equations -- requires the solution of high-dimensional optimal control problems and thus poses immense computational effort. 

We consider systems governed by parametrized parabolic partial differential equations and employ the reduced basis (RB) method as a low-dimensional surrogate model for the finite time optimal control problem. The reduced order optimal control serves as feedback control for the original large-scale system. We analyze the proposed RB-MPC approach by first developing {\it a posteriori} error bounds for the errors in the optimal control and associated cost functional. These bounds can be evaluated efficiently in an offline-online computational procedure and allow us to guarantee asymptotic stability of the closed-loop system using the RB-MPC approach in several practical scenarios. We also propose an adaptive strategy to choose the prediction horizon of the finite time optimal control problem. Numerical results are presented to illustrate the theoretical properties of our approach.

\end{abstract}
%
%
\maketitle
\section*{Introduction}

Model Predictive Control (MPC) -- also called Receding Horizon Control (RHC) -- is an advanced control technique for linear and nonlinear systems. Since the solution of infinite horizon optimal control problems is, in general, computationally intractable, MPC computes an approximate controller by solving finite time horizon optimal control problems in a receding horizon manner. To this end, a finite time optimal control problem is solved over a particular time interval, only the first control input is applied to the system, and the process repeats with an updated initial condition over a shifted time horizon. We refer to the monographs~\cite{RM09,GP11} or survey papers~\cite{GPM89, MRRS00, Le11} for more details.

Even though MPC is widely used in industry~\cite{QB03,FPHG15}, its applicability to high-dimensional systems is limited due to the large compuational cost for solving the associated optimization problems. Thus, reduced order models (ROMs) are often used to lessen the computational effort~\cite{AV15, HK07, HWG06, HGW08, HLGG08, LRHYA14, MV17, THKK07, SFG07, OS09, UKJ12, SBB13, DEMC06}. In this case, the high-dimensional model is typically substituted by a low-dimensional surrogate model and the prediction and computation of the control is based on the ROM. Although the use of ROMs results in significant computational savings, it also introduces an approximation error in the computation of the control which may jeopardize the (asymptotic) stability of MPC. In \cite{Li88}, it was shown that controllers based on ROMs might stabilize only the reduced model but not the original plant. Assuming the control computed from the (original) high-dimensional model resulted in a stable MPC scheme, this property should be inherited by the reduced control computed from the ROM. Robustness of the MPC closed-loop system with respect to the ROM approximation error is thus a crucial ingredient.

The stability of MPC -- without model order reduction (MOR) -- has been analyzed in detail in the literature. Whereas it is known that the optimal controller of the infinite horizon optimal control problem always stabilizes the system, this property may be lost within the MPC scheme~\cite[Chapter 1]{RM09}. It is therefore crucial to ensure (asymptotic) stability of the controller, i.e., to guarantee that the solution trajectory converges to the desired state and stays close to it. In MPC, terminal constraints or terminal costs or a combination of both are often used to guarantee stability, see e.g.~\cite{MM93,FT13, WK03, JYH01, Ma13, LRHYA14}. However, also MPC schemes with neither terminal constraints nor costs have received considerable attention in the literature since they are easier to implement and less restrictive. Without additional terminal constraints or cost the length of the prediction horizon is crucial since the stability and performance of the controller depend on it.  In \cite{PN00}, it was shown that stability of a constrained time-discrete linear system can be guaranteed if the prediction horizon is chosen sufficiently large; the result was generalized to nonlinear problems in \cite{JH05}. However, the computational complexity grows rapidly with the length of the optimization horizon. Therefore, estimates regarding the minimal horizon that guarantees stability are of utmost interest. In \cite{GP08,Gr09,GPSW10,RA11}, an analysis of the minimal stabilizing horizon based on relaxed dynamic programming \cite{LR06, Ra06} was proposed. Here, the main idea is to ensure stability by estimating the degree of suboptimality of the MPC feedback law and subsequently -- since the degree depends on the prediction horizon length -- to estimate the minimal stabilizing horizon. These ideas have been applied to MPC for partial differential equations in~\cite{AGW10,AG12,AG13}.

Although ROMs has been widely used in combination with MPC as mentioned above, the stability of the combined ROM-MPC approach has received less attention~\cite{LRHYA14, HLGG08, AV15,  SBB13, DEMC06}. In \cite{AV15}, proper orthogonal decomposition (POD) was used to derive a reduced order model. The stability analysis is based on the relaxed dynamic programming principle~\cite{Gr09,GPSW10} and does not require terminal costs or constraints. However, it directly applies the analysis from~\cite{Gr09,GPSW10} to the reduced order model, thus showing stability of the ROM but not of the closed-loop systems, i.e., the reduced control applied to the high-dimensional system.  Results from robust MPC are used in~\cite{HLGG08} and~\cite{SBB13}. In \cite{HLGG08}, stability is ensured by designing the cost function parameters so as to satisfy a set of linear matrix inequality (LMI) conditions. Since the approach requires the solution of LMIs involving the full system matrices, it is not applicable to high-dimensional systems. The work in~\cite{SBB13} does not show asymptotic stability of the original but only of a set around the origin. The approach in~\cite{DEMC06} provides asymptotic stability of the origin but is limited to modal decomposition to generate the reduced order model. Also, the stability analysis is quite conservative since a global error bound between the MOR and the high-dimensional system instead of the actual trajectory is used. Finally, the reduced order model considered in~\cite{LRHYA14} is based on a general Petrov-Galerkin projection and not limited to a specific MOR technique. Here, asymptotic stability can be achieved despite the MOR error using a terminal cost and the definition of terminal sets. An error bound is given then by the solution to a scalar ordinary differential equation which is also used in the design of the MOR feedback controller thus accounting for the MOR error.

\subsection*{Main Idea and Contributions}

In this paper, we propose a reduced basis MPC (RB-MPC) approach for parametrized partial differential equations (PDEs) without terminal constraints or costs.  We briefly sketch the main idea and contributions -- precise definitions, formulations, and assumptions will follow in the next section. Let us consider the generic discrete time control system given by
\begin{equation} \label{eq:generic_system}
y(t^{k+1}) = f(y(t^k),u(t^{k+1}))
\end{equation}
with inital condition $y(t^0) = y_0$,  state $y^k = y(t^{k}) \in Y$, and control $u^k = u(t^k) \in U$, where $Y$ and $U$ are appropriate function spaces for the state and control, respectively. In MPC, we first wish to find a control sequence $(u^k)_{k = 1}^K$ which minimizes the finite horizon cost functional
\begin{align}
J(u) = \sum_{k=1}^K \ell(y_u^k(y_0),u^{k})
\end{align}
with running cost $l(\cdot,\cdot): Y \times U \to \xR^+_0$ and prediction horizon $K \in \xN$. Here, $y_u^k(y_0)$ denotes the state solution trajectory for a given initial condition $y_0$ and control sequence $(u^k)_{k \in \xN}$. We denote the optimal control sequence as $u^{*,k}$ and the associated finite horizon optimal value function as $J^\ast(\cdot) : Y \to \xR_0^+$, where $J^\ast(y_{0}) = \min_{u\in U}  J(u)$. In MPC, we then apply only the first control $u^{*,1}$ to the system, compute the state $y_{u^{*}}^{k+1}(y_0)$, shift the horizon by one timestep, and solve the finite horizon problem again with the newly computed state as initial condition. 

Under certain conditions, asymptotic stability of the MPC feedback law can be shown if
\begin{align} \label{eq:stab_check}
J^\ast(y_0) \ge J^\ast(y_{u^\ast}^1) +\omega \ l(y_0, u^{\ast, 1})
\end{align} 
holds for all $y_0 \in Y$ and some $\omega\in (0,1]$. Here, $\omega$ is called the suboptimality degree. The optimal value function $J^{\ast}$ depends on the prediction horizon $K$ and it can be shown that -- if the prediction horizon is chosen sufficiently large -- we have $\omega\in (0,1]$, i.e., stability of the MPC controller can be guaranteed. Unfortunately, checking condition~\eqref{eq:stab_check} requires the solution of two finite horizon optimal control problems.  If system~\eqref{eq:generic_system} is given by a (discretized) partial differential equation of (typically) very large spatial dimension $\cal N$, solving for the optimal control $u^{*,1}$ in each MPC loop is computationally expensive, testing condition~\eqref{eq:stab_check} multiple times for increasing prediction horizons $K$ becomes prohibitive.

To this end, we introduce an RB approximation  $y_N(t^{k}) \in Y_N \subset Y$ to the full-order solution $y(t^{k}) \in Y$ given by
\begin{equation} \label{eq:generic_system_red}
y_N(t^{k+1}) = f_N(y_N(t^k),u_N(t^{k+1}))
\end{equation}
with inital condition $y_N(t^0) = y_{N,0}$. Here, we assume that the control $u_N(t^k)$ lies in $U$ and that the RB space $Y_N$ satisfies $\dim(Y_N) = N \ll {\cal N}$. We then employ the RB approximation to solve the finite time optimal control problems required in each MPC loop. Although generation of the RB approximation requires an up-front (offline) cost, we can then compute the reduced order optimal control $u_N^{*,1}$ efficiently online. To check whether the closed-loop system --- i.e. the original high-dimensional system~\eqref{eq:generic_system} with reduced order control $u_N^{*,1}$ --- is stable, we use condition~\eqref{eq:stab_check}. However, in order to avoid the solution of original (full-order) optimal control problems, we ({\it i}) derive lower and upper bounds for the error between the optimal value function $J^{\ast}$ and its associated reduced order approximation $J_N^{\ast}$, and ({\it ii}) replace $J^{\ast}$ appropriately by its lower and upper bounds. Since the reduced order value function $J_N^{\ast}$ and the associated {\it a posteriori} error bound can be evaluated using an offline-online computational procedure, we can thus also confirm the stability of the RB-MPC approach efficiently online. Finally,  based on the above ideas, we propose an adaptive algorithm to choose the minimal stabilizing prediction horizon of the RB-MPC closed-loop system and discuss possible convergence issues.

We note that although we employ the RB method as a surrogate model in this paper, other model reduction techniques like POD or balanced truncation could also be used. The {\it a posteriori} error bounds as well as the proposed RB-MPC approach remain valid if we replace the RB test and ansatz spaces with  the corresponoding spaces generated by POD or balanced truncation (assuming the reduced order model inherits the stability properties of the original system). In fact, even reduced spaces generated from data could be used in this framework.

After introducing the general problem statement in Section~\ref{sec:PPOCP_problem_definition}, we focus on the RB approximation and associated {\it a posteriori} error bound formulation in Section~\ref{sec:RB_Approx}. We propose the RB-MPC approach and analyze its stability in Section~\ref{sec:Stability_RB_MPC}, before assessing our approach with two numerical examples in Section~\ref{sec:NumExp}.

\section{General Problem Statement} \label{sec:PPOCP_problem_definition}

In this section we introduce the infinite horizon optimal control problem and the associated finite time horizon optimal control problem for the MPC approach. We consider systems governed by parametrized linear parabolic partial differential equations and introduce a finite element truth discretization of the continuous problem. Finally, we state the MPC algorithm and recall existing MPC stability results which are the basis for our analysis. Note that we restrict our attention to state feedback MPC in this work.

\subsection{Preliminaries}

Let $\Omega \subset \xR^d$, $d= 1,2,3$, be a bounded Lipschitz domain with boundary $\Gamma$ and let $Y_e$ with $H^1_0(\Omega)\subset Y_e \subset H^1(\Omega)$ be a Hilbert space over $\Omega$.\footnote{The subscript $e$ refers to "exact".} The associated inner product and norm are given by $(\cdot,\cdot)_{Y_e}$ and $\|\cdot\|_{Y_e}= \sqrt{(\cdot,\cdot)_{Y_e}}$, respectively. Additionally, we assume the norm $\|\cdot\|_{Y_e}$ is equivalent to the $H^1(\Omega)$-norm and denote the dual space of $Y_e$ by $Y_e'$. Next, we recall the Hilbert space  $W(0,T) = \{v \in L^2(0,T;Y_e)\ : \ v_t \in L^2(0,T;Y_e')\} $ for a fixed time $T$ with its standard inner product; for more information see e.g. \cite{QV08, Wl87, Ze90}. We also introduce the control space $U_e(T) = L^2(0,T;\xR^m)$, $m\in \xN$,  with inner product $(u,v)_{U_e} = \int_0^T (u(t),v(t))_{\xR^m} \xdif t$ and induced norm $\| \cdot \|_{U_e} = \sqrt{(\cdot, \cdot)_{U_e}}$; the dual space is denoted by $\left(U_e(T)\right)'$.  Furthermore, we define $ U_{e,ad}(T)$ as the non-empty convex subset of admissible controls, given by $U_{e,ad}(T) = \{ u_e \in U_e(T) \ :\ u_{a,e}(t) \le u_e(t) \le u_{b,e}(t) \} \subset U_e(T)$.
Here, $u_{a,e}(\cdot),u_{b,e}(\cdot)\in  L^2(0,T;\xR^m)$ with $u_{a,e}(t)\le u_{b,e}(t) \textrm{ f.a.a. } t\in [0,T]$ are given lower and upper bounds for the control input (hard constraints). Finally, let $\cD \subset \xR^P$ be a prescribed $P$-dimensional compact parameter set, i.e., an element $\mu \in \cD$ will serve as our $P$-tuple input parameter.

We now introduce the bilinear forms. For the sake of simplicity the first one is parameter independent and defined through $m(w,v) = (w,v)_{L^2(\Omega)}$ for all $ \ w,v \in L^2(\Omega)$. The second parameter-dependent bilinear form is given by $a(\cdot,\cdot;\mu) : Y_e \times Y_e \to \xR$. We assume that $a(\cdot,\cdot;\mu)$ is continuous, 
\begin{align}
	\gamma_e(\mu) = \sup_{w\in Y_e\backslash\{0\}} \sup_{v \in Y_e\backslash\{0\}} \dfrac{a(w,v; \mu)}{\|w\|_{Y_e}\|v\|_{Y_e}} \le \gamma_0 < \infty, \quad \forall \mu \in \cD,
	\end{align}
coercive,
	\begin{align}\label{eq:Exact_ellipticity}
	\alpha_e (\mu) = \inf_{v \in Y_e\backslash\{0\}} \dfrac{a(v,v;\mu)}{\|v\|_{Y_e}^2} \ge \alpha_0 >0 , \quad \forall \mu \in \cD,
  \end{align}	
and admits an affine representation, 
\begin{align}\label{eq:AffineRepr_a}
  	a(w,v; \mu) = \sum_{q=1}^{Q_a} \Theta^q_a(\mu) a^q(w,v),\quad \forall w,v\in Y_e, \quad \forall \mu \in \cD,
\end{align}
for some (preferably) small integer $Q_a$. The coefficient functions $\Theta^q_a :\ \cD \to \xR$ are continuous and depend on the parameter $\mu$, whereas the continuous bilinear forms $a^q$ do \textit{not} depend on $\mu$.  
Last, we also require a continuous bilinear form $b(\cdot,\cdot): \xR^m \times Y_e \to \xR$, defined by 
\begin{align}
 b(u,v) =  \sum_{i=1}^m b_i(v) u_i,
\end{align}
where $b_1, \dots,b_m$ are given bounded linear functionals on $L^2(\Omega)$. We assume the $b_i$ do not depend on the parameter; however, affine parameter dependence is readily admitted~\cite{GP05}.
 In addition, we define the continuous and linear operator $\mathcal{B}_e : U_e(T)\to L^2(0,T; Y_e')$ by 
\begin{align}
 \langle (\mathcal{B}_e u_e)(t), \cdot\rangle_{Y'_e,Y_e} = b(u_e(t), \cdot) = \sum_{i=1}^m b_i(\cdot) u_{e,i}(t). 
\end{align} 
Here, $ \langle \cdot, \cdot \rangle_{Y'_e,Y_e}$ denotes the dual pairing between $Y'_e$ and $Y_e$ and $u_e \in U_e(T)$ is the control input with time-dependent components $u_{e,i} \in L^2(0,T)$, $1\le i \le m$.
The associated adjoint operator (linear and bounded) is denoted by $\mathcal{B}_e^\star : L^2(0,T;Y_e) \to U_e(T) $, where we identify $\left( L^2(0,T;Y_e')\right)'$ with $L^2(0,T, Y_e)$ and $\left(U_e(T)\right)'$ with $U_e(T)$. It can be shown that, for a given $\phi\in L^2(0,T;Y_e)$, the adjoint operator $\mathcal{B}_e^\star \phi$ can be expressed as $(\mathcal{B}_e^\star \phi)_i(t) = b_i(\phi(t)) , \quad 1\le i \le m, \quad t\in [0,T]$~\cite{KG14}.
Finally, we require an observation subdomain $D\subset\Omega$, which is a measurable set with inner product $(u,v)_D= \int_D u v \xdif t$ and an associated semi-norm $|v|_D =\sqrt{(v, v)_{D}} $ for $u,v\in \xLtwo(\Omega)$. We shall assume that the system is LTI (linear time-invariant), this is the case when all linear and bilinear forms are time independent.   

We note that the affine parameter dependence of all of the parameter dependent quantities is crucial for the offline-online dependence and thus the online efficiency of our proposed approach. Although this assumption might seem restrictive, it naturally holds for many material and geometrical parametrizations~\cite{RHP2008}. In the nonaffine case, the Empirical Interpolation Method (EIM) may be used to obtain an approximately affine representation of the involved (bi-)linear forms~\cite{BMN+2004,GMN+2007}.

\subsection{Infinite Horizon Control Problem}

In this paper we consider parametrized parabolic partial differential equations of the form
\begin{align}
\label{eq:PDEConstraint}
\dfrac{\xdif}{\diff t} m(y_e(t),v) + a(y_e(t),v;\mu) = b(u_e(t),v), \quad \forall v\in Y_e, \quad \textrm{f.a.a. } t\in (0,\infty),
\end{align}
where $y_e(t)\in Y_e$ for a given parameter $\mu \in \cD$, control $u_e(t)\in \xR^m$, and initial condition $y_e(0) = y_{0,e}\in L^2(\Omega)$.\footnote{Often, we omit dependencies to simplify the notation, for instance $y_e = y_e(x,t;\mu)$. Whenever we want to point out the underlying control to a trajectory, we write $y_e(u_e)$.  
}  Given all assumptions in the previous section, \eqref{eq:PDEConstraint} is well-posed \cite{QV08}.
Next, we introduce the infinite horizon cost defined by 
\begin{align}\label{eq:Cost_func_infty}
 J_{e,\infty}(y_e, u_e;\mu) = \int_{0}^\infty \left[ \dfrac{\sigma_1}{2}  | y_e(t) - y_{d,e}(t;\mu)|^2_{D} 
	+  \dfrac{\lambda}{2} \| u_e(t) - u_{d,e}(t)\|^2_{R^m} \right] \diff t.  
\end{align}
The integral measures the deviation of the state $y_e$ and the control $u_e$ from the desired state $y_{d,e}(t;\mu)\in L^2(\Omega)$ and the desired control $u_{d,e}(t)\in \xR^m$, respectively. It either has a well-defined nonnegative finite value or it diverges to infinity. 

We assume that the (parameter-dependent) desired state $y_{d,e}(t;\mu)$ has an affine representation 
\begin{align}\label{eq:Affine_DesiredState_infty}
y_{d,e}(x,t;\mu) = \sum_{q=1}^{Q_{y_d}} \Theta_{y_d}^q(t;\mu) y^q_{d,e}(x) ,
\end{align}
where the coefficient functions $\Theta_{y_d}^q: [0,\infty)\times \cD \to \xR$ are parameter and (possibly) time-dependent, whereas the functions $ y^q_{d,e}\in L^2(\Omega)$ are parameter-independent. The control $u_{d,e}(t)$ is assumed to be parameter-independent. However, affine parameter dependence is readily admitted.

Our goal is to steer the system towards the desired trajectories by minimizing $J_{e,\infty}$ subject to the PDE constraint \eqref{eq:PDEConstraint}.  If $y_{d,e} \equiv 0$ this is called a stabilization problem, otherwise we solve a so called tracking problem.

Since infinite horizon optimal control problems are, in general, infeasible, we apply a receding horizon approach. We introduce the associated finite horizon control problem in the next section.
 
\begin{rmrk}
It is well-known that -- under certain assumptions -- optimal control problems with a large time horizon satisfy the so-called turnpike property, i.e., the optimal (state and control) solution trajectories stay exponentially close to the optimal solutions of the corresponding steady-state problem for most of the time. One of these assumptions -- which we do not require -- are time-independent or periodic in time desired states and controls. For more details on the turnpike property we refer e.g.\ to ~\cite{FAULWASSER2022367} and~\cite{doi:10.1137/130907239,doi:10.1137/18M1225811}.
\end{rmrk} 
 
\subsection{Finite Time Horizon Control Problem}

For simplicity - and because the system is LTI - we only present the optimal control system for the time interval $[0,T]$. We consider the parametrized optimal control problem 
\begin{equation}\label{eq:ExactOCP} \tag{$P_e$}
\begin{split}
 & \min_{y_e, u_e} J_e(y_e,u_e;\mu) \quad \textrm{s.t.} \quad (y_e,u_e)\in W(0,T)\times U_{e,ad}(T) \textrm{ solves } \\
& \dfrac{\diff}{\diff t} m(y_e(t),v) + a(y_e(t),v;\mu) = b(u_e(t),v), \quad \forall v\in Y_e, \quad \textrm{f.a.a. } t\in (0,T],
\end{split}
\end{equation}
with initial condition $y_e(0) = y_{0,e}\in L^2(\Omega)$. 
The quadratic cost functional $J_e(\cdot, \cdot;\mu) :  W(0,T)\times U_{e}(T) \to \xR$, is defined by 
\begin{equation}
J_e(y_e, u_e;\mu) = \dfrac{\sigma_1}{2} \int_0^T | y_e(t) - y_{d,e}(t;\mu)|^2_{D} \diff t 
	+ \dfrac{\sigma_2}{2} | y_e(T) - y_{d,e}(T;\mu)|^2_{D} 
	+  \dfrac{\lambda}{2} \| u_e - u_{d,e}\|^2_{U_e}.
\end{equation}
Note, that  $y_{d,e}(\mu)$ acts only locally on the subdomain $D$ (although it is defined on $\Omega$) due to the definition of $J_e$. The desired state $y_{d,e}(\mu)$ again satisfies the affine parameter dependence~\eqref{eq:Affine_DesiredState_infty} with $\Theta_{y_d}^q$ replaced by  $\Theta^q_{y_d \, | \, [0,T]}$, i.e., the restriction of $\Theta_{y_d}^q$ to $[0,T]$. 
The parameters $\lambda>0, \sigma_1,\sigma_2\ge 0 $  are responsible for the weighting between the deviation from the desired control and the desired state (for the whole time interval and/or the final time). They are chosen depending on the application and the desired performance of the optimal controller. 
Under the stated assumptions there exists a unique optimal solution $(y^\ast_e, u_e^\ast)$ to \eqref{eq:ExactOCP} \cite{Li71, Tr05}. 

To solve the optimal control problem we employ a Lagrangian approach. We define the parametrized Lagrangian functional using the parabolic PDE and the definition of the cost functional 
\begin{equation}
\begin{split}\label{eq:ParamLagrangianFunc}
\cL_e(\xi_e;\mu) := & J_e(y_e,u_e;\mu) + \int_{0}^T a(y_e(t), p_e(t);\mu) \diff t + \int_{0}^T  m( \dfrac{\diff}{\diff t} y_e(t), p_e(t)) \diff t \\
& - \int_0^T b(u_e(t),p_e(t)) \diff t + m(y_e(0) - y_{0,e}, p_e(0)),
\end{split} 
\end{equation}
where $\xi_e = (y_e,p_e,u_e)\in  W(0,T)\times W(0,T)\times U_e(T) $ and $p_e$ is the adjoint variable. 
Differentiation of the Lagrangian functional \eqref{eq:ParamLagrangianFunc} leads to the first-order optimality system consisting of the state equation, the adjoint equation and the optimality condition: Given $\mu\in \cD$, the optimal solution (to \eqref{eq:ExactOCP}) $(y_e^\ast, p_e^\ast, u_e^\ast)\in W(0,T)\times W(0,T)\times U_{e,ad}(T) $ satisfies
\begin{subequations}\label{eq:OCPSystem}
 	\begin{alignat}{2}
		 \dfrac{\diff}{\diff t} m(y_e^\ast(t),\phi) + a(y_e^\ast(t),\phi;\mu) &=b(u_e^\ast(t),\phi)  \quad 					
		 		&& \forall \phi \in Y_e \quad \textrm{f.a.a. } t\in (0,T], \label{eq:OCPSystem_StateA}\\
 		y_e^\ast(0) &= y_{0,e}, \label{eq:OCPSystem_StateB}\\
 		-\dfrac{\diff}{\diff t}m(\varphi, p_e^\ast(t)) + a(\varphi, p_e^\ast(t);\mu) &= \sigma_1(y_{d,e}(t;\mu) - y_e^\ast(t),   \varphi)_D   \quad 
 				&& \forall \varphi \in Y_e \quad \textrm{f.a.a. } t\in [0,T), \label{eq:OCPSystem_AdjA}\\
 		m(\varphi, p_e^\ast(T)) &=  \sigma_2(y_{d,e}(T;\mu) - y_e^\ast(T),\varphi)_D \quad 
 				&& \forall \varphi \in Y_e , \label{eq:OCPSystem_AdjB}\\
		 (\lambda (u_e^\ast - u_{d,e}) - \cB_e^\star p_e^\ast, \psi - u_e^\ast)_{U_e} &\ge 0 \quad 
		 		&& \forall \psi\in U_{e,ad}. \label{eq:OCPOptimalityInequality}
  \end{alignat}
\end{subequations}

These conditions \eqref{eq:OCPSystem} are necessary and sufficient for the optimality of $(y_e^\ast, u_e^\ast)$ \cite{Li71, Tr05}.

\begin{rmrk}
Regularization parameters often serve as design parameters which are tuned to achieve a desired performance of the optimal controller.  In terms of our reduced basis approach it is therefore convenient to consider these parameters as additional input parameters. This procedure allows us to efficiently design the optimal controller during the online phase of the algorithm. We will consider $\lambda$ as an additional parameter in the numerical experiments in Section \ref{sec:NumExp}. 
\end{rmrk}

\subsection{Truth Finite Element Approximation}

Since we cannot expect to solve the optimal control problem analytically, we apply a finite element method for the spatial discretization while the backward Euler method\footnote{It is also possible to apply high-order schemes like Crank-Nicolson or a general $\Theta$-scheme instead \cite{GrPhd}.}
is used for the temporal approximation. This yields the so called "truth" approximation. 
To this end, we divide the time interval $[0,T]$ in $K$ subintervals of length $\tau = T/K$ and define $t^k = k \tau$, $0\le k \le K$, and $\K=\{1,\dots , K\}$. 
Furthermore, let $Y\subset Y_e$ denote a finite element space of (typically) large dimension $\cN$. The finite dimensional space $Y$ inherits the inner product and induced norm from $Y_e$, i.e., $(\cdot,\cdot)_{Y}=(\cdot,\cdot)_{Y_e}$ and $\|\cdot\|_{Y} =\|\cdot\|_{Y_e}$ as well as the continuity and coercivity properties, i.e.,
\begin{align}
	\gamma(\mu) = \sup_{w\in Y\backslash\{0\}} \sup_{v \in Y\backslash\{0\}} \dfrac{a(w,v; \mu)}{\|w\|_{Y}\|v\|_{Y}} \le \gamma_e(\mu)  \le \gamma_0 < \infty, \quad \forall \mu \in \cD,
\end{align}
and
	\begin{align} \label{eq:truth_ellipticity}
	\alpha (\mu) = \inf_{v \in Y\backslash\{0\}} \dfrac{a(v,v;\mu)}{\|v\|_{Y}^2} \ge \alpha_e (\mu)  \ge \alpha_0 >0 , \quad \forall \mu \in \cD. 
	  \end{align}
 We define the operator $\cB : U \to (Y')^K$ by
\begin{align}
 \langle (\mathcal{B} u)^k, \cdot\rangle_{Y',Y} = b(u^k, \cdot) = \sum_{i=1}^m b_i(\cdot) u_{i}^k,  \quad k\in \K,
\end{align}
with $U=(\xR^m)^K$ being the discrete control space with inner product $(u,v)_U = \tau\sum_{k=1}^K (u^k,v^k)_{\xR^m}$ and norm $\|\cdot\|_{U} = \sqrt{(\cdot, \cdot)_U}$. A control is denoted by $u = (u^1, \dots, u^K), $ with $ u^k \in \xR^m$. 
We introduce the discretized admissible control set as $ U_{ad} = \{ u \in U \ :\ u^k \in U^k_{ad}, \ k\in \K\}$, where $U^k_{ad}= \{ z \in \xR^m : u_{a}^k \le z \le u_{b}^k \}$ and $u_{a}^k= u_{a,e}(t^k)\in \xR^m, \ u_{b}^k= u_{b,e}(t^k)\in \xR^m$. 
The corresponding truth optimal control problem now reads as follows
\begin{equation}\label{eq:FE_OCP}\tag{$P$}
\begin{split}
 & \min_{y, u} J(y,u;\mu) \quad \textrm{s.t.} \quad (y,u)\in Y^{K}\times U_{ad} \textrm{ solves } \\
& m(y^k,v) + \tau a(y^k,v;\mu) =m(y^{k-1},v) +  \tau b(u^k,v), \quad \forall v\in Y, \quad \forall k\in \K, 
\end{split}
\end{equation}
with initial condition $y^0 = y_0$. Here, $y_0$ is the $L^2$-projection of $y_{0,e}$. Moreover, $y^k$ denotes the truth solution at time $t^k$, i.e., $y^k=y(t^k)$, and $u_i^k\in \xR$ corresponds to the $i$-th control input at time $t^k$, i.e., $u_i^k= u_i(t^k)$. 
The discretized cost functional is given by $J(\cdot, \cdot;\mu) :  Y^{K}\times U\to \xR$, 
\begin{align}
\begin{split}\label{eq:FE_cost_functional}
J(y, u;\mu) &= \dfrac{\sigma_1}{2} \tau \sum_{k=1}^K | y^k - y_{d}^k(\mu)|^2_{D} 
	+ \dfrac{\sigma_2}{2} | y^K - y_{d}^K(\mu)|^2_{D} 
	+  \dfrac{\lambda}{2} \tau \sum_{k=1}^K  \| u^k - u_{d}^k\|^2_{\xR^m}\\
	&= \dfrac{\sigma_1}{2} \tau \sum_{k=1}^{K-1} | y^k - y_{d}^k(\mu)|^2_{D} 
	+\left( \dfrac{\sigma_1}{2} \tau	+ \dfrac{\sigma_2}{2}\right) | y^K - y_{d}^K(\mu)|^2_{D} 
	+  \dfrac{\lambda}{2} \tau \sum_{k=1}^K  \| u^k - u_{d}^k\|^2_{\xR^m}  
	. 
\end{split}
\end{align}
Here, $u_d^k = u_{d,e}(t^k)\in \xR^m$ and $y_{d}^k(\mu)\in Y$ is the $L^2$-approximation of $y_{d,e}(t^k;\mu)\in L^2(\Omega)$.

%

Analogously to the procedure in the previous section, we employ a Lagrangian approach to obtain the first-order optimality system: Given $\mu\in \cD$, the optimal solution (to \eqref{eq:FE_OCP}) $(y^\ast, p^\ast, u^\ast)\in Y^{K}\times Y^{K}\times U_{ad} $ satisfies
\begin{subequations}\label{eq:FE_OCPSystem}
 	\begin{alignat}{2}
		 m(y^{\ast,k}-y^{\ast,k-1},\phi) + \tau a(y^{\ast,k},\phi;\mu) &= \tau b(u^{\ast,k},\phi)  \quad 					
		 		&& \forall \phi \in Y, \quad \forall k\in \K, \\
 		y^{\ast,0} &= y_0, \\
 		m(\varphi, p^{\ast,k}-p^{\ast,k+1}) + \tau a(\varphi, p^{\ast,k};\mu) &= \tau  \sigma_1(y_{d}^k(\mu) - y^{\ast,k},   \varphi)_D   \quad 
 				&& \forall \varphi \in Y, \quad \forall k\in \K\backslash\{K\} , \\
 				m(\varphi, p^{\ast,K}) + \tau a(\varphi, p^{\ast,K};\mu) &= (\tau  \sigma_1 +\sigma_2)(y_{d}^K(\mu) - y^{\ast,K},   \varphi)_D   \quad 
 				&& \forall \varphi \in Y,\\
		 (\lambda (u^{\ast} - u_{d}) - \cB^\star p^{\ast}, \psi - u^{\ast})_{U} &\ge 0 \quad 
		 		&& \forall \psi\in U_{ad}. \label{eq:FE_OCPOptimalityInequality}
  \end{alignat}
\end{subequations}
The adjoint operator $\cB^\star : Y^K \to U $ is given, for $\phi\in Y^K$, by $(\cB^\star\phi)_i^k = b_i(\phi^k)$. 
Note, that we use the same trial and test spaces in the discretization. 

%

Given the definition of the cost functional, we introduce the (running cost) function $l(\cdot, \cdot;\mu) $ measuring the deviation of the state and the control from the desired trajectories at a time $t^k$.  The function $l(\cdot,\cdot;\mu): Y\times \xR^m \to \xR^+_0$ is defined as
\begin{align}
l(y^k,u^{k+1};\mu) = \dfrac{\sigma_1}{2} \tau| y^k - y_{d}^k(\mu)|^2_{D} +   \dfrac{\lambda}{2} \tau  \| u^{k+1} - u_{d}^{k+1}\|^2_{\xR^m}. 
\end{align}

\begin{rmrk}
The system \eqref{eq:FE_OCPSystem} was derived by first discretizing the cost functional and state equation and then optimizing, i.e., the first-discretize-then-optimize (FDTO) approach. However, for the time discretization scheme used here we obtain the same optimality system using the first-optimize-then-discretize (FOTD) approach, i.e., both approaches commute. We note that this is not always the case and depends on the specific time discretization scheme, see e.g.\cite{doi:10.1137/100819333} for more details.
\end{rmrk}

\begin{rmrk}
The discrete infinite horizon problem corresponding to the continuous case \eqref{eq:PDEConstraint} and \eqref{eq:Cost_func_infty} attempts to minimize the infinite horizon cost functional
\begin{align}\label{eq:FE_Cost_func}
J_{\infty}(y,u;\mu)= \sum_{k=0}^\infty \left[ \dfrac{\sigma_1}{2} \tau | y^k - y_{d}^k(\mu)|^2_{D} 
	+  \dfrac{\lambda}{2} \tau  \| u^{k+1} - u_{d}^{k+1}\|^2_{\xR^m} \right]
\end{align}
subject to the system dynamics 
\begin{align}\label{eq:FE_PDEConstraint}
m(y^k,v) + \tau a(y^k,v;\mu) =m(y^{k-1},v) +  \tau b(u^k,v), \quad \forall v\in Y, \quad \forall k\in \xN,
\end{align}
with initial condition $y(0)=y_{0}.$
\end{rmrk}

\subsection{Model Predictive Control Algorithm}\label{sec:MPC}

The classical MPC algorithm computes a feedback control for the plant by solving optimal control problems for a finite time interval and moving the horizon forward: 
At every sampling instance $t^k$, we solve the optimal control problem \eqref{eq:FE_OCP} over the finite interval $[t^k,t^k+T]$  with the state measurement $y(t^k)$ as initial condition. 
Here, $T$ is the finite optimization horizon corresponding to $K = T / \tau $ discrete time-steps. 
Subsequently, we apply the computed optimal control to the system from time $t^k$ until $t^{k+1}= t^k+\tau$. Next, we take a new measurement of the state to update the initial condition and repeat the process for the shifted interval $[t^{k+1}, t^{k+1}+T]$. 
The resulting controller in feedback form $\kappa:Y \to U_{ad}$ with $\kappa(y(t^k))=u^{\ast,1}$ allows us to react to the current deviation of the state $y$ from $y_d$ at time $t^k$. The MPC algorithm is summarized in Algorithm \ref{alg:MPC}.


\begin{algorithm}[t]
\caption{Model Predictive Control}\label{alg:MPC}
\begin{algorithmic}[1]
\State Given: time step $\tau >0$; 
finite optimization horizon $T = K\tau, \ K\in \xN$, parameter $\mu\in \cD$, initial condition $ y_{\kappa}(t^0)= y_0$. 
\State Set $\tilde{y}_d^k(\mu) = y_d^k(\mu)$, $\tilde{u}_d^k = u_d^k$, $k = 1,2,\ldots$.
\For{$i=0,1,2,\dots$} 
\State Set $t^i=i\tau$, $y_0 = y_{\kappa}(t^i)$, and $y_d^k(\mu) = \tilde{y}_d^{k+i}(\mu)$, $u_d^k = \tilde{u}_d^{k+i}$, $k = 1,\ldots, K$.
\State Solve the optimal control problem \eqref{eq:FE_OCP}  
  with initial condition $y(0) = y_{0}$ and obtain the 
  optimal solution $(y^\ast, u^\ast)$ with optimal control $u^\ast=\left(u^{\ast,1}, u^{\ast,2}, \dots, u^{\ast,K}\right)$.
\State Define the feedback value $\kappa(y(t^i)) = u^{\ast,1}$ and use this control to compute the associated state $y_{\kappa}(t^{i+1})$.
\EndFor
\end{algorithmic}
\end{algorithm}

\begin{rmrk}\label{rem:Multistep_Definition}
We shall denote the control horizon of the MPC algorithm by $\taucon$, i.e., the time period over which the control is applied to the plant. So far, we assumed $\taucon = \tau$ for the sake of readability and since this choice corresponds to the "classical" (discrete) approach. However, it is possible to decouple the (discrete) timestep $\tau$ used for the time discretization from the control horizon $\taucon$ resulting in a so-called multi-step feedback law, see e.g. \cite{Gr09}. In particular, if we consider the first $n$ control moves in the MPC algorithm instead of only the first one, we have $\tau = \taucon / n$, $n \in \xN$. Note that $n=1$ corresponds to the classical approach. 
\end{rmrk}

\begin{rmrk} \label{rem:Final_Cost}
There exists a wide range of stability results depending on the MPC scheme, i.e., MPC with stabilizing terminal constrains and/or cost as well as without terminal constraints or costs. We refer the interested reader to the survey paper~\cite{MRRS00} or the monographs~\cite{RM09,GP11}. Although we included a terminal cost in the cost functional in~\eqref{eq:FE_cost_functional}, all of the following results are also valid for problems without terminal cost, i.e., $\sigma_2 = 0$ in~\eqref{eq:FE_cost_functional}.
\end{rmrk}


\begin{dfntn}[Optimal value function]
Let $y_0\in Y$ be a given initial state. We define the infinite optimal value function $J_\infty^\ast(\cdot ;\mu): Y \to \xR_0^+$ as
\begin{align}
J_\infty^\ast(y_{0};\mu) = \inf_{u\in U} \{J_\infty(y, u;\mu)\ | \textrm{ s.t. } \eqref{eq:FE_PDEConstraint}\}. 
\end{align}
Similarly, we define the finite horizon optimal value function $J^\ast(\cdot;\mu) : Y \to \xR_0^+$ as 
\begin{align}
J^\ast(y_{0};\mu) = \min_{u\in U}  \{J(y, u;\mu) \ | \textrm{ s.t. } \eqref{eq:FE_PDEConstraint}\}.
\end{align}
\end{dfntn} 

Given a feedback law $\kappa:Y\to U$, we denote the associated solution by $y_\kappa$, i.e., $y_\kappa$ satisfies 
\begin{align}\label{eq:ClosedLoopSystem}
m(y_\kappa^k,v) + \tau a(y_\kappa^k,v;\mu) =m(y_\kappa^{k-1},v) + \tau b(\kappa(y_\kappa^{k-1}),v), \quad \forall v\in Y ,\quad \forall k\in \xN, 
\quad y_\kappa^0= y_0.
\end{align} 
The corresponding cost functional value is denoted by $J^\kappa_\infty(y_0;\mu):= J_\infty(y_\kappa, \kappa(y_\kappa);\mu )$. 

A feedback law stabilizes the system if it steers the state towards the desired trajectory and holds it there. For definitions and further theory we refer the reader e.g. to \cite{GP11, RM09, Kh15}. In \cite[Chapter 4]{GP11} it is shown that the solution of the optimal control problem over the infinite horizon always stabilizes the system since the optimal value function $J^\ast_\infty$ serves as a Lyapunov function. To ensure the desired property in the receding horzion approach, we rely on results in relaxed dynamic programming \cite{LR06, Ra06} and thus follow the ideas presented e.g. in \cite{Gr09, GPSW10}. Here, we use the following main result, see~\cite{GP11} for a proof.

\begin{prpstn}\label{prop:RelaxedDPP}
Consider a feedback law $\kappa : Y \to U$, the corresponding solution $y_{\kappa}$ to \eqref{eq:ClosedLoopSystem} with $y_{\kappa}^0= y_0$, and a function $\widetilde{V}(\cdot;\mu): Y \to \xR^+_0$ that  satisfies the inequality 
\begin{equation}\label{eq:LyapunovIneq}
\widetilde{V} (y_0;\mu) \ge \widetilde{V}(y_{\kappa}^1;\mu) +\omega(\mu) \ l(y_0, \kappa(y_0);\mu)
\end{equation}
for some $\omega(\mu)\in (0,1]$ and all $y_0\in Y$. Then for all $y\in Y$ the estimate
\begin{equation}
\omega(\mu) J^\ast_\infty(y;\mu) \le \omega(\mu) J^{\kappa}_\infty(y;\mu) \le \widetilde{V}(y;\mu)
\end{equation}
holds. The value $\omega(\mu)$ is called the suboptimality degree of the feedback law $\kappa$. If, in addition, there exist $\rho_1, \ \rho_2, \ \rho_3\in \mathcal{K}_\infty =\{ \rho : \xR^+_0 \to \xR^+_0 \ |  \ \rho \textrm{ is continuous, strictly increasing, unbounded and} \ \rho(0)=0 \}$ satisfying
\begin{align}
\rho_1(\|y_0-y^k_d(\mu)\|_Y)\ge \widetilde{V}(y_0;\mu)\ge \rho_2(\|y_0-y^k_d(\mu)\|_Y)
\quad \textrm{and} \quad 
\min_{u\in U_{ad}} l(y_0,u;\mu) \ge \rho_3(\|y_0-y^k_d(\mu)\|_Y) \label{eq:AsympStabCond}
\end{align}
for all $y_0\in Y$, $k\in \xN$, then the closed-loop system \eqref{eq:ClosedLoopSystem} is globally asymptotically stable.
\end{prpstn}

\begin{rmrk}\label{rem:Stability_by_Detectability}
The assumptions posed in Proposition \ref{prop:RelaxedDPP} can be weakened, see e.g. \cite{GP11,RM09}. 
For example, we introduced the cost functional using the observation subdomain $D\subset \Omega$.  Hence, in general, $|\cdot|_D$ defines only a semi-norm on $L^2(\Omega)$ resulting in a semidefinite function $l(\cdot,\cdot;\mu)$. 
Since the condition $\min_{u\in U_{ad}} l(y_0,u;\mu) \ge \rho_3(\|y_0-y^k_d(\mu)\|_Y)$ is not satisfied for all $y_0\in Y$, an additional assumption is necessary in order to establish asymptotic stability, see \cite{RM09}. 
One appropriate choice for such an requirement is detectability, which essentially says that although not all modes of a system can be observed, all unobservable modes are stable. Another possibility is to check if the system is input/output-to-state-stable (IOSS), see \cite[Section 2.7]{RM09}.
\end{rmrk}

\begin{dfntn}\label{def:FE_MPC_kappa}
Let $u^\ast$ be an admissible minimizing control for \eqref{eq:FE_OCP} with initial value $y_{0}\in Y$. We define the MPC feedback law by 
\begin{align} \label{eq:uast}
\kappa(y_0) = u^{\ast,1}.
\end{align}
\end{dfntn}

To ensure stability of the MPC controller $\kappa$ we utilize the relaxed dynamic programming result from Proposition \ref{prop:RelaxedDPP} for the choice $\widetilde{V} = J^\ast$. Using the definition of $\kappa$, we obtain that satisfying \eqref{eq:LyapunovIneq} is equivalent to ensuring the inequality 
\begin{align} \label{eq:FE_RelaxedDPP_Lyapunov}
J^\ast(y_0;\mu) \ge J^\ast(y_{u^\ast}^1;\mu) +\omega(\mu) \ l(y_0, u^{\ast, 1};\mu)
\end{align} 
holds for all $y_0 \in Y$ and some $\omega(\mu)\in (0,1]$ (preferably close to one). 
We note that the cost functional and the stability condition \eqref{eq:FE_RelaxedDPP_Lyapunov} depend on the prediction horizon $K$. 
It is well known that -- if the prediction horizon is chosen sufficiently large --  stability of the MPC controller can be guaranteed \cite{JH05}. 
Unfortunately, the computational complexity grows rapidly with the length of the optimization horizon. We are thus interested in the minimal horizon that guarantees stability of the controller. 

\begin{rmrk} 
A common assumption in the literature is the following exponential controllability condition: There exist real constants $C\ge 1$ and $o\in (0,1)$ such that for each $y_0\in Y$ there exists a control function $u_{y_0} \in U_{ad}$ satisfying 
\begin{align}
	l(y_{u_{y_0}}^k, u_{y_0}^{k+1};\mu) \le Co^k \min_{u\in U_{ad}^0} l(y_0, u;\mu )
\end{align}
for all $k\in \xN_0$. It was shown in \cite{GPSW10} that~\eqref{eq:LyapunovIneq} holds under this premise if the prediction horizon $K$ is chosen sufficiently large. In particular, one can show that $\omega_K(\mu)\to 1$ if $K\to \infty$. A similar result was presented in~\cite{JH05} for continuous time problems. 
\end{rmrk}

\begin{rmrk}\label{rem:Stability_ExistingWork}
There exist a wide range of stability and suboptimality analysis for MPC schemes based on Proposition \ref{prop:RelaxedDPP}. The goal is to estimate the suboptimality degree based on the inequality \eqref{eq:LyapunovIneq} and, since $\omega(\mu)$ depends on the prediction horizon, in turn to estimate the minimal stabilizing horizon. In~\cite{Gr09,GPSW10}, the controllability assumption mentioned in the last remark is used to express the stability and suboptimality condition in Proposition~\ref{prop:RelaxedDPP} as an optimization problem to bound $\omega(\mu)$. This allows to analyze and explain the qualitative behaviour of the prediction horizon $K$ on the stability and suboptimality of the MPC scheme~\cite{AG12}. Nevertheless, the bound on $\omega(\mu)$ -- and, as a consequence, the bound on $K$ -- tends to be conservative since it presents a global estimate~\cite{Gr09, AGW10}. To actually determine the minimal stabilizing horizon, it is also possible to use an adaptive approach as was suggested in \cite{GP08, GP09} and in~\cite{PW11} for the multi-step approach. Here, the calculation of a global minimal stabilizing horizon is omitted and, hence, less conservatism is expected. 



\end{rmrk}

\begin{rmrk}
In \cite{RA12} the stability of the MPC algorithm for continuous time problems was studied. A comparison of the continuous time and the discrete time suboptimality estimates revealed that the continuous time case is less conservative. To overcome this gap, \cite{WRGA14} employed a multi-step approach and showed that the discrete estimate converges to the continuous result when the sampling period $\tau$ tends to zero. We note that the stability result of Proposition \ref{prop:RelaxedDPP} needs to be adapted for the analysis of the multi-step approach, i.e., inequality \eqref{eq:LyapunovIneq} is replaced by 
\begin{equation}\label{eq:Multistep_Lyapunov}
J^\ast(y_0;\mu) \ge J^\ast(y_{u^\ast}^n;\mu) +\omega(\mu) \ \sum_{k=0}^{n-1} l(y_{u^\ast}^k, u^{\ast, k+1};\mu),
\end{equation}
where we directly considered $\tilde{V}=J^\ast$. The $n$-step MPC feedback law $\kappa : Y \times \{0,\dots,n-1\} \to U $ is then defined through $\kappa(y_0, i) = u^{\ast,i+1},$ $i=0,\dots, n-1,$ for an initial value $y_0\in Y$. If $y_0$ is the current state in the MPC algorithm then $y_{u^\ast}^n$ refers to the state reached after applying the control input, i.e., the starting point in the next loop. We refer to~\cite{Gr09} for details.
\end{rmrk}

In each step of the MPC algorithm, the optimality system~\eqref{eq:FE_OCPSystem} consisting of $2K\cN + Km$ coupled equations and inequalities has to be solved. The associated computational cost for large-scale problems, arising e.g. through discretization of PDEs, is thus considerable and may even be prohibitive in actual practice. Our goal in the remainder of this paper is thus the following: We want to significantly speed up the computation of the feedback law while maintaining the stability of the closed-loop system. To this end, we first substitute 
the PDE constraint in~\eqref{eq:FE_OCP} with a reduced basis model and use the optimal "reduced" control $u_N^\ast$ as feedback in the MPC algorithm. We then derive {\it a posteriori} error estimation procedures which allow us to confirm the stability of the closed-loop system, i.e., the reduced control applied to the original full-order system. Finally, based on the stability condition, we propose an adaptive approach to determine the minimal stabilizing horizon for the reduced feedback controller.

\section{Reduced Basis Method}
\label{sec:RB_Approx}

\subsection{Approximation}

We assume we are given the integrated reduced basis spaces $Y_N = \spn \{\zeta_n, \ 1\le n\le N\}$, $1\le N \le N_{\max}$, where $\zeta_n$ are $(\cdot,\cdot)_Y$-orthogonal basis functions. 
The space $Y_N$ is constructed using the POD/Greedy procedure illustrated in Section \ref{sec:Greedy}. 
 Note that we incorporate both state and the adjoint solutions in $Y_N$ -- thus the term ``integrated'' -- in order to guarantee stability of the reduced basis approximation and to obtain a good approximation for the state and adjoint using a single space.

The reduced basis optimal control problem reads 
\begin{equation}\label{eq:RB_OCP}\tag{$P_N$}
\begin{split}
 & \min_{y_N, u_N} J(y_N,u_N;\mu) \quad \textrm{s.t.} \quad (y_N,u_N)\in Y_N^{K}\times U_{ad} \textrm{ solves } \\
& m(y_N^k,v) + \tau a(y_N^k,v;\mu) =m(y_N^{k-1},v) +  \tau b(u_N^k,v), \quad \forall v\in Y_N, \quad \forall k\in \K, 
\end{split}
\end{equation}
with initial condition $y_N^0 = y_{N,0}$.  Here, $y_{N,0}$ is the projection of $y_0$ onto $Y_N$ and the cost functional is given by \eqref{eq:FE_cost_functional}.

We solve the minimization problem employing the corresponding first-order optimality system: Given $\mu\in \cD$, the optimal solution (to \eqref{eq:RB_OCP}) $(y_N^\ast, p_N^\ast, u_N^\ast)\in Y_N^{K}\times Y_N^{K}\times U_{ad} $ satisfies
\begin{subequations}\label{eq:RB_OCPSystem}
 	\begin{alignat}{2}
		 m(y_N^{\ast,k}-y_N^{\ast,k-1},\phi) + \tau a(y_N^{\ast,k},\phi;\mu) &= \tau  b(u_N^{\ast,k},\phi)  \quad 					
		 		&& \forall \phi \in Y_N, \quad \forall k\in \K, \\
 		m(y_N^{\ast,0},\phi) & = m(y_0,\phi)\quad 	&& \forall \phi \in Y_N, \\
 		m(\varphi, p_N^{\ast,k}-p_N^{\ast,k+1}) + \tau a(\varphi, p_N^{\ast,k};\mu) &= \tau  \sigma_1(y_{d}^k(\mu) - y_N^{\ast,k},   \varphi)_D   \quad 
 				&& \forall \varphi \in Y_N, \quad \forall k\in \K, \\
 				m(\varphi, p_N^{\ast,K}) + \tau a(\varphi, p_N^{\ast,K};\mu) &= (\tau  \sigma_1 +\sigma_2)(y_{d}^K(\mu) - y_N^{\ast,K},   \varphi)_D   \quad 
 				&& \forall \varphi \in Y_N,\\
		 (\lambda (u_N^\ast - u_{d}) - \cB^\star p_N^\ast, \psi - u_N^\ast)_{U} &\ge 0 \quad 
		 		&& \forall \psi\in U_{ad}. \label{eq:RB_OCPOptimalityInequality}
  \end{alignat}
\end{subequations}
We note that the FOTD and FDTO also commute in the reduced basis setting. The idea to combine the MPC algorithm with the reduced optimal control model is straightforward: 
Instead of solving \eqref{eq:FE_OCP} in each MPC loop we compute the reduced basis optimal control $u_N^{\ast}$ from \eqref{eq:RB_OCP} and set $\kappa(y(t^i))=u_N^{\ast,1}$. However, it remains to discuss and show the stability of this approach.

\subsubsection{Computational Procedure - Approximation}
\label{sec:Complexity_OffOnProcedure}

To determine the complexity of the proposed method, we first derive the matrix-vector representation of the reduced optimal control problem. 
To this end, we express $y_N^k(\mu)$ and $p_N^k(\mu)$ using the reduced basis consisting of $N$ elements, i.e., 
\begin{align}
y_N^k(\mu) = \sum_{i=1}^{N} y_{N \ i}^k(\mu)\zeta_i \ \textrm{ and } \
p_N^k(\mu) = \sum_{i=1}^{N} p_{N \ i}^k(\mu)\zeta_i.
\end{align}
The corresponding coefficient vectors are denoted by $\underline{y}_N^k(\mu)= [y_{N \ 1}^k(\mu), \dots , y_{N \ N}^k(\mu) ]^T \in \xR^{N}$ and $\underline{p}_N^k(\mu)= [p_{N \ 1}^k(\mu), \dots ,$ $ p_{N \ N}^k(\mu) ]^T \in \xR^{N}$, respectively. 
Next, we choose $\phi = \zeta_i$, $1\le i\le N$, $\varphi = \zeta_i$, $1\le i\le N$ and $\psi = e_i$, the $i$-th unit vector, $1\le i\le m$, as test functions in \eqref{eq:RB_OCPSystem}, which leads to the matrix-vector representation of the reduced optimal control problem
\begin{subequations}\label{eq:RB_MVRepresentation}
 	\begin{alignat}{3}
\left(M_N + \tau A_N(\mu)\right) \underline{y}_N^k(\mu) & = M_N \underline{y}_N^{k-1}(\mu) + \tau B_N u_N^k , \ &&  &&  k\in \K,\\
M_N\underline{y}_N^0(\mu)&= Z^T M\underline{y}_{0}, &&  && \\
\left(M_N + \tau A_N^T(\mu)\right) \underline{p}_N^k(\mu) & = M_N \underline{p}_N^{k+1}(\mu) + \tau \sigma_1 \left( Y_{d,N}^k(\mu) -  D_N \underline{y}_N^k(\mu)\right), \ && && k\in \K\backslash \{K\},\\
\left(M_N + \tau A_N^T(\mu)\right) \underline{p}_N^K(\mu) & = \left(\tau \sigma_1 +\sigma_2\right) \left( Y_{d,N}^K(\mu) -D_N   \underline{y}_N^K(\mu) \right), && && \\
\left( \lambda \left( u_N^k - u_d^k \right) -  B_N^T \underline{p}_N^{k}(\mu) \right)^T  & \left(\psi^k - u_N^k\right) \ge 0 , &&\forall \psi^k \in U_{ad}^k,  \ &&k\in \K.
	\end{alignat}
\end{subequations}
Here, $Z = [\zeta_1, \dots ,\zeta_{N} ]\in \xR^{\cN \times N}$ is the change of basis matrix, $M\in \xR^{\cN \times \cN}$ denotes the finite element mass  matrix associated with the inner product $(\cdot, \cdot)_{L^2}$ and $\underline{y}_{0}\in \xR^{\cN} $ is the coefficient vector of $y_0\in Y$ in the finite element basis. 
The matrices $M_N\in \xR^{N\times N}$, $D_N\in \xR^{N\times N}$ and $B_N\in \xR^{N\times m}$ are given by $(M_N)_{ij} = m(\zeta_i, \zeta_j)$, $(D_N)_{ij} = (\zeta_i, \zeta_j)_D$, $1 \leq i,j \leq N$, and $(B_N)_{ij}= b_j(\zeta_i)$, $1 \leq i \leq N$, $1 \leq j \leq m$. 
Exploiting, the affine-parameter dependence \eqref{eq:AffineRepr_a} of the bilinear form $a(\cdot, \cdot;\mu)$ for the matrix $A_N\in \xR^{N\times N}$ leads to the representation $A_N(\mu)= \sum_{q=1}^{Q_a} \Theta_a^q(\mu) A^q_N$. 
The matrices $A^q_N\in \xR^{N \times N}$ are parameter-independent and defined by $(A^q_N)_{ij}  =a_q(\zeta_j,\zeta_i)$, $1\le q\le Q_a$,  $1 \leq i,j \leq N$. 
Likewise, using the affine representation of  $y_{d}$ \eqref{eq:Affine_DesiredState_infty} we obtain $Y_{d,N}^k(\mu) = \sum_{q=1}^{Q_{y_d}} \Theta^q_{y_d \, | \, [0,T]}(t^k;\mu)Y_{d,N}^q $. The parameter independent vectors $Y_{d,N}^q\in \xR^{N}$ are given by $(Y_{d,N}^q)_i  = (y_d^q,\zeta_i)_D$, $1\le q\le Q_{y_d}$, $1 \leq i \leq N$. 
In addition, we save the matrix $Y_{d,d}\in \xR^{Q_{y_d}\times Q_{y_d}}$, defined by $(Y_{d,d})_{p,q} = (y_d^p, y_d^q)_D$, $1\le p,q \le Q_{y_d}$ to ensure an efficient evaluation of the cost functional. 

We now summarize the most important steps of the offline-online decomposition: In the offline phase, we compute the parameter independent components, i.e., the matrices $M_N, D_N, B_N$, $A_N^q$, $1\le q\le Q_a$, $Y_{d,N}^q$, $1\le q\le Q_{y_d}$ and $Y_{d,d}$ .  
Clearly, the computational cost depends on the truth finite element dimension $\cN$, but it is only performed once for the reduced basis $Y_N$. During the online phase, we first assemble the parameter-dependent quantities $A_N(\mu)$, $Y_{d, N}^k(\mu) ,\ k \in \K\cup \{0\}$, which can be done in $\cO(Q_aN^2)$ and $\cO(Q_{y_d}KN)$, respectively. 
The optimal control problem is then solved iteratively using the BFGS Quasi-Newton method \cite{NW06}. We state the effort for one iteration of the algorithm. 
By first computing an LU-factorization of the system matrix $M_N + \tau A_N(\mu)$ 
we can solve the reduced basis state and the adjoint equation in  $\cO(N^3+KN^2)$ operations. However, since the initial condition $\underline{y}_{0}$ is in general nonzero and unknown during the offline stage, we note that we must also compute the right hand side of (\ref{eq:RB_MVRepresentation}b) at cost $\cO(N \cN)$ and then solve for $\underline{y}_N^0(\mu)$ at cost $\cO(N^3)$, i.e., we cannot avoid a linear scaling with $\cN$ during the online stage due to the projection of $y_0\in Y$ onto the reduced basis space $Y_N$.

For an efficient evaluation of the cost functional, we rewrite 
\begin{align}
\begin{split}
J(y_N,u_N;\mu) =& \dfrac{\sigma_1}{2} \tau \sum_{k=1}^K 
		\left( \left( \underline{y}_N^k\right)^T D_N  \underline{y}_N^k - 2 \left(Y_{d,N}^k(\mu)\right)^T \underline{y}_N^k + Y_{d,d}^k(\mu)
				\right)\\
				&+ \dfrac{\sigma_2}{2}\left( \left( \underline{y}_N^K\right)^T D_N  \underline{y}_N^K - 2 \left(Y_{d,N}^K(\mu)\right)^T \underline{y}_N^K + Y_{d,d}^K(\mu)
				\right)
				+ \dfrac{\lambda}{2}\tau \sum_{k=1}^K \|u_N^k - u_d^k\|^2_{\xR^m},
\end{split}
\end{align}
with $ Y_{d,d}^k(\mu) = \sum_{p,q=1}^{Q_{y_d}} \Theta_{y_d}^p(t^k;\mu) \Theta_{y_d}^q(t^k;\mu) (Y_{d,d})_{p,q}$. Evaluating the cost functional using this representation then requires $\cO(KN^2 + KQ_{y_d}N + KQ_{y_d}^2 + Km)$ operations. 
Furthermore, we compute the gradient and the optimal step size of $J(y_N(u_N),u_N;\mu)$ in the search direction $u_N$. The former can be done in  $\cO(K mN)$, the latter can be determined in $\cO(N^3 + KN^2)$, see Appendix D in \cite{KaPhd}. 
Finally, the update of the inverse Hessian approximation is computed in $\cO((Km)^2)$ operations, which leads to an overall cost of $\cO(N^3 + KN^2 + KQ_{y_d} + KQ_{y_d}^2 + Km + KmN + (Km)^2)$ for one BFGS iteration. 
 To ensure compliance with the control constraints, we use the primal-dual active set method (PDAS) \cite{HIK03}, implemented as an outer loop around the BFGS iteration. Experience shows that the PDAS method requires only a few iterations to find a solution. In our numerical experiments the algorithm stopped after no more than six iterations.  In summary, all listed calculations -- except for computing the projection of $y_0$ -- can be performed independently of $\cN$. We therefore expect a significant speed-up compared to the solution of \eqref{eq:FE_OCP}, since in general $N\ll \cN$. We present numerical evidence in Section~\ref{sec:NumExp}.

\subsection{A Posteriori Error Estimation}\label{sec:APosterioriErrorBounds}

We now develop reliable and computationally efficient {\it a posteriori} error bounds for the error in the optimal control and the associated cost functional. As shown in the next section, these bounds will allow us to efficiently verify the stability of the proposed reduced basis MPC algorithm. To begin, we require the following definitions.

\begin{dfntn}\label{defi:Residuals}
The residuals for the state and adjoint equation and the optimality condition are defined by
\begin{subequations}
 	\begin{alignat}{2}
r_y^{k} (\phi;\mu) &= b(u_N^{\ast,k}, \phi) - a(y_N^{\ast,k}, \phi;\mu) -\dfrac{1}{\tau} m(y_N^{\ast,k} - y_N^{\ast,k-1}, \phi) , &&\quad \forall \phi\in Y, \ \forall k\in \K, \\
r_p^{k}(\varphi;\mu) &= \sigma_1 (y_d^{k}(\mu) - y_N^{\ast,k}, \varphi)_D - a(\varphi,p_N^{\ast,k} ;\mu) -\dfrac{1}{\tau} m( \varphi,p_N^{\ast,k} - p_N^{\ast,k+1}) ,&& \quad \forall \varphi\in Y, \ \forall k\in \K \backslash \{K\} ,\\
r_p^{K} (\varphi;\mu) &= \left(\sigma_1 + \dfrac{\sigma_2}{\tau}\right) (y_d^K(\mu) - y_N^{\ast,K}, \varphi)_D - a(\varphi,p_N^{\ast,K} ;\mu) -\dfrac{1}{\tau} m( \varphi,p_N^{\ast,K}) ,&& \quad \forall \varphi\in Y, \\
r_{u,i}^{k} (\psi;\mu) & =   \left(\lambda ( u_{N,i}^{\ast,k}-u_{d,i}^k ) - b_i(p_N^{\ast,k})\right) \psi_i^k,&& \quad \forall \psi_i^k\in U^k_i, \ \forall k\in \K, \ i=1,\dots, m.
  \end{alignat}
\end{subequations}
\end{dfntn}

We shall assume that we are given a positive lower bound $\alpha_{\LB}: \cD \to \xR^+$ for the coercivity constant $\alpha(\mu)$ such that
\begin{align}\label{eq:Alpha_LB}
\alpha(\mu) \ge \alphaLB \ge \alpha_0 >0, \quad \forall \mu \in \cD,
\end{align} 
and define the constant  
\begin{align}\label{eq:ConstantC_D}
C_D \equiv \sup_{v \in Y\backslash\{0\}} \frac{\|v\|_{L^2(D)}}{\|v\|_Y}.
\end{align}
We also define, for all $k \in \K$,  the energy norms
\begin{align} \label{eq:st-energy-norm}
\interleave v^k \interleave_y \equiv \left( m(v^k,v^k) + \tau \sum_{k'=1}^k a\left(v^{k'}, v^{k'};\mu\right) \right)^{\frac{1}{2}} \quad \text{and} \quad \interleave v^k \interleave_p \equiv \left( m(v^k,v^k) + \tau \sum_{k'=k}^K a\left(v^{k'}, v^{k'};\mu\right) \right)^{\frac{1}{2}},
\end{align}
as well as the standard "spatio-temporal" energy error bounds
\begin{align} \label{eq:st-energy-norm-bound}
\Delta^{y,k}_N(\mu) \equiv \left( \dfrac{\tau}{\alphaLB} \sum_{k'=1}^k \| r_y^{k'} (\cdot;\mu) \|^2_{Y'}\right)^{\frac{1}{2}} \quad \text{and} \quad \Delta^{p,k}_N(\mu) \equiv \left( \dfrac{\tau}{\alphaLB} \sum_{k'=k}^K \| r_p^{k'} (\cdot;\mu) \|^2_{Y'}\right)^{\frac{1}{2}}.
\end{align}

\subsubsection{Optimal Control Error Bound}

We first present two auxiliary results for the state and adjoint optimality error. The results are based on~\cite{KG14} but taking the nonzero error in the initial condition into account. We present a short sketch of the proofs in the appendix.
\begin{lmm}\label{lem:Hilfslem_OptStateErr}
The state optimality error $e_y^{K}  = y^{\ast, K}(u^\ast) -y_N^{\ast, K}(u_N^\ast) $ satisfies
\begin{align}
\interleave e_y^{K} \interleave_y \le  \left( \dfrac{2\tau}{\alphaLB} \sum_{k=1}^K \| r_y^{k} (\cdot;\mu) \|^2_{Y'}
									+ \dfrac{2}{\alphaLB}\left( \sum_{i=1}^m \|b_i\|_{Y'}^2 \right)  \tau \sum_{k=1}^K \| e_u^k \|_{R^m}^2  + m(e_y^0,e_y^0)
\right)^{\frac{1}{2}}.
\end{align}
Furthermore, we have 
\begin{align}\label{eq:StateOptErrBound_ProofDeltaALT}
\tau  \sum_{k=1}^K \| e_y^k \|_Y^2  \le \dfrac{2 \tau}{\alphaLB^2}   \sum_{k=1}^K \|r_y^k(\cdot;\mu) \|_{Y'}^2 +  \dfrac{2}{\alphaLB^2} \left( \sum_{i=1}^m \|b_i\|_{Y'}^2 \right) \tau \sum_{k=1}^K \| e_u^k \|_{R^m}^2 + \dfrac{\norm{e^0_y}_{L^2(\Omega)}^2}{\alphaLB}  .
\end{align}
\end{lmm}

\begin{lmm}\label{lem:Hilfslem_OptAdjErr}
The adjoint optimality error $e_p^{1}  = p^{\ast, 1}(y^{\ast}(u^\ast)) -p_N^{\ast, 1}(y_N^{\ast}(u_N^\ast) )$ satisfies
\begin{align}\label{eq:AdjOptErrBound_ProofDeltaALT_Enorm}
\interleave e_p^{1} \interleave_p \le \left( \dfrac{2\tau}{\alphaLB} \sum_{k=1}^K \| r_p^{k} (\cdot;\mu) \|^2_{Y'}
									+ \frac{2C_D^2 \sigma_1^2 \tau}{\alphaLB}\sum_{k=1}^K (e_y^k, e_y^k)_D  +  \sigma_2^2 (e_y^K, e_y^K)_D 
\right)^{\frac{1}{2}}.
\end{align}
with $C_D$ defined in \eqref{eq:ConstantC_D}. Furthermore, we have
\begin{align}\label{eq:AdjOptErrBound_ProofDeltaALT}
\tau  \sum_{k=1}^K \|e_p^k \|_{Y}^2
		\le  \dfrac{2 \tau}{\alphaLB^2}   \sum_{k=1}^K \|r_p^k(\cdot;\mu) \|_{Y'}^2 +  \dfrac{2C_D^2\sigma_1^2 \tau}{\alphaLB^2} \sum_{k=1}^K (e_y^k, e_y^k)_D + \dfrac{\sigma_2^2}{\alphaLB} (e_y^K, e_y^K)_D. 
\end{align}
\end{lmm}

We note that the following bound for the error in the optimal control is valid both in the (control) constrained and unconstrained case. The proof follows ideas presented in \cite{KBGV18} and is also sketched in the appendix.

\begin{prpstn}\label{prop:OptControlErrBound}
Let $u^\ast$ and $\uN$ be the optimal solutions of the truth and reduced basis optimal control problems \eqref{eq:FE_OCP} and \eqref{eq:RB_OCP}, respectively. Given $\mu \in \cD$,  the error in the optimal control satisfies
\begin{align}
\|u^\ast - u_N^\ast \|_{U} \le \Delta^{u,\ast}_{N}(\mu),
\end{align}
where $\Delta^{u,\ast}_{N}(\mu):= c_1(\mu) + \sqrt{c_1(\mu)^2 +c_2(\mu)} $ with nonnegative coefficients
\begin{align}
c_1(\mu) &=  \dfrac{1}{\sqrt{2}\alphaLB\lambda} \left( \sum_{i=1}^m \|b_i\|_{Y'}^2 \right)^{\frac{1}{2}} R_p  ,\\
c_2(\mu) &= \dfrac{1}{\lambda}\left[ \left(
 \dfrac{2\sqrt{2}}{\alphaLB}R_y 
	+ \dfrac{1+\sqrt{2}}{\sqrt{\alphaLB}} R_0 \right)R_p    
	 + \left( \dfrac{C_D^2\sigma_1}{\alphaLB} + \dfrac{\sigma_2}{2}\right) R_0^2  
	+ \left( \dfrac{C_D^2\sigma_1}{\alphaLB^2} 
	+ \dfrac{\sigma_2}{2\alphaLB}\right) R_y^2  \label{eq:c_2}
\right]
\end{align}
and the terms $R_p$, $R_y$ and $R_0$ are given by
\begin{align}
R_p = \left( \tau \sum_{k=1}^K  \| r_p^{k} (\cdot;\mu) \|^2_{Y'}\right)^{\frac{1}{2}},
\quad R_y = \left( \tau \sum_{k=1}^K  \| r_y^{k} (\cdot;\mu) \|^2_{Y'}\right)^{\frac{1}{2}},
 \quad R_0  = \norm{y_0 - y_{N,0} }_{L^2(\Omega)} .
\end{align}
\end{prpstn}

\subsubsection{Cost Functional Error Bound}

Before we derive a bound for the optimal cost functional error, we require the following {\it a posteriori} error bounds for the optimal state and adjoint. The following two lemmata extend the results in~\cite{KG14} to nonzero initial conditions. We thus omit the details and refer to~\cite{KG14} for the proof.
\begin{lmm}\label{lem:OptStateErrBound0}
The state optimality error $e_y^{k}  = y^{\ast, k}(u^\ast) -y_N^{\ast, k}(u_N^\ast) $ satisfies
\begin{align}
\interleave e_y^{k} \interleave_y \le \Delta^{y,\ast, k}_{N}(\mu), \quad \forall \mu\in \cD, \ \forall k\in \K,
\end{align}
where  the error bound $\Delta^{y,\ast,k}_{N}(\mu)$ is defined as 
\begin{align}
\Delta^{y,\ast,k}_{N}(\mu) \equiv \left( \dfrac{2\tau}{\alphaLB} \sum_{k'=1}^k \| r_y^{k'} (\cdot;\mu) \|^2_{Y'}
									+ \dfrac{2}{\alphaLB}\left( \sum_{i=1}^m \|b_i\|_{Y'}^2 \right) \left( \Delta^{u,\ast}_{N}(\mu) \right)^2  + 					\|y_0 - y_{N,0} \|^2_{L^2(\Omega)} 
\right)^{\frac{1}{2}}.
\end{align}
\end{lmm}
\begin{lmm}\label{lem:OptAdjointErrBound}
The adjoint optimality error $e_p^{k}  = p^{\ast, k}(y^{\ast}(u^\ast)) -p_N^{\ast, k}(y_N^{\ast}(u_N^\ast) )$ satisfies
\begin{align}
\interleave e_p^{k} \interleave_p \le \Delta^{p,\ast,k}_{N}(\mu), \quad \forall \mu\in \cD, \ \forall k\in \K,
\end{align}
where  the error bound $ \Delta^{p,\ast, k}_{N}(\mu)$ is defined as 
\begin{align}\label{eq:Delta_AdjPredictError0}
\Delta^{p,\ast, k}_{N}(\mu) \equiv \left( \dfrac{2\tau}{\alphaLB} \sum_{k'=k}^K \| r_p^{k'} (\cdot;\mu) \|^2_{Y'}
									+\left(\dfrac{2C_D^4\sigma_1^2}{\alphaLB^2}  + \dfrac{\sigma_2^2}{2} \right)\left( \Delta^{y,\ast,K}_{N}(\mu)\right)^2 
\right)^{\frac{1}{2}},
\end{align}
with $C_D$ as defined in \eqref{eq:ConstantC_D}.
\end{lmm}

The error bound for the cost functional is based on the dual weighted residual approach~\cite{VW08}. At this point we have to distinguish between the control constrained and unconstrained case since we cannot bound the absolute value of the cost functional error in the constrained case.

\begin{prpstn}\label{prop:OptCostInitErrBound}
Let $J^\ast = J(y^\ast, \uast;\mu)$ and $J_N^\ast = J(y_N^\ast, \uN;\mu)$ be the optimal values of the cost functionals of the truth and reduced basis optimal control problems, respectively. In the control constrained case, the cost functional error then satisfies
\begin{multline}\label{eq:OptCostInitErrBound_C}
J^{\ast,c} - J^{\ast,c}_N 
\le \Delta^{J,\ast,c}_{N}(y_0;\mu) 
:=  \dfrac{1}{2} \left( \left(\|e_y^0\|_{L^2(\Omega)} + \Delta^{y,K}_{N}(\mu) \right) \Delta^{p, \ast,1}_{N}(\mu) 	
		+ \Delta^{p,1}_N(\mu) \Delta^{y,\ast,K}_{N}(\mu)
	+ \left( \tau \sum_{i=1}^m \sum_{k=1}^K \left|\tilde{r}_{u,i}^{k}\right|^2 			\right)^{\frac{1}{2}}  \Delta^{u,\ast}_{N}(\mu) \right),\\
		 \forall \mu \in \cD,
\end{multline}
where $\tilde{r}_{u,i}^{k} = \lambda (u_{N,i}^{\ast,k} - u_{d,i}^k) - b_i(p^{\ast,k}_N)$. 
If we consider the problem without control constraints, the bound reduces to
\begin{equation}\label{eq:OptCostInitErrBound_UC}
|J^{\ast,uc} - J^{\ast,uc}_N |
\le \Delta^{J,\ast,uc}_{N}(y_0;\mu) 
:= \dfrac{1}{2} \left( \left(\|e_y^0\|_{L^2(\Omega)} + \Delta^{y,K}_{N}(\mu) \right) \Delta^{p, \ast,1}_{N}(\mu) 	
		+ \Delta^{p,1}_N(\mu) \Delta^{y,\ast,K}_{N}(\mu)\right) ,
		 \forall \mu \in \cD.
\end{equation}

\end{prpstn}
\begin{proof} 
We use the result in \cite{VW08} to estimate the cost functional error by   
\begin{equation}\label{eq:Estimate_result_VW08}
J(y^\ast,u^\ast;\mu) - J(y_N^\ast,u_N^\ast;\mu)
 \le  
\dfrac{1}{2} \left( m(e_y^0, e_p^0) - \tau \sum_{k=1}^K r_y^{k}(e_p^k;\mu) - \tau \sum_{k=1}^K r_p^{k}(e_y^k;\mu) + \tau \sum_{i=1}^m \sum_{k=1}^K r_{u,i}^{k}(e_u^k;\mu)\right) .
\end{equation}
Next, we derive from the Cauchy-Schwarz inequality and $e^{0}_p = e^{1}_p$ that 
\begin{multline}
J(y^\ast,u^\ast;\mu) - J(y_N^\ast,u_N^\ast;\mu)
 \le  \dfrac{1}{2}   \|e_y^0\|_{L^2(\Omega)} \|e_p^0\|_{L^2(\Omega)} 
		+ \dfrac{1}{2}  \left( \tau \sum_{k=1}^K \|r_y^{k}(\cdot;\mu)\|^2_{Y'}   \right)^{\frac{1}{2}} \left( \tau \sum_{k=1}^K  \|e_p^k\|_Y^2  \right)^{\frac{1}{2}} \\
		+ \dfrac{1}{2}  \left( \tau \sum_{k=1}^K \|r_p^{k}(\cdot;\mu)\|^2_{Y'}   \right)^{\frac{1}{2}} \left( \tau \sum_{k=1}^K  \|e_y^k\|_Y^2  \right)^{\frac{1}{2}} 
		+ \dfrac{1}{2}  \left( \tau \sum_{i=1}^m \sum_{k=1}^K |r_i^{u,k}(\cdot;\mu)|^2   \right)^{\frac{1}{2}} \|e_u\|_U  \\
\le \dfrac{1}{2} \left( \left(\|e_y^0\|_{L^2(\Omega)} + \Delta^{y,K}_{N}(\mu) \right) \Delta^{p, \ast,1}_{N}(\mu) 	
		+ \Delta^{p,1}_N(\mu) \Delta^{y,\ast,K}_{N}(\mu)  
		+ \left( \tau \sum_{i=1}^m \sum_{k=1}^K \left|\tilde{r}_i^{u,k}\right|^2 			\right)^{\frac{1}{2}}  \Delta^{u,\ast}_{N}(\mu) \right) . 
\end{multline}
We obtain the bound in \eqref{eq:OptCostInitErrBound_UC} as the estimate in \eqref{eq:Estimate_result_VW08} simplifies to  
\begin{equation}
J(y^\ast,u^\ast;\mu) - J(y_N^\ast,u_N^\ast;\mu) =  
\dfrac{1}{2} \left( m(e_y^0, e_p^0) - \tau \sum_{k=1}^K r_y^{k}(e_p^k;\mu) - \tau \sum_{k=1}^K r_p^{k}(e_y^k;\mu) \right) 
\end{equation}
in the unconstrained case.
\end{proof}

\begin{rmrk}
If $y_0 \in Y_N$ the error in the initial condition is zero. The error bounds derived in the last section simplify accordingly, i.e., all of the terms containing the error $e_y^0= y_0 - y_{N,0}$ and $R_0 = \norm{e^0_y}_{L^2(\Omega)}$ vanish. Also, the last term in $c_2(\mu)$ in \eqref{eq:c_2} is multiplied by $1/2$. 
\end{rmrk}

\subsubsection{Computational Procedure - Error Estimation}\label{sec:Complexity_Bounds}

The offline-online decomposition of the error bounds -- mainly the dual norm of the residuals -- is fairly standard by now. We therefore only summarize the computational costs involved. We assume here that the projection of the initial condition onto the reduced basis space is given since it was already computed during the solution of the optimal control problem. 

For the evaluation of the error bounds, the following components have to be computed: the dual norms of the state and adjoint residual equations $\|r_y^k(\cdot;\mu)\|_{Y'}$ and $\|r_p^k(\cdot;\mu)\|_{Y'}$ for $k\in \K$, the dual norms $\|b_i\|_{Y'}$ for $1\le i\le m$, the constant $C_D$ and the lower bound $\alphaLB$ for the coercivity lower bound $\alpha(\mu)$. In addition, we need to compute  all terms related to the initial error $e_y^0 = y_0 - y_{N,0}$, i.e., the norms $\|e_y^0\|_{L^2(\Omega)}$ and $\|e_y^0\|_{L^2(D)}$. For the dual norms $\|r_y^k(\cdot;\mu)\|_{Y'}$ and $\|r_p^k(\cdot;\mu)\|_{Y'}$ we employ the standard offline-online decomposition~\cite{GP05}. In the online phase we need $\cO((Q_a N +m)^2)$ operations for the evaluation of $\|r_y^k(\cdot;\mu)\|_{Y'}$ and $\cO((Q_a N +Q_{y_d})^2)$ for $ \|r_p^k(\cdot;\mu)\|_{Y'}$, $k\in \K$. Likewise we calculate the dual norms $\|b_i\|_{Y'}$ during the offline phase utilizing the Riesz representation theorem. The constant $C_D$ can be computed offline solving a generalized eigenvalue problem, see Appendix E in \cite{KaPhd}. There are several possibilities to determine the coercivity lower bound $\alphaLB$. For instance, the min-$\Theta$ approach is a simple method in case of parametrically coercive problems. In general, the successive constraint method~\cite{HRSP07} can be applied.  Overall, for a given parameter $\mu\in \cD$ and corresponding solution $(y^\ast_N,p^\ast_N,u^\ast_N)$ to \eqref{eq:RB_OCP} the online computational cost to compute the error bounds $\Delta^{u,\ast}_{N}(\mu)$ and $\Delta^{J,\ast}_{N}(y_0;\mu)$ scales with $\cO(K(Q_aN +Q_{y_d})^2 +K(Q_a N +m)^2)$ and is thus independent of $\cN$.


\subsection{Greedy Algorithm}\label{sec:Greedy}

The reduced basis spaces $Y_N$ are constructed using the POD/Greedy approach for optimal control problems introduced in \cite{KG14}; the approach is summarized in Algorithm \ref{alg:Greedy}. Here, the operator $\textrm{POD}_{Y}(\{ v_k \ : \ k\in \K\})$ returns the largest POD mode for a given time history $v_k, \ k\in \K$ with respect to the $(\cdot, \cdot)_Y$-inner product. 
The state and adjoint projection errors are given by $\{ e^{y,k}_{\textrm{proj}, N}(\mu^\ast) = y^{\ast,k}(\mu) - y^{\ast,k}_{\textrm{proj}, N}(\mu) \ : \ k\in \K\}$ and $\{ e^{p,k}_{\textrm{proj}, N}(\mu^\ast) = p^{\ast,k}(\mu) - p^{\ast,k}_{\textrm{proj}, N}(\mu) \ : \ k\in \K\}$, where $v^{k}_{\textrm{proj}, N}(\mu)$ denotes the $Y$-orthogonal projection of $v_k(\mu)$ onto $Y_N$. 

Given a desired error tolerance $\epsilon>0$ and a (sufficiently fine) training sample $\Xi_{train}\subset \cD$, we initialize the reduced basis space $Y_N$ as follows: If the initial condition of the state equation is 
nonzero, we set $Y_1 \leftarrow \spn\{y_0\}$ and thus obtain $e^0_y = 0$, cf.~\cite{GP05}.  In the case of a zero initial condition, we initialize $Y_N$ with the largest POD mode of the optimal state solution at a parameter value  $\mu^1\in \Xi_{train}$ (usually $\mu_{\rm min}$ or $\mu_{\rm max}$) and then append the largest POD mode of the adjoint projection error at the same parameter value.  We then successively expand $Y_N$ by repeating the following two steps: first, we search the parameter space for a parameter value $\mu^\ast$ that maximizes the error metric $\Delta_N(\mu)$; and second, we enlarge the reduced basis with the largest POD modes of the state and the adjoint projection errors at $\mu^\ast$, see steps \ref{alg:PODState} and \ref{alg:PODAdjoint}. The algorithm finishes as soon as the error metric is below the desired error tolerance for all parameters in $\Xi_{train}$. As the error metric $\Delta_N(\mu)$ in the unconstrained case we propose to use either the relative cost functional error bound $\Delta^{J,\ast}_N(\mu)/ J^\ast_N(\mu)$ or the relative control error bound $\Delta^{u,\ast}_N(\mu)/ \|u^\ast_N(\mu)\|_{U}$.  In the constrained case we only use the relative control error bound since we cannot bound the absolute value of the cost functional error in this case.

We note that the POD/Greedy approach requires the solution of several truth optimal control problems and may thus be fairly expensive. However, reduced order models are usually beneficial in two contexts: the many-query context and the real-time context, see~\cite{RHP2008} for a discussion. The MPC framework clearly belongs to the latter, our primary interest is in a minimal marginal cost. Reduced order models may allow us to apply the MPC framework to problems which might not be feasible with classical discretization approaches. Finally, we note that -- if the problem is LTI -- it is possible to avoid the solution of the truth optimal control problem during the POD/Greedy procedure by training the reduced basis on an impulse response, see~\cite{KG14} for more details.

\begin{algorithm}[t]
\caption{Optimal POD/Greedy Sampling Procedure}\label{alg:Greedy}
\begin{algorithmic}[1]
\State Given: $\Xi_{train}\subset \cD$, an arbitrary $\mu^1\in \Xi_{train}$ and a tolerance $\epsilon>0$.
\If  {$y_0\neq0$}
	\State Set $N \leftarrow 1$, $\zeta_1=y_0 , \ Y_1 \leftarrow \spn\{\zeta_1\}$
\Else 
	\State Set $N \leftarrow 1$, $\zeta_1=\textrm{POD}_Y(\{ y^{\ast,k}(\mu^1)\ :\ k\in\K  \}) 	, \ Y_1 \leftarrow \spn\{\zeta_1\}$
	\State Set $N \leftarrow 2$, $\zeta_2=\textrm{POD}_Y(\{ e^{p,k}_{\textrm{proj}, N-1}			(\mu^1)\ :\ k\in\K  \}) , \ Y_2 \leftarrow Y_1 \oplus  \spn\{\zeta_2\}$
\EndIf
\State $\mu^\ast \leftarrow \argmax_{\mu\in \Xi_{train}} \Delta_N(\mu)$
\While{$ \Delta_N(\mu^\ast)> \epsilon$} 
\State Set $N \leftarrow N+1$, $\zeta_N=\textrm{POD}_Y(\{ e^{y,k}_{\textrm{proj}, N-1}(\mu^\ast)\ :\ k\in\K  \}) , \ Y_N \leftarrow Y_{N-1} \oplus  \spn\{\zeta_N\}$ \label{alg:PODState}
\State Set $N \leftarrow N+1$, $\zeta_N=\textrm{POD}_Y(\{ e^{p,k}_{\textrm{proj}, N-1}(\mu^\ast)\ :\ k\in\K  \}) , \ Y_N \leftarrow Y_{N-1} \oplus  \spn\{\zeta_N\}$ \label{alg:PODAdjoint} 
\State $\mu^\ast \leftarrow \argmax_{\mu\in \Xi_{train}} \Delta_N(\mu)$
\EndWhile
\State $N_{\max} \leftarrow N$
\end{algorithmic}
\end{algorithm}

\section{Stability of the RB-MPC controller}
\label{sec:Stability_RB_MPC}

We introduce our RB-MPC approach by combining the MPC stability results from Section \ref{sec:MPC} with the certified RB results for optimal control problems derived in the last section. We can thus not only efficiently evaluate the RB-MPC controller but also guarantee its stability. Whenever necessary we add a superscript $c$ in the control constrained or $uc$ in the unconstrained case.

\begin{dfntn}
Let $u_N^\ast$ be an admissible minimizing control for \eqref{eq:RB_OCP} with initial value $y_{N,0}$, the projection of $y_0$ onto $Y_N$. We define the RB-MPC feedback law by 
\begin{align} \label{eq:kappa_N}
\kappa_N(y_0) = u_N^{\ast,1}.
\end{align}
\end{dfntn}

\begin{dfntn}[Optimal value function $J_N^\ast$]
Let $y_0\in Y$ be a given initial state and $y_{N,0}$ be its projection onto $Y_N$. We define the finite horizon optimal value function $J_N^\ast(\cdot;\mu) : Y_N \to \xR_0^+$ by 
\begin{align}
J_N^\ast(y_{N,0};\mu) = \min_{u_{N}\in U}  \{J(y_N, u_N;\mu) \ | \textrm{ s.t. the PDE constraint in }  \eqref{eq:RB_OCP}\}.
\end{align}
\end{dfntn} 

We address stability of the controller by giving estimates related to its degree of suboptimality. We utilize the relaxed dynamic programming result from Proposition \ref{prop:RelaxedDPP} for the choice $\kappa = \kappa_N$ from~\eqref{eq:kappa_N} and $\widetilde{V} = J^\ast$. To guarantee stability, we thus need to ensure that the inequality 
\begin{align} \label{eq:RelaxedDPP_Lyapunov}
J^\ast(y_0;\mu) \ge J^\ast(y_{\kappa_N}^1;\mu) +\omega(\mu) \ l(y_0, \kappa_N(y_0);\mu) = J^\ast(y_{\uN}^1;\mu) +\omega(\mu) \  l(y_0, u_N^{\ast,1};\mu) 
\end{align} 
holds for all $y_0 \in Y$ and some $\omega(\mu)\in (0,1]$ (preferably close to one). 
We note that $y_{\uN}^1$ is the state reached by plugging $u_N^{\ast,1}$ into the full-order system and $J^\ast(y_{\uN}^1;\mu) $ is the cost function value corresponding to the solution of the full-order optimization problem  \eqref{eq:FE_OCP} with initial state $y_{\uN}^1$.  

Although we can compute the RB-MPC feedback law $\kappa_N(y_0)$ efficiently as shown in the last section, checking the validity of \eqref{eq:RelaxedDPP_Lyapunov} for a given $\kappa_N(y_0)$ requires the solution of two full-order optimal control problems and is thus not online-efficient. However, given the {\it a posteriori} error bounds derived in the last section, we can bound the optimal value functions of the full order problem with their respective reduced basis approximations. The latter can be computed online-efficient, thus allowing us to derive an estimate for $\omega(\mu)$ that can also be computed online-efficient and used to enforce the validity of \eqref{eq:RelaxedDPP_Lyapunov}. We recall that variants of algorithms computing $\omega(\mu)$ online were already considered in \cite[Chapter 7]{GP11}, but these approaches focused only on the high-fidelity model and did not consider model order reduction.


\subsection{Unconstrained Case} 

We first consider the problem without control constraints. We present the following result only for the single-step approach, i.e., where $\taucon = \tau$ as discussed in Remark~\ref{rem:Multistep_Definition}. The extension to the multi-step approach based on~\eqref{eq:Multistep_Lyapunov} is straightforward and thus omitted.

\begin{prpstn}\label{prop:Alpha_UC}
Given $y_0\in Y$, let $y_{N,0}$ be the projection of $y_0$ onto $Y_N$, $y_N^\ast$ and $u_N^\ast$ be the optimal state and optimal control solution to \eqref{eq:RB_OCP} with initial condition $y_{N,0}$ and optimization horizon $T = K\tau$. Furthermore, let $y^1_{u^\ast_N}$ be the state reached by plugging $u^{\ast,1}_N$ into the full-order system and $y^1_{u^\ast_N, N}$ be its projection onto $Y_N$. 
We define 
\begin{align} \label{eq:omegatilde_uc}
\tilde{\omega}_{N,K}(\mu) = \dfrac{J^{\ast,uc}_N(y_{N,0};\mu)- \Delta^{J,\ast,uc}_{N}(y_0;\mu)  -  J^{\ast,uc}_N(y^1_{u^\ast_N, N};\mu)  -  \Delta^{J,\ast,uc}_{N}(y^1_{u^\ast_N};\mu) }{\dfrac{\sigma_1}{2} \tau| y_0 - y_{d}^0(\mu)|^2_{D} +   \dfrac{\lambda}{2} \tau  \| u_N^{\ast,1} - u_{d}^{1}\|^2_{\xR^m}}.
\end{align}
If there exists a horizon $K$ such that $\tilde{\omega}_{N,K}(\mu)>0$ then \eqref{eq:RelaxedDPP_Lyapunov} holds for an $\omega(\mu)\in (0,1]$. 
\end{prpstn}
\begin{proof}
Since $J^{\ast,uc}(y_0;\mu)$ refers to the optimal cost functional value of \eqref{eq:FE_OCP} over the time horizon $[0,T]$ with initial condition $y_0\in Y$, we invoke the error bound \eqref{eq:OptCostInitErrBound_UC} to obtain 
\begin{align}\label{eq:Alpha_UC_Proof1}
J^{\ast,uc}(y_0;\mu) \ge J^{\ast,uc}_N(y_{N,0};\mu)- \Delta^{J,\ast,uc}_{N}(y_0;\mu) .
\end{align}
By definition of $\tilde{\omega}_{N,K}(\mu)$ we have
\begin{align}\label{eq:Alpha_UC_Proof2}
 J^{\ast,uc}_N(y_{N,0};\mu)- \Delta^{J,\ast,uc}_{N}(y_0;\mu)
 = J^{\ast,uc}_N(y^1_{u^\ast_N, N};\mu)  +  \Delta^{J,\ast,uc}_{N}(y^1_{u^\ast_N};\mu)
 +  \tilde{\omega}_{N,K}(\mu) \left(\dfrac{\sigma_1}{2} \tau| y_0 - y_{d}^0(\mu)|^2_{D} +   \dfrac{\lambda}{2} \tau  \| u_N^{\ast,1} - u_{d}^{1}\|^2_{\xR^m}\right).
\end{align} 
We apply \eqref{eq:OptCostInitErrBound_UC} a second time and employ the definition of $l(\cdot, \cdot;\mu)$ to get
 \begin{multline}\label{eq:Alpha_UC_Proof3}
 J^{\ast,uc}_N(y^1_{u^\ast_N, N};\mu)  +  \Delta^{J,\ast,uc}_{N}(y^1_{u^\ast_N};\mu)
 +  \tilde{\omega}_{N,K}(\mu) \left(\dfrac{\sigma_1}{2} \tau| y_0 - y_{d}^0(\mu)|^2_{D} +   \dfrac{\lambda}{2} \tau  \| u_N^{\ast,1} - u_{d}^{1}\|^2_{\xR^m}\right)
 \\ \ge J^{\ast,uc}(y_{\uN}^1;\mu) +\tilde{\omega}_{N,K}(\mu) \  l(y_0, u_N^{\ast,1};\mu).
\end{multline}
It thus follows from \eqref{eq:Alpha_UC_Proof1}, \eqref{eq:Alpha_UC_Proof2}, and \eqref{eq:Alpha_UC_Proof3} that 
\begin{equation}
\begin{aligned}
J^{\ast,uc}(y_0;\mu) 
 \ge& J^{\ast,uc}(y_{\uN}^1;\mu) +\tilde{\omega}_{N,K}(\mu) \  l(y_0, u_N^{\ast,1};\mu)\\
 \ge& J^{\ast,uc}(y_{\uN}^1;\mu) +\omega(\mu)  \  l(y_0, u_N^{\ast,1};\mu),
\end{aligned}
\end{equation}
with $\omega(\mu):=\min(\tilde{\omega}_{N,K}(\mu),1)$. 
If there exists a $K$ such that $\tilde{\omega}_{N,K}(\mu)>0$ then \eqref{eq:RelaxedDPP_Lyapunov} holds for $\omega(\mu)\in (0,1]$.
\end{proof}

\subsection{Constrained Case}

Unfortunately, in the presence of control constraints the cost functional error bound cannot serve as a lower bound. However, for a given input parameter $\mu$ the optimal value function of the constrained problem is an upper bound for the optimal value function of the unconstrained case since $U_{ad}\subset U$. Hence, we again invoke the corresponding error bound for the unconstrained case since
\begin{align}\label{eq:Alpha_C_add_est}
 J^{\ast, c}(y_0;\mu) \ge J^{\ast, uc}(y_0;\mu) \ge J_N^{\ast, uc}(y_{N,0};\mu) - \Delta^{J,\ast,uc}_{N}(y_0;\mu).
\end{align}
The proof of the following result directly follows from the proof of Proposition \ref{prop:Alpha_UC} and \eqref{eq:Alpha_C_add_est} and is thus omitted.
\begin{prpstn}\label{prop:Alpha_C}
Given $y_0\in Y$, let $y_{N,0}$ be the projection of $y_0$ onto $Y_N$, $y_N^\ast$ and $u_N^\ast$ be the optimal state and optimal control solution to \eqref{eq:RB_OCP} with initial condition $y_{N,0}$ and optimization horizon $T = K\tau$.  Furthermore, let $y^1_{u^\ast_N}$ be the state reached by plugging $u^{\ast,1}_N$ into the full-order system and $y^1_{u^\ast_N, N}$ be its projection onto $Y_N$. 
We define 
\begin{align} \label{eq:omegatilde_c}
\tilde{\omega}_{N,K}(\mu) = \dfrac{J^{\ast,uc}_N(y_{N,0};\mu)- \Delta^{J,\ast,uc}_{N}(y_0;\mu)  -  J^{\ast,c}_N(y^1_{u^\ast_N, N};\mu)  -  \Delta^{J,\ast,c}_{N}(y^1_{u^\ast_N};\mu) }{\dfrac{\sigma_1}{2} \tau| y_0 - y_{d}^0(\mu)|^2_{D} +   \dfrac{\lambda}{2} \tau  \| u_N^{\ast,1} - u_{d}^{1}\|^2_{\xR^m}}.
\end{align}
If there exists a horizon $K$ such that $\tilde{\omega}_{N,K}(\mu)>0$ then \eqref{eq:RelaxedDPP_Lyapunov} holds for an $\omega(\mu)\in (0,1]$. 
\end{prpstn}

\subsection{RB-MPC Algorithm}\label{sec:RB-MPC_Alg}

In the RB-MPC approach we do not aim to find a global minimal stabilizing horizon for all $y\in Y$ since this usually results in very conservative (worst case) estimates for the prediction horizon $K$ and associated higher computational cost. Instead, we propose an adaptive online computation of the reduced basis feedback controller. To this end, we first note that -- for a given parameter $\mu \in \cD$, state $y_0$, and prediction horizon $K$ -- computing the RB-MPC feedback law~\eqref{eq:kappa_N}  and evaluating the associated suboptimality degree $\tilde{\omega}_{N,K}(\mu)$ in \eqref{eq:omegatilde_uc} (resp.\ \eqref{eq:omegatilde_c}) only requires online computations: we need to solve the reduced basis optimal control problem and subsequently evaluate the {\it a posteriori} error bounds and the term $\dfrac{\sigma_1}{2} \tau| y_0 - y_{d}^0(\mu)|^2_{D} +   \dfrac{\lambda}{2} \tau  \| u_N^{\ast,1} - u_{d}^{1}\|^2_{\xR^m}$. The computational cost is summarized in Sections~\ref{sec:Complexity_OffOnProcedure} and~\ref{sec:Complexity_Bounds}. Although the online cost still depends linearly on $\cN$ due to the projection of $y_0$ and evaluation of $|y_0 - y_{d}^0(\mu)|^2_{D}$, we expect and observe -- see Section \ref{sec:NumExp} -- considerable computational savings compared to a full order MPC approach, i.e., evaluation of $\omega(\mu)$ based on~\eqref{eq:FE_RelaxedDPP_Lyapunov}.

To ensure stability of the feedback controller $\kappa_N$, it is sufficient to guarantee that~\eqref{eq:RelaxedDPP_Lyapunov} holds for all states along the closed-loop trajectory $y_{\kappa_N}$. The following approach attempts to find the minimum stabilizing horizon for a given reduced basis size $N$ by adaptively changing the prediction horizon, see Algorithm~\ref{alg:Adaptive_RB_MPC} for a summary: Given the state measurement $y(t^k)$ at a sampling instance $t^k$ we compute $\tilde{\omega}_{N,K}(\mu)$ as indicated in Proposition~\ref{prop:Alpha_UC} for the unconstrained case (resp.\ Proposition~\ref{prop:Alpha_C} for the constrained case) for increasing horizons $K=1,2,\dots$. For each $K$ (starting with the smallest), this involves solving the reduced optimal control problem~\eqref{eq:RB_OCP} twice. First, we optimize over $[t^k, t^k+K\tau]$ with initial condition $\textrm{proj}_{Y_N}(y(t^k))$, the projection of $y(t^k)$ onto $Y_N$, to obtain the reduced basis control candidate $u_N^{\ast,1}$. Afterwards, we solve the optimization problem over $[t^{k+1}, t^{k+1}+K\tau]$ with initial condition $\textrm{proj}_{Y_N}(y_{\uN}^1)$ and evaluate the suboptimality degree $\tilde{\omega}_{N,K}(\mu)$ from \eqref{eq:omegatilde_uc} (resp.\ \eqref{eq:omegatilde_c}). Once $\tilde{\omega}_{N,K}(\mu)$ becomes positive for a horizon $K$, we are guaranteed that the control candidate $u_N^{\ast,1}$ stabilizes the original high-fidelity system. We thus define $\kappa_N(y(t^k))= u_N^{\ast,1}$, update the state estimate $y_{\kappa_N}(t^{k+1})$ and proceed to the next timestep. 

We first note that if Algorithm~\ref{alg:Adaptive_RB_MPC} finishes without reaching $K_{\max}$, the computed reduced basis feedback controller guarantees stability of the original system. However, whether the suboptimality degree $\tilde{\omega}_{N,K}(\mu)$ computed in the while-loop in Algorithm~\ref{alg:Adaptive_RB_MPC} becomes positive before reaching $K_{\max}$ strongly depends on the choice of $K_{\max}$ and the accuracy of the RB approximation and thus $N$. To this end, we first comment on the unconstrained case and distinguish between two scenarios: ({\it i\/}) the suboptimality degree $\omega(\mu)$ based on~\eqref{eq:FE_RelaxedDPP_Lyapunov} becomes positive for $K = K_{\max}$ sufficiently large, i.e., the full order MPC approach results in a stabilizing control for $K_{\max}$; and ({\it ii\/}) $\omega(\mu)$ remains negative for all $K$. Concerning the former scenario, we know that $y_N^k \to y^k$ and $u_N^k \to u^k$ for $1 \leq k \leq K$ and all $\mu \in \cD$  as $N \to {\cal N}$; the reduced basis control will thus also stabilize the original dynamics for $K = K_{\max}$ as $N \to {\cal N}$. In order to adjust the accuracy of the RB approximation we may modify Algorithm~\ref{alg:Adaptive_RB_MPC} as follows: after line 10 we check if $\tilde{\omega}_{N,K}(\mu) \le 0$ and $K = K_{\mathrm{max}}$ and $N \leq {\cal N}$ is true, and then set $N = N+2$ and $K = 0$. Once the while-loop finishes, we additionally reset $N$ to the initial basis size $N_{\min}$ in line 4. The modification would gradually improve the accuracy of the RB approximation and theoretically guarantee stability of the RB-MPC approach for scenario ({\it i\/}). However, we note that considering the limit $N \to {\cal N}$ in practice is not only inefficient but also unnecessary from a computational point of view, since the {\it a posteriori} error bound is sensitive to round-off errors and converges to approximately square root of machine precision usually for $N \ll {\cal N}$~\cite{CEL2014}. Unfortunately, we cannot determine a maximum $N^*$ which guarantees sufficient accuracy of the RB approximation due to a lack of an {\it a priori} error bound for the optimal cost and the unknown minimal decrease in the cost functional in~\eqref{eq:omegatilde_uc}. Considering the second scenario, let us know assume that $\omega(\mu)$ in (31) remains negative for all $K$, i.e., the full-order MPC approach does not result in a stabilizing control even for $K \to \infty$. Since the RB approximation is build on the full-order system, we of course cannot expect the RB-MPC approach to deliver a stabilizing control (even as $N \to {\cal N}$). In this case the problem is inherent to the full-order approach, e.g., the high-fidelity discretization, and not a failure of the RB approach. To summarize, for scenario ({\it i\/}) we can theoretically guarantee stability of the RB-MPC approach using the modified Algorithm~\ref{alg:Adaptive_RB_MPC}, but remaining issus are round-off errors and efficiency as $N \to {\cal N}$. For scenario ({\it ii\/}) we cannot guarantee stability of the RB-MPC approach since even the standard (full-order) MPC approach would fail.

In fact, numerical/round-off errors are an issue not only for the RB-MPC approach, but also for the full-order MPC approach. Even if the full-order system is -- in theory -- stabilizable using MPC, the suboptimality degree $\omega(\mu)$ computed from~\eqref{eq:FE_RelaxedDPP_Lyapunov} may never become positive due to numerical errors. Especially in a neighborhood of the desired trajectory this will always be true since numerical errors start to dominate at one point~\cite{GR08,GP09}. We also observe this behaviour in our second numerical example in Section~\ref{sec:NumExp_WP} once the actual state is close to the desired state. This issue has been addressed in the literature under the notion ``practical optimiality'', see e.g.\cite{GR08,GP09}. The idea here is drive the system towards a small $\varepsilon$-neighborhood of the desired trajectory in order to account for numerical errors. In the RB context it may be possible to link the size of the $\varepsilon$-neighborhood to the accuracy of the RB approximation through the {\it a posteriori} error bounds. We may thus be able to account for round-off in the offline-online decomposition and in turn guarantee convergence to the $\varepsilon$-neighborhood of the desired trajectory. 

It remains to comment on the choice of $K_{\max}$. In~\cite{AGW10,AG12,AG13}, the authors derived global estimates for the minimal stabilizing horizon $K_{\infty}$ specifically for MPC of partial differential equations. The dependence of $K_{\infty}$  on system and regularization parameters was shown and confirmed in numerical studies. Since the global bounds for $K_{\infty}$ are usually quite pessimistic, a valid choice for $K_{\max}$ would be to set $K_{\max} = K_{\infty}$ if such an estimate exists for the problem at hand. In fact, we used the results in~\cite{AG13} to choose a suitable $K_{\max}$ for our first numerical example in Section \ref{sec:1DHeat}. Otherwise, we can only resort to (offline) parameter studies for various initial conditions and desired trajectories in order to find a suitable and conservative $K_{\max}$.

Finally, we note that our discussion above only holds for the unconstrained case. In the control constrained case we cannot guarantee that the suboptimality degree $\tilde{\omega}_{N,K}(\mu)$ eventually becomes positive even if we gradually improve the RB approximation quality, i.e., by modifying Algorithm~\ref{alg:Adaptive_RB_MPC}. The reason is that we need to use the unconstrained solution and cost in order to obtain a lower bound for the cost of the constrained problem. If the constrained cost always exceeds the unconstrained one, the suboptimality degree is negative even if the reduced model is perfect.

\begin{algorithm}[t]
\caption{Adaptive Reduced Basis Model Predictive Control}\label{alg:Adaptive_RB_MPC}
\begin{algorithmic}[1]
\State Given: time step $\tau >0$; 
maximal prediction horizon $K_{\mathrm{max}}\in \xN$, parameter $\mu\in \cD$, reduced basis space $Y_N$, initial condition $ y_{\kappa_N}(t^0)= y_0$. 
\State Set $\tilde{y}_d^k(\mu) = y_d^k(\mu)$, $\tilde{u}_d^k = u_d^k$, $k = 1,2,\ldots$.
\For{$i=0,1,2,\dots$} 
	\State $K \leftarrow 0$; $\tilde{\omega}_{N,K}(\mu) \leftarrow -1$
	\While{$\tilde{\omega}_{N,K}(\mu) \le 0$ \textbf{ and } $K\le K_{\mathrm{max}}$}
	\State $K \leftarrow K+1$
		\State Set $t^i=i\tau$, $y_0 = y_{\kappa_N}(t^i)$, and $y_d^k(\mu) = \tilde{y}_d^{k+i}(\mu)$, $u_d^k = \tilde{u}_d^{k+i}$, $k = 1,\ldots, K$.

		\State Solve the optimal control problem \eqref{eq:RB_OCP} 
 			with 
 			$y_N^0 = \textrm{proj}_{Y_N}(y_0)$ and obtain the 
 			 optimal solution $(y_N^\ast, u_N^\ast)$.
 		\State Solve the optimal control problem  \eqref{eq:RB_OCP}  
 			with 
 			$y_N^0 = \textrm{proj}_{Y_N}(y^1_{u_N^\ast})$ and obtain the 
 			 optimal solution $(\tilde{y}_N^\ast, \tilde{u}_N^\ast)$.
 		\State  Compute $\tilde{\omega}_{N,K}(\mu)$ according to Proposition \ref{prop:Alpha_UC} or Proposition \ref{prop:Alpha_C}.
	\EndWhile
	\State Define the feedback value $\kappa_N(y(t^i)) = u_N^{\ast,1}$ and compute the associated state $y_{\kappa_N}(t^{i+1})$.
\EndFor
\end{algorithmic}
\end{algorithm}


\begin{rmrk}\label{rem:omega_min}
In each loop of the algorithm, we compute the smallest $K$ with $\tilde{\omega}_{N,K}>0$. Although this approach guarantees the stability of the RB-MPC controller, the controller may not meet desired performance criteria due to the low suboptimality degree. As a result, the actual state of the closed-loop system can converge rather slowly towards the desired state. In case we want to enhance the suboptimality of the controller and thus achieve faster stabilization, we propose to incorporate a threshold $\omega_{\min}>0$ for the suboptimality degree. The prediction horizon in the RB-MPC algorithm is then only accepted as soon as $\tilde{\omega}_{N,K}(\mu)>\omega_{\min}>0$ holds. 
\end{rmrk}

\begin{rmrk}\label{rem:K_adaptive}
The proposed approach starts from the smallest possible prediction horizon in each loop, which is obviously not very efficient. In fact, we expect the minimal prediction horizon to vary rather slowly along trajectories. In practice, we would thus suggest to start each loop with an initial guess for $K$ depending on previous horizons. Then, starting from this initialization, the stability condition can be verified and $K$ may be accepted, decreased or increased accordingly. If the goal is to determine the minimal stabilizing horizon, we may apply a strategy to reduce $K$ if the initial guess directly results in a stabilizing controller. However, investigating different adaptive strategies is beyond the scope of the paper and we thus revert to Algorithm~\ref{alg:Adaptive_RB_MPC} for the numerical experiments in Section~\ref{sec:NumExp}.


Finally, we would like to note that the suboptimality degree may not increase monotonically in $K$, also see Section \ref{sec:Num_Exp_Alpha}. This should be taken into account in the adaptive choice of $K$.
\end{rmrk}

\section{Application to Weakly Coercive Problems}
\label{sec:Weakly_Coercive}

The proposed analysis is not restricted to systems with a spatially coercive bilinear form $a(\cdot, \cdot, \mu)$. We can extend our RB-MPC approach and treat optimal control systems where the governing PDE is only weakly coercive, thus making the proposed approach applicable to a wider range of problems. To this end, we introduce a temporal transformation to again obtain a spatially coercive but now time-variant optimal control problem. We first state the problem and then show how the reduced basis approximation and especially the {\it a posteriori} error bounds need to be adapted for the class of weakly coercive problems. 

\subsection{Problem Statement}

We now assume that the bilinear form $a(\cdot,\cdot;\mu)$ in the parametrized PDE~\eqref{eq:PDEConstraint} --- the PDE constraint of~\eqref{eq:ExactOCP} --- is only weakly coercive. More precisely, instead of satisfying~\eqref{eq:Exact_ellipticity} we assume that the bilinear form $a(\cdot,\cdot;\mu)$ only satisfies the weaker G\r{a}rding inequality, i.e., there exist constants $\delta_e(\mu)\ge0$ and $\alpha_e(\mu)\ge \alpha_0>0$ such that 
\begin{align}\label{eq:Garding_inequality}
a(v,v;\mu) + \delta_e(\mu) \|v\|^2_{L^2(\Omega)} \ge \alpha_e(\mu) \|v\|^2_{Y_e}, \quad \forall v\in Y_e, \ \forall \mu \in \cD.
\end{align}
If we introduce the standard change of variables $\hat{y}_e(x,t):= e^{-\delta_e t} y_e(x,t)$, see e.g.~\cite{QV08}, it is well-known that the new unknown $\hat{y}_e(x,t)$ satisfies the variational problem
\begin{equation}
\dfrac{\diff}{\diff t} m(\hat{y}_e(t),\phi) + a_\delta(\hat{y}_e(t),\phi;\mu) = e^{-\delta_e t} b(u_e(t),\phi),  \quad \forall \phi \in Y_e ,\quad \textrm{f.a.a. } t\in (0,\infty),
\end{equation}
where the bilinear form $a_{\delta}(\cdot,\cdot;\mu) : Y_e \times Y_e \to \xR$ is defined by $a_\delta(w,v;\mu):= a(w,v;\mu) + \delta_e\left(w,v\right)_{L^2(\Omega)}$. It  thus follows that $a_\delta(\cdot,\cdot;\mu)$ is coercive, i.e., in particular it holds that
\begin{align}
a_\delta(v,v;\mu) \ge \alpha_e(\mu) \|v\|^2_{Y_e}, \quad \forall v\in Y_e, \ \forall \mu \in \cD,
\end{align}
with $\alpha_e(\mu)$ as in \eqref{eq:Garding_inequality}. 

In the optimal control context, we additionally introduce a change of variables for the adjoint and the control given by $\hat{p}_e(x,t):= e^{\delta_e t} p_e(x,t)$ and $\hat{u}_e(t):= e^{-\delta_e t} u_e(t)$, respectively. Instead of solving the optimal control problem~\eqref{eq:ExactOCP} where the bilinear form $a(\cdot,\cdot;\mu)$ is only weakly coercive, we can solve the transformed problem with the coercive bilinear form $a_\delta(\cdot,\cdot;\mu)$. More precisely, for a given parameter $\mu\in \cD$, we consider the parametrized optimal control problem~\footnote{Often we omit the dependence of $\mu$ to simplify the notation, for instance $\delta_e = \delta_e(\mu)$. }  
 \begin{equation}\label{eq:TransformedOCP}\tag{$\hat{P}_e$}
 \begin{split}
 & \min_{\hat{y}_e, \hat{u}_e} \hat{J}_e(\hat{y}_e, \hat{u}_e;\mu) \quad \textrm{s.t.} \quad (\hat{y}_e,\hat{u}_e)\in W(0,T)\times U_{e,ad}(T) \textrm{ solves } \\
& \dfrac{\diff}{\diff t} m(\hat{y}_e(t),\hat{v}) + a_{\delta}(\hat{y}_e(t),\hat{v};\mu) = b(\hat{u}_e(t),\hat{v}), \quad \forall \hat{v}\in Y_e, \quad \textrm{f.a.a. } t\in (0,T], 
\end{split}
 \end{equation}
with initial condition $\hat{y}_e(0) = y_{0,e}\in L^2(\Omega)$. 
The quadratic cost functional is $\hat{J}_e(\cdot, \cdot;\mu) :  W(0,T)\times U_{e}(T) \to \xR$,  defined by 
$$ \hat{J}_e(\hat{y}_e, \hat{u}_e;\mu) = \dfrac{1}{2} \int_0^T \hat{\sigma}_1(t) | \hat{y}_e(t) - \hat{y}_{d,e}(t;\mu)|^2_{D} \diff t 
	+ \dfrac{\hat{\sigma}_2(T)}{2} | \hat{y}_e(T) - \hat{y}_{d,e}(T;\mu)|^2_{D} 
	+  \dfrac{1}{2}\int_0^T\hat{\lambda}(t) \| \hat{u}_e(t) - \hat{u}_{d,e}(t)\|^2_{\xR^m} 
$$
with time dependent optimization parameters $\hat{\sigma}_1(t):= e^{2\delta_e t} \sigma_1$, $\hat{\sigma}_2(T):= e^{2\delta_e T} \sigma_2$ and $\hat{\lambda}(t):= e^{2\delta_e t} \lambda$. The desired state and control are defined through $\hat{y}_{d,e}(t;\mu):= e^{-\delta_e t} y_{d,e}(t;\mu),$ $ \hat{y}_{d,e}(\cdot;\mu) \in L^2(0,T;L^2(\Omega))$ and $\hat{u}_{d,e}(t):= e^{-\delta_e t} u_{d,e}(t),$ $\hat{u}_{d,e}\in U_e(T)$, respectively. 
Hence, the desired state $\hat{y}_{d,e}(\mu)$ admits an affine representation based on the decomposition of $y_{d,e}$ in \eqref{eq:Affine_DesiredState_infty}, 
\begin{align}\label{eq:Affine_DesiredStateTrafo}
\hat{y}_{d,e}(x,t;\mu) = \sum_{q=1}^{Q_{y_d}} \hat{\Theta}_{y_d}^q(t;\mu) y^q_{d,e}(x) ,
\end{align}
with $\hat{\Theta}_{y_d}^q(t;\mu) = e^{-\delta_e t}\Theta_{y_d}^q(t;\mu). $
The now parameter dependent desired control $\hat{u}_{d,e}$ can be decomposed in the obvious way. In the presence of control constraints we set $\hat{u}_{a,e}(t):= e^{-\delta_e t} u_{a,e}(t)$, $\hat{u}_{b,e}(t):= e^{-\delta_e t} u_{b,e}(t)$. We denote the associated control space by $\widehat{U}_{e,ad}(T) = \{ u_e \in U_e(T) \ :\ \hat{u}_{a,e}(t) \le u_e(t) \le \hat{u}_{b,e}(t) \} \subset U_e(T) $. 

The problem \eqref{eq:TransformedOCP} is well-defined for every $\mu\in \cD$, i.e., there exists a unique solution $(\hat{y}^\ast_e,\hat{u}^\ast_e)$. 
The associated first-order optimality system is derived analogously to \eqref{eq:OCPSystem} and reads as follows: 
Given $\mu\in \cD$, the optimal solution to \eqref{eq:TransformedOCP}, $(\hat{y}_e^\ast, \hat{p}_e^\ast, \hat{u}_e^\ast)\in W(0,T)\times W(0,T)\times \widehat{U}_{e,ad}(T) $ satisfies
\begin{subequations}\label{eq:TransformedOCPSystem}
 	\begin{alignat}{2}
		 \dfrac{\diff}{\diff t} m(\hat{y}_e^\ast(t),\phi) + a_\delta(\hat{y}_e^\ast(t),\phi;\mu) &=b(\hat{u}_e^\ast(t),\phi)  \quad 					
		 		&& \forall \phi \in Y_e \quad \textrm{f.a.a. } t\in (0,T], \label{eq:TransformedOCPSystem_StateA}\\
 		\hat{y}_e^\ast(0) &= y_{0,e}, \label{eq:TransformedOCPSystem_StateB} \\
 		-\dfrac{\diff}{\diff t}m(\varphi, \hat{p}_e^\ast(t)) + a_\delta(\varphi, \hat{p}_e^\ast(t);\mu) &= \hat{\sigma}_1(t)(\hat{y}_{d,e}(t;\mu) - \hat{y}_e^\ast(t),   \varphi)_D   \quad 
 				&& \forall \varphi \in Y_e \quad \textrm{f.a.a. } t\in [0,T), \label{eq:TransformedOCPSystem_AdjA}\\
 		m(\varphi, \hat{p}_e^\ast(T)) &=  \hat{\sigma}_2(T)(\hat{y}_{d,e}(T;\mu) - \hat{y}_e^\ast(T),\varphi)_D \quad 
 				&& \forall \varphi \in Y_e , \label{eq:TransformedOCPSystem_AdjB}\\
		 (\hat{\lambda} (\hat{u}_e^\ast - \hat{u}_{d,e}) - \cB_e^\star \hat{p}_e^\ast, \psi - \hat{u}_e^\ast)_{U_e} &\ge 0 \quad 
		 		&& \forall \psi\in \widehat{U}_{e,ad}(T). \label{eq:TransformedOCPOptimalityInequality}
  \end{alignat}
\end{subequations}

It directly follows that $(y_e,p_e,u_e)$ solves \eqref{eq:ExactOCP} if and only if $(\hat{y}_e,\hat{p}_e,\hat{u}_e)$ solves \eqref{eq:TransformedOCP}.

\subsection{Reduced Basis Approximation and {\it A Posteriori} Error Estimation}

Since the transformed optimal control problem~\eqref{eq:TransformedOCP} satisfies the assumptions of Section~\ref{sec:PPOCP_problem_definition}, we can directly apply the reduced basis approximation and associated {\it a posteriori} error bounds from Section~\ref{sec:RB_Approx} as well as the RB-MPC approach from Section~\ref{sec:Stability_RB_MPC}. We note, however, that -- given e.g.\ estimates for the primal (coercive) solution $\hat{y}_e(x,t)$ -- the corresponding estimates for the (noncoercive) solution $y_e(x,t)$ exhibit an additional multiplicative factor $e^{\delta_e t}$. We refer the interested reader to~\cite{UP14}, where certified reduced basis methods for space-time variational formulations of noncoercive parabolic PDEs have been analyzed. The exponential growth of the {\it a posteriori} error bounds is certainly not suitable in long-time integration. In the MPC context, on the other hand, the prediction horizon is often quite small and the exponential growth --  although not desirable -- thus very limited. We consider a numerical example in the next section showing that the growth poses no significant detriment to our approach.

Since the derivation of the following results is analogous to the derivations in Sections~~\ref{sec:RB_Approx} and \ref{sec:Stability_RB_MPC} we omit the details. Note that, due to the time-dependence of the optimization parameters, we use lower and upper bounds for $\hat{\lambda}(t)$ and $\hat{\sigma}_1(t)$ in the error bounds, i.e., we introduce
\begin{align}
\hat{\lambda}_{\min} := \lambda \le e^{2\delta t} \lambda = \hat{\lambda}(t),  \quad 
\hat{\sigma}_1(t) = e^{2\delta t} \sigma_1 \le e^{2\delta T}\sigma_1 =: \hat{\sigma}_{1,\max} ,\quad \forall t\in [0,T].
\end{align}
We only present the results for the error bound in the optimal control and the optimal cost functional. Concerning the former, we directly obtain from the definition of the control energy norm and the temporal transformation for $t = K \, \tau$ that
\begin{align}
\|u\|_U \le e^{\delta \, K \, \tau } \|\hat{u}\|_U, 
\end{align}
which leads us directly to an estimate for the optimal control error bound for the weakly coercive case by invoking Proposition~\ref{prop:OptControlErrBound}. 
\begin{prpstn}\label{prop:ExpOptControlErrBound}
Let $u^\ast$ and $\uN$ be the optimal solutions of the truth and reduced basis optimal control problems \eqref{eq:FE_OCP} and \eqref{eq:RB_OCP} and let $\hat{u}^\ast$ and $\hat{u}_N^\ast$ be the solutions to the corresponding transformed problems, respectively. Given $\mu \in \cD$,  the error in the optimal control satisfies
\begin{align}
\|u^\ast - u_N^\ast \|_{U} \le e^{\delta \, K \, \tau} \|\hat{u}^\ast-\hat{u}^\ast_N\|_U \le e^{\delta \, K \, \tau} \hat{\Delta}^{u,\ast}_{N}(\mu),
\end{align}
where $\hat{\Delta}^{u,\ast}_{N}(\mu):= c_1(\mu) + \sqrt{c_1(\mu)^2 +c_2(\mu)} $ with nonnegative coefficients
\begin{align}
c_1(\mu) &=  \dfrac{1}{\sqrt{2}\alphaLB\hat{\lambda}_{\min}} \left( \sum_{i=1}^m \|b_i\|_{Y'}^2 \right)^{\frac{1}{2}} R_p  ,\\
c_2(\mu) &= \dfrac{1}{\hat{\lambda}_{\min}}\left[ \left(
 \dfrac{2\sqrt{2}}{\alphaLB}R_y 
	+ \dfrac{1+\sqrt{2}}{\sqrt{\alphaLB}} R_0 \right)R_p    
	  + \left( \dfrac{C_D^2\hat{\sigma}_{1,\max}}{\alphaLB} + \dfrac{\hat{\sigma}_2(T)}{2}\right) R_0^2  
	+ \left( \dfrac{C_D^2\hat{\sigma}_{1,\max}}{\alphaLB^2} 
	+ \dfrac{\hat{\sigma}_2(T)}{2\alphaLB}\right) R_y^2 
\right].
\end{align}
The shortcuts $R_p$ $R_y$ and $R_0$ are given by
\begin{align}
R_p = \left( \tau \sum_{k=1}^K  \| \hat{r}_p^{k} (\cdot;\mu) \|^2_{Y'}\right)^{\frac{1}{2}},
\quad R_y = \left( \tau \sum_{k=1}^K  \| \hat{r}_y^{k} (\cdot;\mu) \|^2_{Y'}\right)^{\frac{1}{2}},
 \quad R_0  = \norm{y_0 - y_{N,0} }_{L^2(\Omega)} .
\end{align}
\end{prpstn}
 
Considering the transformations $\hat{y}^{\ast,k}:= e^{-\delta k \tau} y^{\ast,k}$, $\hat{u}^{\ast,k}:= e^{-\delta k \tau} u^{\ast,k}$ and the definition of the cost functionals $J$ and $\hat{J}$ we notice that
\begin{align}
J(y^\ast, u^\ast;\mu) = \hat{J}(\hat{y}^\ast,\hat{u}^\ast;\mu).
\end{align}
The same result holds for the cost functional values of the corresponding reduced basis optimal solutions. Hence, we can avoid an additional exponential factor in the cost functional error bounds. Here, we just invoke the estimates of Proposition \ref{prop:OptCostInitErrBound}. 
 
\begin{prpstn}
Let $J^\ast = J(y^\ast, \uast;\mu)$ and $J_N^\ast = J(y_N^\ast, \uN;\mu)$ be the optimal values of the cost functionals of the truth and reduced basis optimal control problems, respectively. The error then satisfies
\begin{equation}
\begin{split}
J^\ast - J^\ast_N 
&\le \hat{\Delta}^{J,\ast}_{N}(y_0;\mu) 
		, \quad \forall \mu \in \cD,
\end{split}
\end{equation}
where $\hat{\Delta}^{J,\ast}_{N}(y_0;\mu)$ as in \eqref{eq:OptCostInitErrBound_C}.
If we consider the problem without control constraints it holds 
\begin{equation}
\begin{split}
|J^\ast - J^\ast_N |
\le \hat{\Delta}^{J,\ast}_{N}(y_0;\mu) 
		, \quad \forall \mu \in \cD.
\end{split}
\end{equation}
with $\hat{\Delta}^{J,\ast}_{N}(y_0;\mu)$ as in \eqref{eq:OptCostInitErrBound_UC}.
\end{prpstn} 

\section{Numerical Results}\label{sec:NumExp}

\subsection{One-dimensional Reaction-Diffusion Equation}\label{sec:1DHeat}

We first test our approach on a model problem that was already studied in \cite{AG12,AG13} with the classical MPC approach. We consider the one dimensional linear reaction-diffusion equation with Neumann boundary control on $\Omega=(0,1)$ 
given by 
\begin{equation}
\begin{aligned}\label{eq:heat_eq_dirichlet}
y_t(x,t) &= y_{xx}(x,t) + \mu_1 y(x,t)  \quad &&\textrm{on }\Omega\times (0,\infty),\\
y(x,0) & = y_0(x) 		  &&\textrm{on } \Omega,\\
y(0,t)&=0  				  &&\textrm{on }(0,\infty),\\
y_x(1,t) &=u(t)				 &&\textrm{on }(0,\infty),
\end{aligned}
\end{equation}
where $y_0\in H^1_0(\Omega)$; in particular, we set $ y_0(x)= \frac{1}{5}\sin(\pi x)$. Here, we impose a homogeneous Dirichlet condition on the left boundary and the control acts only on the right boundary. It is well-known that the uncontrolled equation, i.e., setting $u(t)\equiv 0$, is unstable if $\mu_1\ge \pi^2/4$ for $\Omega=(0,1)$~\cite{KS08}.



Our goal is to stabilize the system using our RB-MPC approach, i.e., we wish to steer the system towards the equilibrium point $y_d\equiv 0$ using the RB-MPC feedback law from Section \ref{sec:RB-MPC_Alg}. We also compare the computed minimal stabilizing horizons with numerical and theoretical results for the classical MPC controller~\cite{GPSW10, AG12, AG13} which serves as our reference.
To this end, we briefly state the truth approximation of the problem and all necessary definitions for the reduced basis method. 

We use Euler-Backward for time integration with time step $\tau=0.01$ and linear finite elements for the truth approximation subspace $Y\subset Y_e=\{v\in H^1(\Omega) \ | \ v(0)=0\}$ over a triangulation of $\Omega$ (dimension $\cN=199$). The field variable $y^k(\mu)\in Y$ satisfies the governing equation in \eqref{eq:FE_OCP}. The bilinear and linear forms are 
\begin{align}
m(w,v) = \int_0^1 w\ v \diff x, 
\quad a(w,v;\mu) =  \int_0^1   w_x\   v_x \diff x - \mu_1 \int_0^1 w \ v \diff x ,
\quad b(v) = v(1)  .
\end{align}
The bilinear form $a(\cdot, \cdot;\mu)$ admits the affine representation \eqref{eq:AffineRepr_a} with $\Theta_a^1(\mu)=1$, $\Theta_a^2(\mu)=-\mu_1$ and $Q_a=2$. Moreover, we note that $a(\cdot, \cdot;\mu)$ is only weakly coercive. In particular, the G\r{a}rding inequality \eqref{eq:Garding_inequality} holds with $\delta(\mu)=\mu_1\ge 0$ and $\alpha(\mu)=1 \left( = \alphaLB \right)>0$.  
Hence, we utilize the transformation and the error bounds from Section \ref{sec:Weakly_Coercive}. For the cost functional \eqref{eq:FE_cost_functional} we choose $D = \Omega$ and fix the regularization parameters $\sigma_1 = 1$ and $\sigma_2 = 0$. Since we want to investigate the performance of the RB-MPC approach for different values of $\lambda$ we consider the four values $\lambda = 10^{-1}, \, 10^{-2}, \, 10^{-3}$, and $10^{-4}$. Furthermore, we set $u^k_d\equiv 0$ and the desired state $y^k_d \equiv 0$, $k\in \xN$, resulting in the quadratic cost functional
\begin{align}
J(y, u;\mu) 
 = \frac{1}{2} \tau \sum_{k=1}^K  \|y^k\|^2_{L^2(\Omega)} + \frac{\lambda}{2} \tau \sum_{k=1}^K  |u^k|^2 
\end{align}
for an optimization horizon $K\ge 1$. Besides the regularization parameter $\lambda$ we allow the parameter $\mu_1$ in  the reaction term to vary. Our parameter and parameter domain is thus $\mu = (  \mu_1, \lambda) \in \cD^i = [1,15] \times \{10^{-i}\}$, $i = 1, \ldots, 4$.
 
We consider the control unconstrained case in this section and construct the reduced basis space $Y_N$ using the POD/Greedy procedure proposed for optimal control problems in Section \ref{sec:Greedy}.
During the sampling we assume $K=20$, corresponding to the finite time interval $[0,T]=[0,0.2]$, and we sample on $\Delta_N(\mu)= \hat{\Delta}^{J,\ast}_N(\mu)/J_N^\ast(\mu)$.  

\pgfplotstableset{
	columns/N/.style={fixed, column type/.add={}{|},column name=$N$},
	every last row/.style={after row= \hline}
}

\subsubsection{Effectivity}
\label{sec:NumExp_Effectivity}


We present numerical results for the maximum relative errors, bounds and average effectivities for the cost functional. We recall that the optimal value function and the {\it a posteriori} error bound are the crucial ingredients entering  the RB-MPC stability analysis. We first investigate the influence of the regularization parameter $\lambda$ on the effectivity of the error bound. To this end, we define the set $\Xi_{\mu_1, train}\subset [1,15]$ consisting of $|\Xi_{\mu_1, train}|=20$ uniformly distributed elements in $\mu_1$ and specify training sets $\Xi^i_{train}= \Xi_{\mu_1, train} \times \{10^{-i}\}$ for $i = 1, 2, 3, 4$. We thus only vary the first parameter in this setting while keeping the second parameter fixed. We use $y_0\neq0$ to initialize the Greedy sampling procedure and construct a reduced basis space for each training sample. Next, we introduce parameter test sets $\Xi^i_{test}=\Xi_{\mu_1, test} \times \{10^{-i}\}$, $i = 1,2,3,4$, with $\Xi_{\mu_1, test}\subset [1,15]$, $|\Xi_{\mu_1, test}|=30$, and uniformly distributed parameter points in $\mu_1$. We define the maximum relative optimal errors and bounds as
\begin{align}
\epsilon^{J,\ast,i}_{N,\textrm{max},\textrm{rel}} \equiv \max_{\mu\in \Xi^i_{test}}\dfrac{|J^\ast(\mu)-J^\ast_N(\mu)|}{J^\ast_{\max}}, \quad
\Delta^{J,\ast,i}_{N,\textrm{max},\textrm{rel}} \equiv \max_{\mu\in \Xi^i_{test}}\dfrac{\hat{\Delta}^{J,\ast}_{N}(\mu)}{J^\ast_{\max}},
\end{align} 
where $J^\ast_{\max} = \max_{i=1,2,3,4} \max_{\mu\in \Xi^i_{test}} J^\ast(\mu) = 9.29\cdot 10^{-3}$, 
\begin{align}
\epsilon^{u, \ast,i}_{N,\textrm{max},\textrm{rel}} \equiv \max_{\mu\in \Xi^i_{test}}\dfrac{\|u^\ast(\mu)-u^\ast_N(\mu)\|_{U}}{u^\ast_{\max}}, \quad
\Delta^{u, \ast,i}_{N,\textrm{max},\textrm{rel}} \equiv \max_{\mu\in \Xi^i_{test}}\dfrac{e^{\delta K \tau}\hat{\Delta}^{u,\ast}_{N}(\mu)}{u^\ast_{\max}}.
\end{align}   
where $u^\ast_{\max} =  \max_{i=1,2,3,4} \max_{\mu\in \Xi^i_{test}}  \|u^\ast(\mu)\|_{U} = 0.785$, 
and 
\begin{align}
\epsilon^{y, \ast,i}_{N,\textrm{max},\textrm{rel}} \equiv \max_{\mu\in \Xi^i_{test}}\dfrac{\interleave y^\ast(\mu)-y^\ast_N(\mu)\interleave _{y}}{y^\ast_{\max}}, \quad
\Delta^{y, \ast,i}_{N,\textrm{max},\textrm{rel}} \equiv \max_{\mu\in \Xi^i_{test}}\dfrac{e^{\delta K \tau}\hat{\Delta}^{y,\ast,K}_{N}(\mu)}{y^\ast_{\max}},
\end{align}   
where $y^\ast_{\max} =  \max_{i=1,2,3,4} \max_{\mu\in \Xi^i_{test}} \interleave y^{\ast,K}(\mu)\interleave _{y}  = 0.368$. 
Here, error bounds denoted with a hat refer to the quantities of the transformed problem, see Section \ref{sec:Weakly_Coercive}. Note that no additional exponential factor appears in the error bound of the cost functional.
We also define the average effectivity $\eta^{J,\ast,i}_{\textrm{ave}}$ as the average over $\Xi^i_{test}$ of $\hat{\Delta}^{J,\ast}_{N}(\mu)/|J^\ast(\mu)-J^\ast_N(\mu)|$, the average effectivities for the state and control are defined similarly.
The defined quantities are presented in Table \ref{tab:Eta_CostFunc} as a function of $N$ for each test set $\Xi^i_{test}$. 
We observe a very fast decay in the maximum relative errors and bounds as $N$ increases, but also a significant overestimation of the actual error. Whereas the effectivities remain constant for the state and cost as $N$ increases, the effectivities for the control grow considerably with $N$. Also, the effectivities for the cost functional deteriorate as $\lambda$ decreases, see~\cite{KTGV18} for a discussion of this behaviour in the elliptic case. Although not visible from the table, we note that for fixed $\lambda$ the effectivity increases with $\mu_1$, i.e., as the problem becomes more unstable. Nevertheless, we do overestimate the cost functional error significantly even for larger $\lambda$. However, there are two remedies which help mitigate the large effectivities: the growth of the effectivity with the prediction horizon $K$ and the dependence of the minimal stabilizing horizon on $\lambda$.


\begin{table}[t]
\pgfplotstableread{tables/data_eta_lambda1_update.dat}\dataA
\pgfplotstableread{tables/data_eta_lambda2_update.dat}\dataB
\pgfplotstableread{tables/data_eta_lambda3_update.dat}\dataC
\pgfplotstableread{tables/data_eta_lambda4_update.dat}\dataD

\pgfplotstablecreatecol[copy column from table={\dataB}{[index] 1}] {Dave2} {\dataA}
\pgfplotstablecreatecol[copy column from table={\dataB}{[index] 2}] {errave2} {\dataA}
\pgfplotstablecreatecol[copy column from table={\dataB}{[index] 3}] {etaave2} {\dataA}
\pgfplotstablecreatecol[copy column from table={\dataC}{[index] 1}] {Dave3} {\dataA}
\pgfplotstablecreatecol[copy column from table={\dataC}{[index] 2}] {errave3} {\dataA}
\pgfplotstablecreatecol[copy column from table={\dataC}{[index] 3}] {etaave3} {\dataA}
\pgfplotstablecreatecol[copy column from table={\dataD}{[index] 1}] {Dave4} {\dataA}
\pgfplotstablecreatecol[copy column from table={\dataD}{[index] 2}] {errave4} {\dataA}
\pgfplotstablecreatecol[copy column from table={\dataD}{[index] 3}] {etaave4} {\dataA}

\begin{center}
\scalebox{0.95}{
	\pgfplotstabletypeset[font=\footnotesize, 
		every head row/.style={
			before row={
				\multicolumn{1}{c}{} & \multicolumn{3}{c}{$\lambda= 10^{-1}$} & \multicolumn{3}{c}{$\lambda= 10^{-2}$} 
				& \multicolumn{3}{c}{$\lambda= 10^{-3}$}
				& \multicolumn{3}{c}{$\lambda= 10^{-4}$}
				\\
				\hline & & & & &   & & & & &   & & \\[-7pt]
			},
			after row= \hline & & & & &   & & & & &   & & \\[-7pt]
		},
		columns={N, err-rel-max1, DeltaJ-rel-max1, eta-ave1, errave2, Dave2, etaave2, errave3, Dave3, etaave3, errave4, Dave4, etaave4
		},
	 columns/N/.style={fixed, column type/.add={}{|},column name=$N$},
	 columns/DeltaJ-rel-max1/.style={column name = $\Delta^{J,\ast,1}_{N,\textrm{rel,max}}$, sci , sci E, sci zerofill},
	columns/err-rel-max1/.style={column name = $\epsilon^{J,\ast,1}_{N,\textrm{rel,max}}$, sci, sci E,
	sci zerofill},
	columns/eta-ave1/.style={fixed,column type/.add={}{|}, column name = $\eta^{J,\ast,1}_{\textrm{ave}}$, sci, sci E,
	sci zerofill},
		columns/Dave2/.style={column name = $\Delta^{J,\ast,2}_{N,\textrm{rel,max}}$, sci, sci E, sci zerofill},
	columns/errave2/.style={column name = $\epsilon^{J,\ast,2}_{N,\textrm{rel,max}}$, sci, sci E, sci zerofill},
	columns/etaave2/.style={column name = $\eta^{J,\ast,2}_{\textrm{ave}}$, sci, sci E,
	sci zerofill, 
	column type/.add={}{|}},
%
	columns/Dave3/.style={column name = $\Delta^{J,\ast,3}_{N,\textrm{rel,max}}$, sci , sci E, sci zerofill},
	columns/errave3/.style={column name = $\epsilon^{J,\ast,3}_{N,\textrm{rel,max}}$, sci, sci E, sci zerofill},
	columns/etaave3/.style={column name = $\eta^{J,\ast,3}_{\textrm{ave}}$, sci, sci E,
	sci zerofill, 
	column type/.add={}{|}},
	columns/Dave4/.style={column name = $\Delta^{J,\ast,4}_{N,\textrm{rel,max}}$, sci , sci E, sci zerofill},
	columns/errave4/.style={column name = $\epsilon^{J,\ast,4}_{N,\textrm{rel,max}}$, sci, sci E,
	sci zerofill},
	columns/etaave4/.style={column name = $\eta^{J,\ast,4}_{\textrm{ave}}$, sci, sci E,
	sci zerofill, 
	},
	]{\dataA}	
}	
\end{center}
\caption{Maximum relative optimal cost functional errors and error bounds and average effectivities over the test samples $\Xi^i_{test}$ as a function of $N$ for various regularization parameters.}
\label{tab:Eta_CostFunc}
\end{table}

\begin{table}[t]
\pgfplotstableread{tables/data_etaU_lambda1_update.dat}\dataA
\pgfplotstableread{tables/data_etaU_lambda2_update.dat}\dataB
\pgfplotstableread{tables/data_etaU_lambda3_update.dat}\dataC
\pgfplotstableread{tables/data_etaU_lambda4_update.dat}\dataD

\pgfplotstablecreatecol[copy column from table={\dataB}{[index] 1}] {Dave2} {\dataA}
\pgfplotstablecreatecol[copy column from table={\dataB}{[index] 2}] {errave2} {\dataA}
\pgfplotstablecreatecol[copy column from table={\dataB}{[index] 3}] {etaave2} {\dataA}
\pgfplotstablecreatecol[copy column from table={\dataC}{[index] 1}] {Dave3} {\dataA}
\pgfplotstablecreatecol[copy column from table={\dataC}{[index] 2}] {errave3} {\dataA}
\pgfplotstablecreatecol[copy column from table={\dataC}{[index] 3}] {etaave3} {\dataA}
\pgfplotstablecreatecol[copy column from table={\dataD}{[index] 1}] {Dave4} {\dataA}
\pgfplotstablecreatecol[copy column from table={\dataD}{[index] 2}] {errave4} {\dataA}
\pgfplotstablecreatecol[copy column from table={\dataD}{[index] 3}] {etaave4} {\dataA}

\begin{center}
\scalebox{0.95}{
	\pgfplotstabletypeset[font=\footnotesize, 
		every head row/.style={
			before row={
				\multicolumn{1}{c}{} & \multicolumn{3}{c}{$\lambda= 10^{-1}$} & \multicolumn{3}{c}{$\lambda= 10^{-2}$} 
				& \multicolumn{3}{c}{$\lambda= 10^{-3}$}
				& \multicolumn{3}{c}{$\lambda= 10^{-4}$}
				\\
				\hline & & & & &   & & & & &   & &  \\[-7pt]
			},
			after row= \hline & & & & &   & & & & &   & &   \\[-7pt]
		},
		columns={N, errU-rel-max1, DeltaU-rel-max1, etaU-ave1, errave2, Dave2, etaave2, errave3, Dave3, etaave3, errave4, Dave4, etaave4
		},
	 columns/N/.style={fixed, column type/.add={}{|},column name=$N$},
	 columns/DeltaU-rel-max1/.style={column name = $\Delta^{u,\ast,1}_{N,\textrm{rel,max}}$, sci , sci E, sci zerofill},
	columns/errU-rel-max1/.style={column name = $\epsilon^{u,\ast,1}_{N,\textrm{rel,max}}$, sci, sci E,
	sci zerofill},
	columns/etaU-ave1/.style={fixed,column type/.add={}{|}, column name = $\eta^{u,\ast,1}_{\textrm{ave}}$, sci, sci E,
	sci zerofill},
		columns/Dave2/.style={column name = $\Delta^{u,\ast,2}_{N,\textrm{rel,max}}$, sci, sci E, sci zerofill},
	columns/errave2/.style={column name = $\epsilon^{u,\ast,2}_{N,\textrm{rel,max}}$, sci, sci E, sci zerofill},
	columns/etaave2/.style={column name = $\eta^{u,\ast,2}_{\textrm{ave}}$, sci, sci E,
	sci zerofill, 
	column type/.add={}{|}},
	columns/Dave3/.style={column name = $\Delta^{u,\ast,3}_{N,\textrm{rel,max}}$, sci , sci E, sci zerofill},
	columns/errave3/.style={column name = $\epsilon^{u,\ast,3}_{N,\textrm{rel,max}}$, sci, sci E, sci zerofill},
	columns/etaave3/.style={column name = $\eta^{u,\ast,3}_{\textrm{ave}}$, sci, sci E,
	sci zerofill, 
	column type/.add={}{|}},
	columns/Dave4/.style={column name = $\Delta^{u,\ast,4}_{N,\textrm{rel,max}}$, sci , sci E, sci zerofill},
	columns/errave4/.style={column name = $\epsilon^{u,\ast,4}_{N,\textrm{rel,max}}$, sci, sci E,
	sci zerofill},
	columns/etaave4/.style={column name = $\eta^{u,\ast,4}_{\textrm{ave}}$, sci, sci E,
	sci zerofill, 
	},
	]{\dataA}	
}	
\end{center}
\caption{Maximum relative optimal control errors and error bounds and average effectivities over the test samples $\Xi^i_{test}$ as a function of $N$ for various regularization parameters.}
\label{tab:Eta_Control}
\end{table}


\begin{table}[t]
\pgfplotstableread{tables/data_etaY_lambda1_update.dat}\dataA
\pgfplotstableread{tables/data_etaY_lambda2_update.dat}\dataB
\pgfplotstableread{tables/data_etaY_lambda3_update.dat}\dataC
\pgfplotstableread{tables/data_etaY_lambda4_update.dat}\dataD

\pgfplotstablecreatecol[copy column from table={\dataB}{[index] 1}] {Dave2} {\dataA}
\pgfplotstablecreatecol[copy column from table={\dataB}{[index] 2}] {errave2} {\dataA}
\pgfplotstablecreatecol[copy column from table={\dataB}{[index] 3}] {etaave2} {\dataA}
\pgfplotstablecreatecol[copy column from table={\dataC}{[index] 1}] {Dave3} {\dataA}
\pgfplotstablecreatecol[copy column from table={\dataC}{[index] 2}] {errave3} {\dataA}
\pgfplotstablecreatecol[copy column from table={\dataC}{[index] 3}] {etaave3} {\dataA}
\pgfplotstablecreatecol[copy column from table={\dataD}{[index] 1}] {Dave4} {\dataA}
\pgfplotstablecreatecol[copy column from table={\dataD}{[index] 2}] {errave4} {\dataA}
\pgfplotstablecreatecol[copy column from table={\dataD}{[index] 3}] {etaave4} {\dataA}

\begin{center}
\scalebox{0.95}{
	\pgfplotstabletypeset[font=\footnotesize, 
		every head row/.style={
			before row={
				\multicolumn{1}{c}{} & \multicolumn{3}{c}{$\lambda= 10^{-1}$} & \multicolumn{3}{c}{$\lambda= 10^{-2}$} 
				& \multicolumn{3}{c}{$\lambda= 10^{-3}$}
				& \multicolumn{3}{c}{$\lambda= 10^{-4}$}
				\\
				\hline & & & & &   & & & & &   & &  \\[-7pt]
			},
			after row= \hline & & & & &   & & & & &   & & \\[-7pt]
		},
		columns={N, errY-rel-max1, DeltaY-rel-max1, etaY-ave1, errave2, Dave2, etaave2, errave3, Dave3, etaave3, errave4, Dave4, etaave4
		},
	 columns/N/.style={fixed, column type/.add={}{|},column name=$N$},
	 columns/DeltaY-rel-max1/.style={column name = $\Delta^{y,\ast,1}_{N,\textrm{rel,max}}$, sci , sci E, sci zerofill},
	columns/errY-rel-max1/.style={column name = $\epsilon^{y,\ast,1}_{N,\textrm{rel,max}}$, sci, sci E,
	sci zerofill},
	columns/etaY-ave1/.style={fixed,column type/.add={}{|}, column name = $\eta^{y,\ast,1}_{\textrm{ave}}$, sci, sci E,
	sci zerofill},
		columns/Dave2/.style={column name = $\Delta^{y,\ast,2}_{N,\textrm{rel,max}}$, sci, sci E, sci zerofill},
	columns/errave2/.style={column name = $\epsilon^{y,\ast,2}_{N,\textrm{rel,max}}$, sci, sci E, sci zerofill},
	columns/etaave2/.style={column name = $\eta^{y,\ast,2}_{\textrm{ave}}$, sci, sci E,
	sci zerofill, 
	column type/.add={}{|}},
	columns/Dave3/.style={column name = $\Delta^{y,\ast,3}_{N,\textrm{rel,max}}$, sci , sci E, sci zerofill},
	columns/errave3/.style={column name = $\epsilon^{y,\ast,3}_{N,\textrm{rel,max}}$, sci, sci E, sci zerofill},
	columns/etaave3/.style={column name = $\eta^{y,\ast,3}_{\textrm{ave}}$, sci, sci E,
	sci zerofill, 
	column type/.add={}{|}},
	columns/Dave4/.style={column name = $\Delta^{y,\ast,4}_{N,\textrm{rel,max}}$, sci , sci E, sci zerofill},
	columns/errave4/.style={column name = $\epsilon^{y,\ast,4}_{N,\textrm{rel,max}}$, sci, sci E,
	sci zerofill},
	columns/etaave4/.style={column name = $\eta^{y,\ast,4}_{\textrm{ave}}$, sci, sci E,
	sci zerofill, 
	},
	]{\dataA}	
}
\end{center}
\caption{Maximum relative optimal state errors and error bounds and average effectivities over the test samples $\Xi^i_{test}$ as a function of $N$ for various regularization parameters.}
\label{tab:Eta_State}
\end{table}

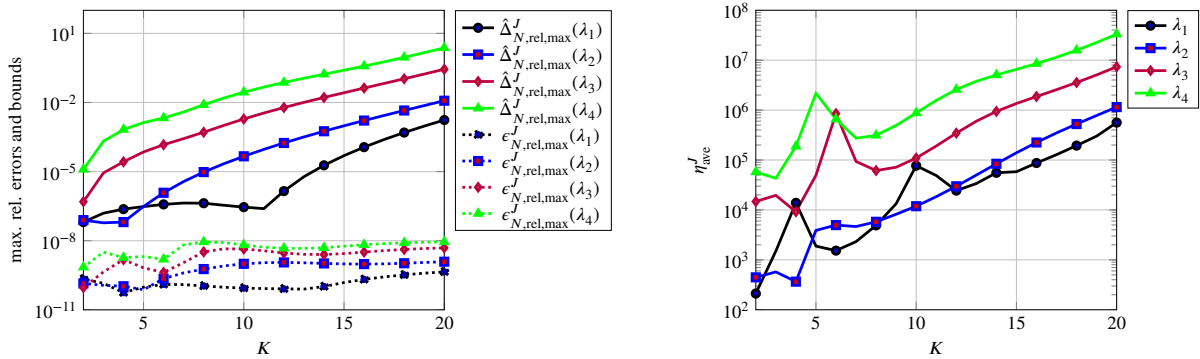
\begin{figure}
\centering
\begin{minipage}[b]{0.49 \textwidth}
\centering
	\begin{tikzpicture}[scale=0.7]
	\begin{axis}[xlabel=$K$, ylabel=max. rel. errors and bounds, xmin=2, xmax=20, ymode=log, ymin=1e-11, ymax=1e+2,
	xmajorgrids,
	ymajorgrids,
	yminorticks=true, 
	axis background/.style={fill=white},
		legend pos=outer north east]
		\addplot+[color=black, line width=1.5pt, mark=*, mark repeat=2, mark phase = 1]
			table[x=K,y=DeltaJ-rel-max, mark repeat=2, mark phase = 1] {tables/data_etaK_lambda1_N7_update.dat};
			\addlegendentry{$\hat{\Delta}^J_{N,\textrm{rel,max}}(\lambda_1)$}
		\addplot+[color=blue, line width=1.5pt,mark=square*, mark repeat=2, mark phase = 1]
			table[x=K,y=DeltaJ-rel-max] {tables/data_etaK_lambda2_N7_update.dat};
			\addlegendentry{$\hat{\Delta}^J_{N,\textrm{rel,max}}(\lambda_2)$}
		\addplot+[color=purple, line width=1.5pt,mark=diamond*, mark repeat=2, mark phase = 1]
			table[x=K,y=DeltaJ-rel-max] {tables/data_etaK_lambda3_N7_update.dat};
			\addlegendentry{$\hat{\Delta}^J_{N,\textrm{rel,max}}(\lambda_3)$}
		\addplot+[color=green,line width=1.5pt,mark= triangle*, mark repeat=2, mark phase = 1]
			table[x=K,y=DeltaJ-rel-max] {tables/data_etaK_lambda4_N7_update.dat};
			\addlegendentry{$\hat{\Delta}^J_{N,\textrm{rel,max}}(\lambda_4)$}
		\addplot+[dotted,color=black,line width=1.5pt,mark=*, mark repeat=2, mark phase = 1]
			table[x=K,y=err-rel-max] {tables/data_etaK_lambda1_N7_update.dat};
			\addlegendentry{$\epsilon^J_{N,\textrm{rel,max}}(\lambda_1)$}
		\addplot+[dotted,color=blue, line width=1.5pt,mark=square*, mark repeat=2, mark phase = 1]
			table[x=K,y=err-rel-max] {tables/data_etaK_lambda2_N7_update.dat};
			\addlegendentry{$\epsilon^J_{N,\textrm{rel,max}}(\lambda_2)$}
		\addplot+[dotted,color=purple, line width=1.5pt,mark=diamond*, mark repeat=2, mark phase = 1]
			table[x=K,y=err-rel-max] {tables/data_etaK_lambda3_N7_update.dat};
			\addlegendentry{$\epsilon^J_{N,\textrm{rel,max}}(\lambda_3)$}
		\addplot+[dotted,color=green,line width=1.5pt,mark= triangle*, mark repeat=2, mark phase = 1]
			table[x=K,y=err-rel-max] {tables/data_etaK_lambda4_N7_update.dat};
			\addlegendentry{$\epsilon^J_{N,\textrm{rel,max}}(\lambda_4)$}
	\end{axis}
	\end{tikzpicture}
\end{minipage}
\vspace*{0.1cm}
\begin{minipage}[b]{0.49 \textwidth}
\centering
	\begin{tikzpicture}[scale=0.7]
	\begin{axis}[xlabel=$K$, ylabel=$\eta^J_{\textrm{ave}}$, xmin=2, xmax=20, ymode=log, ymin=1e+2, ymax=1e+8,
	xmajorgrids,
	ymajorgrids,
	yminorticks=true, 
	axis background/.style={fill=white},
		legend pos=outer north east]
		\addplot+[color=black, line width=1.5pt, mark=*, mark repeat=2, mark phase = 1]
			table[x=K,y=eta-ave] {tables/data_etaK_lambda1_N7_update.dat};
			\addlegendentry{$\lambda_1$}
		\addplot+[color=blue, line width=1.5pt,mark=square*, mark repeat=2, mark phase = 1]
			table[x=K,y=eta-ave] {tables/data_etaK_lambda2_N7_update.dat};
			\addlegendentry{$\lambda_2$}
		\addplot+[color=purple, line width=1.5pt,mark=diamond*, mark repeat=2, mark phase = 1]
			table[x=K,y=eta-ave] {tables/data_etaK_lambda3_N7_update.dat};
			\addlegendentry{$\lambda_3$}
		\addplot+[color=green,line width=1.5pt,mark= triangle*, mark repeat=2, mark phase = 1]
			table[x=K,y=eta-ave] {tables/data_etaK_lambda4_N7_update.dat};
			\addlegendentry{$\lambda_4$}
	\end{axis}
	\end{tikzpicture}
\end{minipage}
\caption{Maximum relative optimal cost functional error and error bound (left) and average effectivity (right) as a function of $K$ for fixed basis size $N=7$.}
\label{fig:Effectivity_overK_fixedN}
\end{figure}

We first comment on the former remedy. Since we are interested in the performance of the bounds for varying prediction horizons, we plot the maximum relative optimal cost functional errors and bounds for fixed basis size $N$ as a function of the optimization horizon $K$ in Figure~\ref{fig:Effectivity_overK_fixedN} on the left. In the right plot we show the corresponding average effectivties. We again consider the four different bases with regularization parameters $\lambda_i = 10^{-i}$ for $i = 1, 2, 3, 4$. Depending on the parameter $\lambda$, we observe exponential growth of the estimator in $K$ resulting in the high effectivities of Table \ref{tab:Eta_CostFunc}.  However, the effectivities are approximately $\cO(10^3)$ smaller --- and the error bounds hence considerably sharper --- for the smallest prediction horizon $K \gtrapprox 1$ compared to the largest one. The fact that the effectivities are larger for smaller $\lambda$ also holds for $K$ small. We discuss the minimal stabilizing horizon and the second remedy in the next section.

\subsubsection{Minimal Stabilizing Horizon}
\begin{table}%
  \centering
  
\pgfplotstableread{tables/data_minh_lambda1_I100.dat}\dataA
\pgfplotstableread{tables/data_minh_lambda2_I100.dat}\dataB
\pgfplotstableread{tables/data_minh_lambda3_I100.dat}\dataC
\pgfplotstableread{tables/data_minh_lambda4_I100.dat}\dataD

\pgfplotstablecreatecol[copy column from table={\dataB}{[index] 1}] {2ave2} {\dataA}
\pgfplotstablecreatecol[copy column from table={\dataB}{[index] 4}] {2ave5} {\dataA}
\pgfplotstablecreatecol[copy column from table={\dataB}{[index] 7}] {2ave8} {\dataA}
\pgfplotstablecreatecol[copy column from table={\dataB}{[index] 10}] {2ave11} {\dataA}
\pgfplotstablecreatecol[copy column from table={\dataB}{[index] 13}] {2ave14} {\dataA}

\pgfplotstablecreatecol[copy column from table={\dataC}{[index] 1}] {3ave2} {\dataA}
\pgfplotstablecreatecol[copy column from table={\dataC}{[index] 4}] {3ave5} {\dataA}
\pgfplotstablecreatecol[copy column from table={\dataC}{[index] 7}] {3ave8} {\dataA}
\pgfplotstablecreatecol[copy column from table={\dataC}{[index] 10}] {3ave11} {\dataA}
\pgfplotstablecreatecol[copy column from table={\dataC}{[index] 13}] {3ave14} {\dataA}

\pgfplotstablecreatecol[copy column from table={\dataD}{[index] 1}] {4ave2} {\dataA}
\pgfplotstablecreatecol[copy column from table={\dataD}{[index] 4}] {4ave5} {\dataA}
\pgfplotstablecreatecol[copy column from table={\dataD}{[index] 7}] {4ave8} {\dataA}
\pgfplotstablecreatecol[copy column from table={\dataD}{[index] 10}] {4ave11} {\dataA}
\pgfplotstablecreatecol[copy column from table={\dataD}{[index] 13}] {4ave14} {\dataA}

\scalebox{0.95}{
\pgfplotstabletypeset[font=\footnotesize, 
		every row no 3/.style={before row=\hline \hline},
		every head row/.style={
			before row={
				\multicolumn{1}{c}{} & \multicolumn{5}{c}{$\lambda= 10^{-1}$}& \multicolumn{5}{c}{$\lambda= 10^{-2}$} 
				& \multicolumn{5}{c}{$\lambda= 10^{-3}$} & \multicolumn{5}{c}{$\lambda= 10^{-4}$}
				\\
				\hline  & & &  & & &  & & &  & & &  & & &  & & & & & \\[-7pt]
				$\mu_1$ & 2& 5 & 8 & 11& 14 &  2& 5 & 8 & 11& 14 
				&  2& 5 & 8 & 11& 14 &  2& 5 & 8 & 11& 14
				 \\
				\hline   & & &  & & &  & & &  & & &  & & &  & & & & &   \\[-7pt]
				$N$ & \multicolumn{5}{c|}{$K_{\textrm{ave}}$} & \multicolumn{5}{c|}{$K_{\textrm{ave}}$} 
				& \multicolumn{5}{c|}{$K_{\textrm{ave}}$} &  \multicolumn{5}{c}{$K_{\textrm{ave}}$} 
				\\
				\hline  & & &  & & &  & & &  & & &  & & &  & & & & &  \\[-7pt]
			},
			output empty row
		},
			columns={N, K-ave2, K-ave5, K-ave8, K-ave11, K-ave14, 2ave2, 2ave5, 2ave8, 2ave11, 2ave14, 3ave2, 3ave5, 3ave8, 3ave11, 3ave14, 4ave2, 4ave5, 4ave8, 4ave11, 4ave14
			 },
			columns/K-ave14/.append style={string type,column type={c|}},
			columns/2ave14/.append style={string type,column type={c|}},
			columns/3ave14/.append style={string type,column type={c|}},
			columns/N/.style={fixed, column type/.add={}{|},column name=$N$},
        	columns/N/.append style={string replace={-1}{}},
        	       empty cells with={$\kappa$}
		]{\dataA}
} 
  \caption{Average minimal stabilizing horizons over $[0,1]$ for different parameters $\mu=(\mu_1,\lambda)$.}%
  \label{tab:MinH_lambda}%
\end{table}

We now investigate the performance of our proposed RB-MPC algorithm. We follow the approach presented in Algorithm~\ref{alg:Adaptive_RB_MPC} with the (known) initial condition $y_0$ and $K_{\textrm{max}} = 20$. We keep the four distinct parameter sets ${\cal D}^i = [1,15] \times \{10^{-i}\}$, $i = 1, \ldots, 4$, and employ the reduced basis approximations generated for these sets in the previous section. We consider the time interval $[0,1]$ corresponding to $100$ iterations in the RB-MPC algorithm. Note that we generated the reduced bases for the shorter interval $[0, 0.2]$. 
We present the average minimal stabilizing horizons (over the $100$ iterations) for different representative parameter values --- $\mu = (\mu_1,\lambda)\in \{2,  \, 5, \, 8, \, 11, \, 14\} \times \{10^{-i}\}$, $i = 1, \ldots, 4$ --- and reduced basis dimensions $N = 5, \, 7,$ and  $9$ in Table~ \ref{tab:MinH_lambda}. Note that we do not present results for $N = 1$ and $3$ because the approximations were not sufficiently accurate for most parameter values, i.e.,  the RB-MPC algorithm did not result in a stabilizing control before reaching the maximal prediction horizon $K_{\textrm{max}} = 20$.  For comparison, we also show the corresponding values for the feedback controller $\kappa$ computed with the classical (high-fidelity) MPC approach in the last row of the table.  

Concerning the feedback law $\kappa$, we observe that the minimal horizon $K=1$ is stabilizing for the (already) stable problem $\mu_1=2$ for all values of $\lambda$.  In the unstable regime $\mu_1\ge \pi^2/4$ the minimal stabilizing horizon decreases with decreasing regularization parameter $\lambda$. These observations are consistent with the results in~\cite{AG13}. For the RB-MPC approach we observe that the average minimal stabilizing horizon is close to the reference values for all $N$  --- especially for $\mu_1$ small --- and approaches the reference value with increasing RB size $N$. For the maximum basis size ($N = 9$) the RB-MPC stabilizing horizons are identical to or only slightly larger than the reference value. 
 
Returning to the second remedy, we recall that we obtained the largest effectivities of the cost functional error bound  for small $\lambda$  and large $K$. We observe from Table~ \ref{tab:MinH_lambda},  however, that for small values of $\lambda$ the minimal stabilizing horizon is also very small, i.e., $K$ is often  close to 1. Due to the correlation between $\lambda$ and the prediction horizon $K$, the worst case scenario --- small $\lambda$ {\it and} large $K$ --- does not appear and thus mitigates the exponential growth of the effectivity. 

To summarize, we note that the RB-MPC algorithm performs very well thanks to the exponential convergence of the cost functional error bound despite the fairly large effectivities. We can replicate the results of the classical MPC approach as soon as the reduced basis approximation is sufficiently accurate, in the present example already for fairly small $N$.


\subsubsection{Suboptimality Degree}
\label{sec:Num_Exp_Alpha}

The suboptimality degree $\tilde{\omega}_{N,K}(\mu)$ defined in~\eqref{eq:omegatilde_c} is the key ingredient in order to gauge stability of the RB-MPC controller. We therefore study the behaviour of $\tilde{\omega}_{N,K}(\mu)$ as a function of $K$ and compare it to the suboptimality degree of the FE-MPC feedback law $\kappa(y_0) = u^{\ast,1}$. To this end, we define the suboptimality degree
\begin{align}
\tilde{\omega}_{K}(\mu) = \dfrac{J^{\ast}(y_{0};\mu)  -  J^{\ast}(y^{\ast,1};\mu) }{\dfrac{\sigma_1}{2} \tau| y_0 - y_{d}^0(\mu)|^2_{D} +   \dfrac{\lambda}{2} \tau  \| u^{\ast,1} - u_{d}^{1}\|^2_{\xR^m}},
\end{align}
for $y_0 \in Y$. We use the same RB approximations as in the last sections.

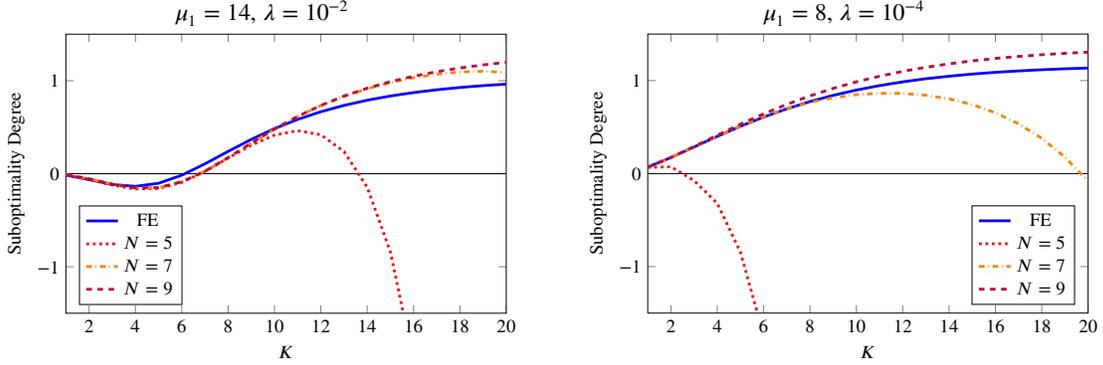
\begin{figure}[t]
\centering
\begin{minipage}[b]{0.45 \textwidth}
	\begin{center}
		\small
		$\mu_1 = 14$, $\lambda = 10^{-2}$
	\end{center}
	\centering
	\begin{tikzpicture}[scale=0.7]
	\begin{axis}[xlabel=$K$, ylabel=Suboptimality Degree, xmin=1, xmax=20, ymin=-1.5, 		 		ymax=1.5,
	width=1.3\textwidth,
    height=0.9\textwidth,
		legend pos=south west]
		\addplot[mark=none, black, forget plot] coordinates {(0,0) (20,0)};
		\addplot+[color=blue,line width=1.5pt,mark=none]
			table[x=K,y=alphaFE] {figures/plot_alpha_lambda2_delta14.dat};
			\addlegendentry{$\textrm{FE}$}
		\addplot+[dotted,color=red,line width=1.5pt,mark=none]
			table[x=K,y=alpha5] {figures/plot_alpha_lambda2_delta14.dat};
			\addlegendentry{$N=5$}
		\addplot+[dashdotted,color=orange,line width=1.5pt,mark=none]
			table[x=K,y=alpha7] {figures/plot_alpha_lambda2_delta14.dat};
			\addlegendentry{$N=7$}
		\addplot+[dashed,color=purple,line width=1.5pt,mark=none]
			table[x=K,y=alpha9] {figures/plot_alpha_lambda2_delta14.dat};
			\addlegendentry{$N=9$}
	\end{axis}
	\end{tikzpicture}
\end{minipage}
\begin{minipage}[b]{0.45 \textwidth}
	\begin{center}
		\small
		$\mu_1 = 8$, $\lambda = 10^{-4}$
	\end{center}
	\centering
	\begin{tikzpicture}[scale=0.7]
	\begin{axis}[xlabel=$K$, ylabel=Suboptimality Degree, xmin=1, xmax=20, ymin=-1.5, 		 		ymax=1.5,
	width=1.3\textwidth,
    height=0.9\textwidth,
		legend pos=south east]
		\addplot[mark=none, black, forget plot] coordinates {(0,0) (20,0)};
		\addplot+[color=blue,line width=1.5pt,mark=none]
			table[x=K,y=alphaFE] {figures/plot_alpha_lambda4_delta8.dat};
			\addlegendentry{$\textrm{FE}$}
		\addplot+[dotted,color=red,line width=1.5pt,mark=none]
			table[x=K,y=alpha5] {figures/plot_alpha_lambda4_delta8.dat};
			\addlegendentry{$N=5$}
		\addplot+[dashdotted,color=orange,line width=1.5pt,mark=none]
			table[x=K,y=alpha7] {figures/plot_alpha_lambda4_delta8.dat};
			\addlegendentry{$N=7$}
		\addplot+[dashed,color=purple,line width=1.5pt,mark=none]
			table[x=K,y=alpha9] {figures/plot_alpha_lambda4_delta8.dat};
			\addlegendentry{$N=9$}
	\end{axis}
	\end{tikzpicture}
\end{minipage}
\caption{Suboptimality degrees $\tilde{\omega}_{K}(\mu)$ (blue) and $\tilde{\omega}_{N,K}(\mu)$ as a function of $K$ for the initial condition $y_0=y^{10}(u^\ast)$.}
\label{fig:Subopt_lambda2_4}
\end{figure}

To illustrate the behavior of the suboptimality estimates, we choose two representative parameter values from Table~ \ref{tab:MinH_lambda} and $y_0=y^{10}(u^\ast)$. In Figure \ref{fig:Subopt_lambda2_4}, we plot the suboptimality estimates $\tilde{\omega}_{K}(\mu)$ and $\tilde{\omega}_{N,K}(\mu)$ as a function of $K$ for different basis sizes. In the left plot we observe that the suboptimality degree is not necessarily a monotonically increasing function in $K$. The functions intersect the x-axis between $K=6$ and $K=7$. Hence, $K = 7$ is the minimal stabilizing horizon at this point.  For the parameter pair in the right plot $\tilde{\omega}_{K}(\mu)$ and $\tilde{\omega}_{N,K}(\mu)$ are positive for all $K$, i.e., the minimal stabilizing horizon is $K = 1$. We also observe that $\tilde{\omega}_{N,K}(\mu)$ approaches $\tilde{\omega}_{K}(\mu)$ with increasing $N$ in both cases.


\subsection{Control of a Welding Process}\label{sec:NumExp_WP}

\begin{figure}[htb]
\begin{center}
\scalebox{1.5}{
\begin{tikzpicture}
\coordinate (v1) at (3.6,1);
\coordinate (v2) at (3.55,0.5);
\coordinate (v3) at (3,0.5);
\coordinate (v3a) at (1.7,0.4);
\coordinate (v3b) at (0.5,0.7);
\coordinate (v4) at (0,0.45);
\filldraw[gray!20] (v1) .. controls (v2) .. (v3) .. controls  (v3a) and (v3b) .. (v4) -- (0,1) -- cycle;

\draw[->] (0,0) -- (5.5,0) node[right] {\tiny{$x_1$}};
\draw[->] (0,0) -- (0,1.5) node[above] {\tiny{$x_2$}};

\draw[blue,thick] (0,0)--(5,0);
\draw[blue,thick] (0,0)--(0,1);
\draw[thick] (5,0)--(5,1);
\draw[blue, thick] (0,1)--(5,1);
\node (A) at (1,1.2) [blue] {\tiny{$\Gamma_N$}};
\node (B) at (5.3,0.5) [black] {\tiny{$\Gamma_D$}};

\draw (5,-.1) -- (5,.1) node[below=4pt] {\tiny{$\scriptstyle5$}};
\draw (3.5,-.1) -- (3.5,.1) node[below=4pt] {\tiny{$\scriptstyle3.5$}};
\draw (1.5,-.1) -- (1.5,.1) node[below=4pt] {\tiny{$\scriptstyle1.5$}};
\draw (-.1,1) -- (.1,1) node[left=4pt] {\tiny{$\scriptstyle1$}};
\draw (-.1,0.5) -- (.1,0.5) node[left=4pt] {\tiny{$\scriptstyle0.5$}};
\draw (-.1,0) -- (.1,0) ;
\draw (0,-.1) -- (0,.1) ;
\node[below left]{\tiny{$\scriptstyle0$}};

\filldraw[draw=black, fill=black] (1.42,0.47) rectangle (1.58,0.53);		
\node (s) at (1.8,0.32) [black!80] {\tiny{$\Omega_s$}};
\draw [densely dashed] (1.5,0.5) -- (0, 0.5);
\draw [densely dashed] (1.5,0.5) -- (1.5, 0);

\coordinate (d1) at (3.5, 1); 
\coordinate (d2) at (3.4, 1.4); 
\coordinate (d3) at (3.6, 1.4); 
\filldraw[draw=black, fill=red!60] (d1) -- (d2) -- (d3) -- cycle;
\draw [densely dashed] (d1) -- (3.5, 0);

\coordinate (m1) at (3.7, 1.2); 
\coordinate (m2) at (4.3, 1.2);
\draw [->] (m1) -- (m2);
\node (V) at (4.2, 1.4) [black] {\tiny{Pe}};
\end{tikzpicture}
}
\end{center}
\caption{Sketch of the domain, the position of the laser and the observation domain.}
\label{fig:WeldingPlate_Sketch}
\end{figure}
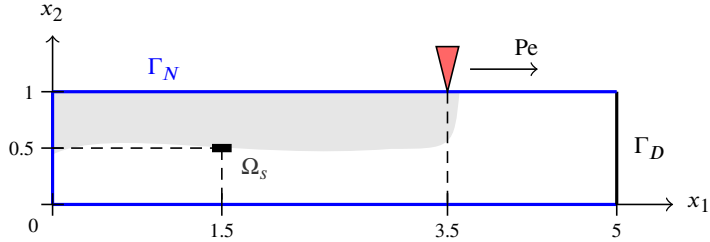

Based on \cite{KGV12}, we consider the problem of controlling the heat flow in a welding process. The non-dimensionalized temperature distribution within the workpiece is governed by the non-dimensionalized unsteady convection-diffusion equation 
\begin{equation} \label{eq:PDE_WeldingPlate}
\frac{\partial}{\partial t} y(x,t) + \textrm{Pe}\frac{\partial}{\partial x_1} y(x, t) - \mu_1 \nabla^2 y(x, t) = q(x)u(t), 		\quad  x\in \Omega, \ t\in [0,\infty)
\end{equation} 
and initial condition $y(x, 0) = y_{e,0} \equiv 0$,  i.e., we consider the start-up of the welding process. A sketch of the domain and setup is illustrated in Figure \ref{fig:WeldingPlate_Sketch}. The spatial domain is $\Omega= (0,5)\times(0,1)\subset\xR^2$ and a point shall be denoted by $x=(x_1,x_2)\in \Omega$. The heat flux input is modelled as a a Gaussian distribution centered at the torch position $x^T = (3.5,1)$ given by
\begin{equation}
q(x) = \frac{1}{2\pi} e^{-0.5((x_1-x_1^T)^2+ (x_2-x_2^T)^2)}.
\end{equation} 
We assume homogeneous Neumann boundary conditions on $\Gamma_N$ and homogeneous Dirichlet boundary conditions on $\Gamma_D$. The torch moves in the $x_1$-direction with fixed velocity $\textrm{Pe}=3.65$ which enters as a convective term in the equation. We consider the thermal diffusivity $\mu_1$ as a parameter which is allowed to vary in the range ${\cal D} = [0.5,2]$. We use Euler-Backward with a step-size $\tau = 0.02$ for time integration  and linear finite elements for the truth approximation subspace (dimension $\cN=3647$). The bilinear and linear forms are 
\begin{equation}
m(u,v) = \int_\Omega u \ v \diff x, \quad
a(u,v;\mu) = \mu_1 \int_\Omega  \nabla u\  \nabla v \diff x +\int_\Omega  v\  ( w \cdot \nabla u)\diff x  , \quad
b(v, q(x) ) = \int_\Omega q(x)\  v\diff x,
\end{equation}
where $w = [\textrm{Pe}, 0]^T$. Note, that the bilinear form $a(\cdot,\cdot)$ depends affinely on the parameter with $\Theta_a^1(\mu)=\mu_1$, $\Theta_a^2(\mu)=1$ and $Q_a=2$. We define the $Y$-inner product $(u,v)_Y = \frac{1}{2}\left( a(u,v,\bar{\mu}_1) + a(v,u,\bar{\mu}_1) \right)$, $\forall u,v\in Y$, with reference parameter $\bar{\mu}_1=1$. A lower bound for the coercivity constant is thus given by $\alphaLB= \min(\mu_1,1)$. 

We use a multi-step approach to control the heat flux of the welding process, see Remark~\ref{rem:Multistep_Definition}. More precisely, we choose an $n$-step controller with $n=5$ resulting in the control horizon $\taucon = n \, \tau = 0.1$. Since the quality of the weld is measured by means of the welding depth, we define the output of interest as the average temperature over the measurement domain $\Omega_{s} = [1.42, 1.58]\times [0.47, 0.53]$ given by $s(t^k,\mu) = |\Omega_{s}|^{-1} \int_{\Omega_{s}} y(x,t^k;\mu) \diff x$.  We choose $D = \Omega_s$ and the regularization parameters are set to $\sigma_1 = 10$ and $\sigma_2 = 10$. The parameter $\lambda$ is allowed to vary and considered an additional input parameter. The parameter input and domain are hence given by $\mu = ( \mu_1, \lambda) \in \cD = [0.5,2] \times [10^{-6}, 10^{-4}]$.  Furthermore, we set $u_d\equiv 0$ and the desired state $y^k_d(\mu)$ is chosen such that the output $s\equiv 1$ for all $k\in \xN$ independent of $\mu$.  The quadratic cost functional is thus given by
\begin{align}
J(y, u;\mu) = \frac{\sigma_1}{2} \tau \sum_{k=1}^K  |y^k- y^k_{d}(\mu)|^2_{D} 
+ \dfrac{\sigma_2}{2} | y^K- y^K_{d}(\mu)|^2_{D} 
 + \frac{\lambda}{2} \tau \sum_{k=1}^K  |u^{k}|^2,
\end{align}
with an optimization horizon $K\ge n$. Since we observe the state only on a part of the domain $D=\Omega_s$, the associated running cost $l(\cdot,\cdot;\mu)$ is semidefinite. However, the governing equation~\eqref{eq:PDE_WeldingPlate} is stable for all $\mu_1$ and hence we obtain asymptotic stability based on the discussion in Remark~\ref{rem:Stability_by_Detectability}. Due to physical limitations of the laser we incorporate control constraints: a (natural) lower bound is given by $u_a \equiv 0$,  and $u_b \equiv 80 $ serves as an upper bound. To construct the reduced basis space we use the POD/Greedy sampling procedure. We choose the training sample $\Xi_{train}\subset  \cD$ consisting of $|\Xi_{train}|= 10 \cdot 20  = 200$ elements logarithmically distributed in $\mu_1$ and $\lambda$, $\mu^{1} = (0.5, 10^{-6})\in \Xi_{train}$ as initial parameter value, the time interval $[0,T] = [0, 1.5]$ ($K = 75$), and sample on $\Delta_N(\mu)= \Delta^{u,\ast}_N(\mu)/ \|u^\ast_N(\mu)\|_{U}$.

\subsubsection{Effectivity}\label{sec:NumExp_WP_Effectivity}

We first present results for the maximum relative errors, bounds, and average effectivities. We define
\begin{align}
\epsilon^{J, \ast}_{N,\textrm{max},\textrm{rel}} \equiv \max_{\mu\in \Xi_{test}}\dfrac{|J^\ast(\mu)-J^\ast_N(\mu)|}{J^\ast_{\max}}, \quad
\Delta^{J, \ast}_{N,\textrm{max},\textrm{rel}} \equiv \max_{\mu\in \Xi_{test}}\dfrac{\Delta^{J,\ast}_{N}(\mu)}{J^\ast_{\max}}
\end{align}   
and
\begin{align}
\epsilon^{u, \ast}_{N,\textrm{max},\textrm{rel}} \equiv \max_{\mu\in \Xi_{test}}\dfrac{\|u^\ast(\mu)-u^\ast_N(\mu)\|_{U}}{u^\ast_{\max}}, \quad
\Delta^{u, \ast}_{N,\textrm{max},\textrm{rel}} \equiv \max_{\mu\in \Xi_{test}}\dfrac{\Delta^{u,\ast}_{N}(\mu)}{u^\ast_{\max}}.
\end{align}   
where $J^\ast_{\max} =  \max_{\mu\in \Xi_{test}} J^\ast(\mu)$ and $u^\ast_{\max} = \max_{\mu\in \Xi_{test}} \|u^\ast(\mu)\|_{U}$. 
In addition, we define 
\begin{align}
\epsilon^{y, \ast}_{N,\textrm{max},\textrm{rel}} \equiv \max_{\mu\in \Xi_{test}}\dfrac{\interleave y^\ast(\mu)-y^\ast_N(\mu)\interleave_{y}}{y^\ast_{\max}}, \quad
\Delta^{y, \ast}_{N,\textrm{max},\textrm{rel}} \equiv \max_{\mu\in \Xi_{test}}\dfrac{\Delta^{y,\ast,K}_{N}(\mu)}{y^\ast_{\max}},
\end{align} 
with $y^\ast_{\max} = \max_{\mu\in \Xi_{test}} \interleave y^\ast(\mu) \interleave_{y}$.
The parameter test set $\Xi_{test}\subset \cD$ is of size  $|\Xi_{test}|=6 \cdot 10 = 60$ with logarthmically distributed parameter points in $\mu_1$ and $\lambda$. The effectivity $\eta^{J, \ast}_{\textrm{ave}}$ 
is the average over $\Xi_{test}$ of $\Delta^{J,\ast}_{N}(\mu)/|J^\ast(\mu)-J^\ast_N(\mu)|$. 
The average effectivities regarding the optimal state and control are defined analogously. 
The defined quantities are presented in Table \ref{tab:Eta_BasisDeltaU_UC} as a function of $N$ for the unconstrained case, where $J^\ast_{\max} = 6.31\cdot 10^{-2}$,  $u^\ast_{\max} = 73.5$, and $y^\ast_{\max} = 3.03$.  We again observe exponential convergence of the state, control, and cost functional errors and bounds but also a considerable overestimation of the actual error. Whereas the effectivities for the cost functional remain fairly constant over $N$, the effectivities for the state and control deteriorate significantly with increasing basis size $N$.

\begin{table}[t]
\begin{center}
\pgfplotstableread{tables/data_etaU_welding_DU_uc_lambda6_4_K75_MS5_fine_update.dat}\dataA
\pgfplotstableread{tables/data_etaJ_welding_DU_uc_lambda6_4_K75_MS5_fine_update.dat}\dataB
\pgfplotstableread{tables/data_etaY_welding_DU_uc_lambda6_4_K75_MS5_fine_update.dat}\dataC

\pgfplotstablecreatecol[copy column from table={\dataB}{[index] 1}] {DeltaJ-rel-max} {\dataA}
\pgfplotstablecreatecol[copy column from table={\dataB}{[index] 2}] {errJ-rel-max} {\dataA}
\pgfplotstablecreatecol[copy column from table={\dataB}{[index] 3}] {etaJ-ave} {\dataA}
\pgfplotstablecreatecol[copy column from table={\dataB}{[index] 4}] {etaJ-min} {\dataA}
\pgfplotstablecreatecol[copy column from table={\dataB}{[index] 5}] {etaJ-max} {\dataA}

\pgfplotstablecreatecol[copy column from table={\dataC}{[index] 1}] {DeltaY-rel-max} {\dataA}
\pgfplotstablecreatecol[copy column from table={\dataC}{[index] 2}] {errY-rel-max} {\dataA}
\pgfplotstablecreatecol[copy column from table={\dataC}{[index] 3}] {etaY-ave} {\dataA}
\pgfplotstablecreatecol[copy column from table={\dataC}{[index] 4}] {etaY-min} {\dataA}
\pgfplotstablecreatecol[copy column from table={\dataC}{[index] 5}] {etaY-max} {\dataA}


\pgfplotstabletypeset[font=\footnotesize, 
		every head row/.style={
			before row={
				\multicolumn{1}{c}{} & \multicolumn{3}{c}{State}
				& \multicolumn{3}{c}{Control}& 
				\multicolumn{3}{c}{Cost functional}
				\\
				\hline  & & &   & & &  & & &\\[-7pt]
			},
			after row= \hline & & &  & & &  & & &\\[-7pt],
		},
			columns={N, errY-rel-max, DeltaY-rel-max, etaY-ave, errU-rel-max, DeltaU-rel-max, etaU-ave,  errJ-rel-max, DeltaJ-rel-max, etaJ-ave},
	columns/N/.style={fixed, column type/.add={}{|},column name=$N$},
		columns/DeltaY-rel-max/.style={column name = $\Delta^{y,\ast}_{N,\textrm{rel,max}}$, sci , sci E, sci zerofill},
	columns/errY-rel-max/.style={column name = $\epsilon^{y,\ast}_{N,\textrm{rel,max}}$, sci, sci E,	sci zerofill},
	columns/etaY-ave/.style={column name = $\eta^{y,\ast}_{\textrm{ave}}$, sci, sci E,
	sci zerofill,
	column type/.add={}{|}},
	columns/DeltaU-rel-max/.style={column name = $\Delta^{u,\ast}_{N,\textrm{rel,max}}$, sci , sci E, sci zerofill},
	columns/errU-rel-max/.style={column name = $\epsilon^{u,\ast}_{N,\textrm{rel,max}}$, sci, sci E,	sci zerofill},
	columns/etaU-ave/.style={column name = $\eta^{u,\ast}_{\textrm{ave}}$, sci, sci E,
	sci zerofill,
	column type/.add={}{|}},
	columns/DeltaJ-rel-max/.style={column name = $\Delta^{J,\ast}_{N,\textrm{rel,max}}$, sci , sci E, sci zerofill},
	columns/errJ-rel-max/.style={column name = $\epsilon^{J,\ast}_{N,\textrm{rel,max}}$, sci, sci E,	sci zerofill},
	columns/etaJ-ave/.style={column name = $\eta^{J,\ast}_{\textrm{ave}}$, sci, sci E,
	sci zerofill},
		]{\dataA}
\end{center}
\caption{Maximum relative errors and error bounds with corresponding average effectivities over $\Xi_{test}$ as a function of $N$.}
\label{tab:Eta_BasisDeltaU_UC}
\end{table}

\subsubsection{Minimal Stabilizing Horizon}
\label{sec:NumExp_WP_MinH}

We test our proposed RB-MPC algorithm for the start-up of the welding process. We consider the initial condition $y_0(x,t=0) = 0$, choose the time interval $I=[0,3]$ corresponding to $30$ iterations of the RB-MPC algorithm and set the largest possible prediction horizon to $K_{max}=75$. Note that we apply the first $5$ controls in each iteration of the RB-MPC algorithm since $\taucon = 5 \tau = 0.1$. We present results for the constrained and unconstrained case for various representative parameter values $\mu = (\mu_1, \lambda)$ with $\mu_1 = 0.75,  1.75$, and $\lambda = 2\cdot10^{-6}, 5\cdot10^{-6},8\cdot10^{-6}, 5\cdot10^{-5}$. 


\medskip
 
\paragraph{\bf Unconstrained Case.}

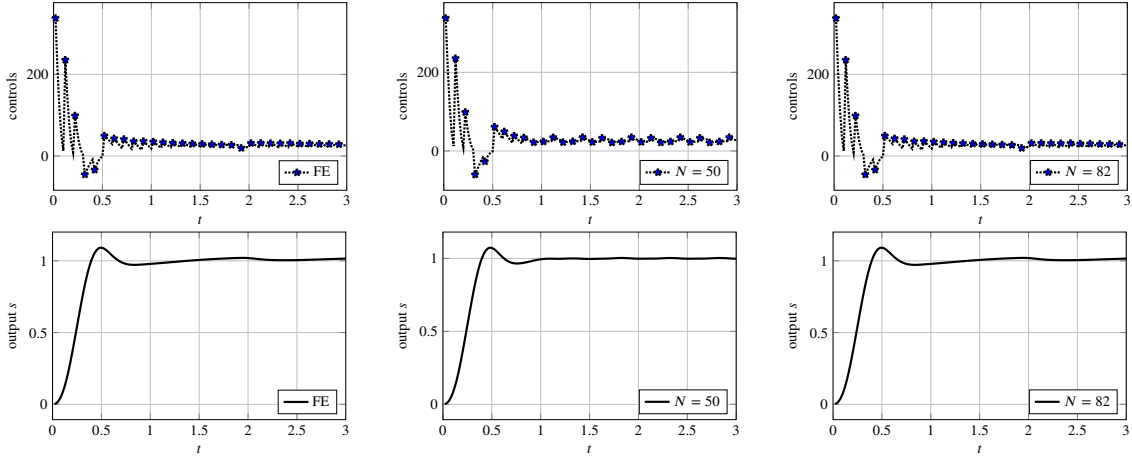
\begin{figure}
\centering
\begin{minipage}[b]{0.3 \textwidth}
\centering
	\begin{tikzpicture}[scale=0.55]
	\begin{axis}[xlabel=$t$, ylabel=controls, xmin=0, xmax=3, 
	width=1.7\textwidth,
    height=1.2\textwidth,
	xmajorgrids,
	ymajorgrids,
	yminorgrids,
	yminorticks=true, 
	axis background/.style={fill=white},
		legend pos=south east]
		\addplot+[densely dotted, color=black, line width=1.5pt, mark=*, mark repeat=5, mark phase = 1]
			table[x=t,y=u-FE] {figures/plotu_welding_DU_uc_lambda6_4_K75_MS5_fine_minh_lambda2e-06_kappa1.75_I30_omega0.0_simple.dat};
			\addlegendentry{FE}
	\end{axis}
	\end{tikzpicture}
\end{minipage}
\begin{minipage}[b]{0.3 \textwidth}
\centering
	\begin{tikzpicture}[scale=0.55]
	\begin{axis}[xlabel=$t$, ylabel=controls, xmin=0, xmax=3, 
	width=1.7\textwidth,
    height=1.2\textwidth,
	xmajorgrids,
	ymajorgrids,
	yminorgrids,
	yminorticks=true, 
	axis background/.style={fill=white},
		legend pos= south east]
		\addplot+[densely dotted, color=black, line width=1.5pt, mark=*, mark repeat=5, mark phase = 1]
			table[x=t,y=uN-2] {figures/plotu_welding_DU_uc_lambda6_4_K75_MS5_fine_minh_lambda2e-06_kappa1.75_I30_omega0.0_simple.dat};
			\addlegendentry{$N=50$}
	\end{axis}
	\end{tikzpicture}
\end{minipage}
\begin{minipage}[b]{0.3 \textwidth}
\centering
	\begin{tikzpicture}[scale=0.55]
	\begin{axis}[xlabel=$t$, ylabel=controls, xmin=0, xmax=3, 
	width=1.7\textwidth,
    height=1.2\textwidth,
	xmajorgrids,
	ymajorgrids,
	yminorgrids,
	yminorticks=true, 
	axis background/.style={fill=white},
		legend pos= south east]
		\addplot+[densely dotted, color=black, line width=1.5pt, mark=*, mark repeat=5, mark phase = 1]
			table[x=t,y=uN-4] {figures/plotu_welding_DU_uc_lambda6_4_K75_MS5_fine_minh_lambda2e-06_kappa1.75_I30_omega0.0_simple.dat};
			\addlegendentry{$N=82$}
	\end{axis}
	\end{tikzpicture}
\end{minipage}
%
\begin{minipage}[b]{0.3 \textwidth}
\centering
	\begin{tikzpicture}[scale=0.55]
	\begin{axis}[xlabel=$t$, ylabel=output $s$, xmin=0, xmax=3, 
	width=1.7\textwidth,
    height=1.2\textwidth,
	xmajorgrids,
	ymajorgrids,
	yminorgrids,
	yminorticks=true, 
	axis background/.style={fill=white},
		legend pos=south east]
		\addplot+[color=black, line width=1.5pt, mark=none]
			table[x=t,y=s-FE] {figures/plotu_welding_DU_uc_lambda6_4_K75_MS5_fine_minh_lambda2e-06_kappa1.75_I30_omega0.0_simple.dat};
			\addlegendentry{FE}
	\end{axis}
	\end{tikzpicture}
\end{minipage}
\begin{minipage}[b]{0.3 \textwidth}
\centering
	\begin{tikzpicture}[scale=0.55]
	\begin{axis}[xlabel=$t$, ylabel=output $s$, xmin=0, xmax=3, 
	width=1.7\textwidth,
    height=1.2\textwidth,
	xmajorgrids,
	ymajorgrids,
	yminorgrids,
	yminorticks=true, 
	axis background/.style={fill=white},
		legend pos= south east]
		\addplot+[color=black, line width=1.5pt,mark=none]
			table[x=t,y=s-2] {figures/plotu_welding_DU_uc_lambda6_4_K75_MS5_fine_minh_lambda2e-06_kappa1.75_I30_omega0.0_simple.dat};
			\addlegendentry{$N=50$}
	\end{axis}
	\end{tikzpicture}
\end{minipage}
\begin{minipage}[b]{0.3 \textwidth}
\centering
	\begin{tikzpicture}[scale=0.55]
	\begin{axis}[xlabel=$t$, ylabel=output $s$, xmin=0, xmax=3, 
	width=1.7\textwidth,
    height=1.2\textwidth,
	xmajorgrids,
	ymajorgrids,
	yminorgrids,
	yminorticks=true, 
	axis background/.style={fill=white},
		legend pos= south east]
		\addplot+[color=black, line width=1.5pt,mark=none]
			table[x=t,y=s-4] {figures/plotu_welding_DU_uc_lambda6_4_K75_MS5_fine_minh_lambda2e-06_kappa1.75_I30_omega0.0_simple.dat};
			\addlegendentry{$N=82$}
	\end{axis}
	\end{tikzpicture}
\end{minipage}
\caption{Unconstrained case: controls and outputs for $\mu_1=1.75$, $\lambda= 2\cdot10^{-6}$, FE-MPC on the left, RB-MPC with $N = 50$ (middle) and $N = 82$ (right).}
\label{fig:WP_MinH_uc_mu1_1.75}
\end{figure}

We first consider the unconstrained case and focus on the minimal stabilizing horizons. In Table~\ref{tab:WP_MinH_uc}, we present the average stabilizing horizons $K_{\textrm{ave}}$ over the RB-MPC iterations of the feedback control $\kappa_N$ for basis dimensions $N=34, \, 50, \, 66,$ and $82$. For comparison, we also present the corresponding result for the feedback control $\kappa$ resulting from the classical MPC approach in the last row. In Figure~\ref{fig:WP_MinH_uc_mu1_1.75}, we plot the control and output as a function of time for $(\mu_1, \lambda) = (1.75, 2\cdot10^{-6})$. The FE-MPC result is shown on the left, the RB-MPC results for $N = 50$ and $N = 82$ is shown in the middle and on the right, respectively. The dots in the control plots have a distance of $\taucon$ and represent the starting points of the MPC loops.

\begin{table}[t]
\begin{center}
\pgfplotstableread{tables/data_welding_DU_uc_lambda6_4_K75_MS5_fine_minh_lambda2e-06_I30_omega0.0_simple.dat}\dataA
\pgfplotstableread{tables/data_welding_DU_uc_lambda6_4_K75_MS5_fine_minh_lambda5e-06_I30_omega0.0_simple.dat}\dataB
\pgfplotstableread{tables/data_welding_DU_uc_lambda6_4_K75_MS5_fine_minh_lambda8e-06_I30_omega0.0_simple.dat}\dataC
\pgfplotstableread{tables/data_welding_DU_uc_lambda6_4_K75_MS5_fine_minh_lambda5e-05_I30_omega0.0_simple.dat}\dataD

\pgfplotstablecreatecol[copy column from table={\dataB}{[index] 1}] {56ave1} {\dataA}
\pgfplotstablecreatecol[copy column from table={\dataB}{[index] 4}] {56ave2} {\dataA}

\pgfplotstablecreatecol[copy column from table={\dataC}{[index] 1}] {86ave1} {\dataA}
\pgfplotstablecreatecol[copy column from table={\dataC}{[index] 4}] {86ave2} {\dataA}

\pgfplotstablecreatecol[copy column from table={\dataD}{[index] 1}] {55ave1} {\dataA}
\pgfplotstablecreatecol[copy column from table={\dataD}{[index] 4}] {55ave2} {\dataA}

	\pgfplotstabletypeset[font=\footnotesize, 
		every row no 4/.style={before row=\hline \hline},
		every head row/.style={
			before row={
				 $\lambda$  & \multicolumn{2}{c}{$2\cdot10^{-6}$ } & \multicolumn{2}{c}{$5\cdot10^{-6}$ }  & \multicolumn{2}{c}{$8\cdot10^{-6}$ } & \multicolumn{2}{c}{$5\cdot10^{-5}$ }\\
				\hline & &  & &  & &  & &  \\[-7pt]
				$\mu_1$ & 0.75& 1.75  & 0.75& 1.75 & 0.75& 1.75  & 0.75& 1.75\\
				\hline & &   & &  & &  & &   \\[-7pt]
				$N$ & \multicolumn{2}{c|}{$K_{\textrm{ave}}$} & \multicolumn{2}{c|}{$K_{\textrm{ave}}$} & \multicolumn{2}{c|}{$K_{\textrm{ave}}$} & \multicolumn{2}{c}{$K_{\textrm{ave}}$} \\
				\hline & &  & &  & &  & &    \\[-7pt]
			},
			output empty row
		},
		columns={N, K-ave1, K-ave2, 56ave1, 56ave2, 86ave1, 86ave2, 55ave1, 55ave2},
        columns/N/.append style={string replace={-1}{}},
		columns/K-ave2/.append style={column type={c|}},
		columns/56ave2/.append style={column type={c|}},
		columns/86ave2/.append style={column type={c|}},
        empty cells with={$\kappa$}
        ]{\dataA}
\end{center} 
\caption{Unconstrained case: average minimal stabilizing horizons over $[0,3]$.}
\label{tab:WP_MinH_uc}
\end{table}

\begin{figure}[t]
\centering
\begin{minipage}[b]{0.9 \textwidth}
	\centering
	\begin{tikzpicture}[scale=0.55]
	\begin{axis}[xlabel=$t$, ylabel=$K$, xmin=0, xmax=3, ymin=0,	ymax=80,
	width=1.5\textwidth,
    height=0.5\textwidth,
		legend pos=outer north east]
		\addplot+[only marks, color=red,line width=1.5pt,mark=+, mark options={fill=red,scale=4.5}]
			table[x=time, y=K-FE] {figures/data_Kmax_welding_DU_uc_lambda6_4_K75_MS5_fine_minh_lambda2e-06_kappa1.75_I30_omega0.0_simple.dat};
			\addlegendentry{$\textrm{FE}$}
\addplot+[only marks, color=orange,line width=1.5pt,mark=square*, mark options={fill=orange,scale=3.5}]
			table[x=time, y=K-1] {figures/data_Kmax_welding_DU_uc_lambda6_4_K75_MS5_fine_minh_lambda2e-06_kappa1.75_I30_omega0.0_simple.dat};
			\addlegendentry{$N=34$}
\addplot+[only marks, color=black,line width=1.5pt,mark=diamond*, mark options={fill=black!30,scale=4.5}]
			table[x=time, y=K-2] {figures/data_Kmax_welding_DU_uc_lambda6_4_K75_MS5_fine_minh_lambda2e-06_kappa1.75_I30_omega0.0_simple.dat};
			\addlegendentry{$N=50$}
\addplot+[only marks, color=purple,line width=1.8pt,mark=o, mark options={fill=purple!30,scale=3}]
			table[x=time, y=K-3] {figures/data_Kmax_welding_DU_uc_lambda6_4_K75_MS5_fine_minh_lambda2e-06_kappa1.75_I30_omega0.0_simple.dat};
			\addlegendentry{$N=66$}
\addplot+[only marks, color=blue,line width=1.5pt,mark=x, mark options={fill=blue,scale=4.5}]
			table[x=time, y=K-4] {figures/data_Kmax_welding_DU_uc_lambda6_4_K75_MS5_fine_minh_lambda2e-06_kappa1.75_I30_omega0.0_simple.dat};
			\addlegendentry{$N=82$}
	\end{axis}
	\end{tikzpicture}
\end{minipage}
\caption{Computed minimal stabilizing horizons using FE-MPC and RB-MPC for $\mu_1=1.75$, $ \lambda = 2\cdot 10^{-6}$.}
\label{fig:allK_min_stab_horizon}
\end{figure}
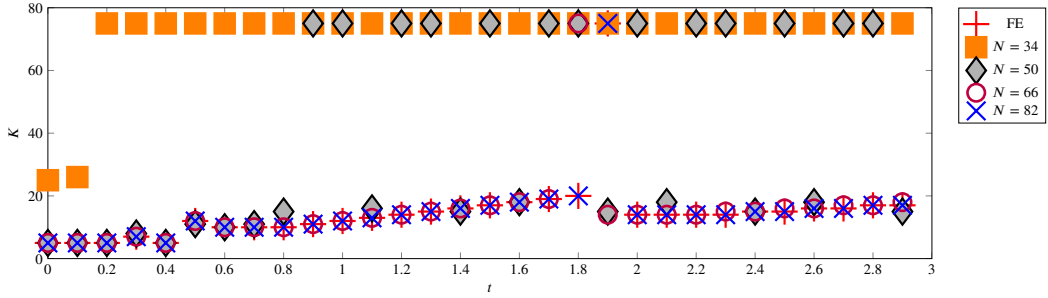

We again observe that the average stabilizing horizons decrease with increasing $N$ and are almost identical to the reference value for $N = 66$ and $82$. Recalling that  $K_{\textrm{max}} = 75$ is the largest possible prediction horizon, we conclude that the RB-MPC feedback control $\kappa_N$ obviously did not meet the desired stability criterion for $\mu_1=0.75$, $\lambda= 2\cdot10^{-6}$, and $N = 34$. To further analyse the stability of the proposed RB-MPC approach -- in view of our discussion in Section~\ref{sec:RB-MPC_Alg} -- we present the  maximum prediction horizon computed along each trajectory in Table~\ref{tab:WP_MinH_uc_Kmax}. We observe that the maximal horizons also decrease with increasing reduced basis size and that the values coincide with the reference value at the latest for $N=82$.  The only exception to this is the control problem for $\mu_1=1.75$ and $\lambda= 2\cdot10^{-6}$. Here, the maximal prediction horizon reaches $K_{\textrm{max}}$ in all cases, including the FE reference. Since this parameter value seems to be the most critical, we show all computed prediction horizons in Figure \ref{fig:allK_min_stab_horizon} over time. The instances where $K_{\textrm{max}}$ is reached decreases with increasing $N$, and  $K_{\textrm{max}}$ is in fact only reached once for $N = 66$ and $N = 82$. For the FE reference $K_{\textrm{max}}$ is also reached once at the same timestep as for $N = 82$ (for $N = 66$ one timestep earlier). These outliers are not due to a lack of accuracy of the RB approximation or the {\it a posteriori} error bounds, but correspond to instances where the actual state is close to the desired state and numerical errors dominate.

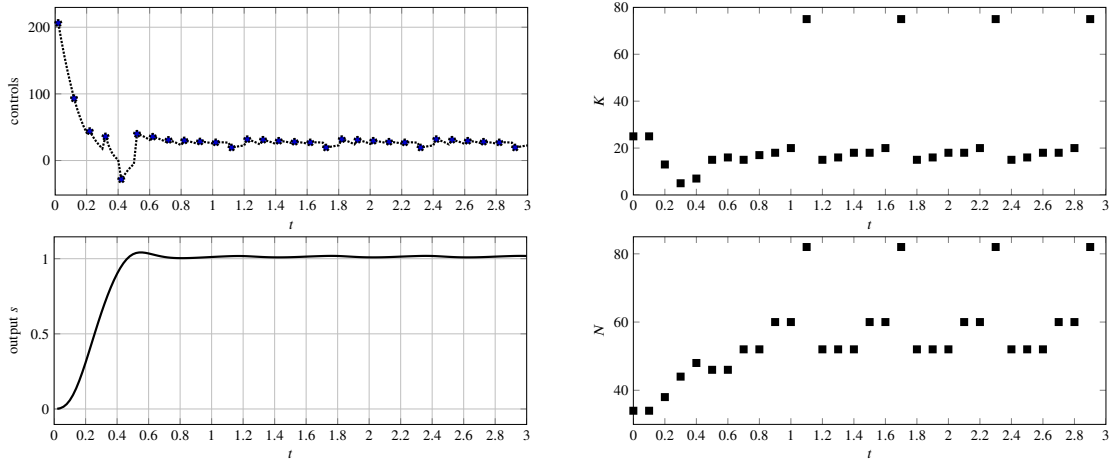
\begin{figure}
\centering
\begin{minipage}[b]{0.45\textwidth}
\centering
	\begin{tikzpicture}[scale=0.55]
	\begin{axis}[xlabel=$t$, ylabel=controls, xmin=0, xmax=3, 
	width=1.7\textwidth,
    height=0.8\textwidth,
	xmajorgrids,
	ymajorgrids,
	yminorgrids,
	yminorticks=true, 
	axis background/.style={fill=white},
		legend pos=south east]
		\addplot+[densely dotted, color=black, line width=1.5pt, mark=*, mark repeat=5, mark phase = 1]
			table[x=t,y=u-varN] {figures/plotu_varN_welding_DU_uc_lambda6_4_K75_MS5_fine_minh_lambda2e-06_kappa1.75_I30_omega0.0_Nmin34_Nmax82_simple.dat};
	\end{axis}
	\end{tikzpicture}
\end{minipage}
\begin{minipage}[b]{0.45 \textwidth}
	\centering
	\begin{tikzpicture}[scale=0.55]
	\begin{axis}[xlabel=$t$, ylabel=$K$,xmin=0, xmax=3, ymin=0,	ymax=80,
	width=1.7\textwidth,
    height=0.8\textwidth,
		legend pos=outer north east]
		\addplot+[only marks, color=black,line width=1.5pt,mark=square*, mark options={fill=black,scale=1}]
			table[x=t, y=K] {figures/data_varN_welding_DU_uc_lambda6_4_K75_MS5_fine_minh_lambda2e-06_I30_omega0.0_Nmin34_Nmax82_simple.dat};
	\end{axis}
	\end{tikzpicture}
\end{minipage}
%
\begin{minipage}[b]{0.45 \textwidth}
\centering
	\begin{tikzpicture}[scale=0.55]
	\begin{axis}[xlabel=$t$, ylabel=output $s$, xmin=0, xmax=3, 
	width=1.7\textwidth,
    height=0.8\textwidth,
	xmajorgrids,
	ymajorgrids,
	yminorgrids,
	yminorticks=true, 
	axis background/.style={fill=white},
		legend pos=south east]
		\addplot+[color=black, line width=1.5pt, mark=none]
			table[x=t,y=s-varN] {figures/plotu_varN_welding_DU_uc_lambda6_4_K75_MS5_fine_minh_lambda2e-06_kappa1.75_I30_omega0.0_Nmin34_Nmax82_simple.dat};
	\end{axis}
	\end{tikzpicture}
\end{minipage}
\begin{minipage}[b]{0.45 \textwidth}
	\centering
	\begin{tikzpicture}[scale=0.55]
	\begin{axis}[xlabel=$t$, ylabel=$N$, xmin=0, xmax=3, ymin=30,	ymax=85,
	width=1.7\textwidth,
    height=0.8\textwidth,
		legend pos=outer north east]
		\addplot+[only marks, color=black,line width=1.5pt,mark=square*, mark options={fill=black,scale=1}]
			table[x=t, y=N] {figures/data_varN_welding_DU_uc_lambda6_4_K75_MS5_fine_minh_lambda2e-06_I30_omega0.0_Nmin34_Nmax82_simple.dat};
	\end{axis}
	\end{tikzpicture}
\end{minipage}
\caption{Modified Algorithm~\ref{alg:Adaptive_RB_MPC} for parameter value $\mu_1=1.75$, $\lambda= 2\cdot10^{-6}$: control and output over time (left column), computed minimal stabilizing horizon $K$ and corresponding reduced basis size (right column) with $N_{\min} = 34$ and $N_{\max} = 82$.}
\label{fig:WP_MinH_uc_mu1_1.75_VarN_Nmin34}
\end{figure}

Obviously, the approximation quality of the reduced model with basis size $N=34$ is too poor to produce a stabilizing RB controller. As discussed in Section~\ref{sec:RB-MPC_Alg}, we thus initialize the RB size with $N_{\min}=34$ and test the modified Algorithm~\ref{alg:Adaptive_RB_MPC}, i.e., where we also adaptively increase $N$. We do not consider the limit $N \to {\cal N}$ but choose a maximum RB size of $N=82$. We show the RB control, output, minimal stabilizing horizon, and reduced basis size over time in Figure \ref{fig:WP_MinH_uc_mu1_1.75_VarN_Nmin34}. We observe that the required RB size increases as the output approaches the desired output and that the system seems to reach a periodic behaviour. There are four instances where the algorithm reaches the maximum prediction horizon despite having choosen the maximum RB size. Since the reason is not a lack of accuracy of the reduced model but numerical errors, further increasing $N$ would not alleviate this issue. Based on our discussion in Section~\ref{sec:RB-MPC_Alg}, we believe that a further improvement could be achieved by introducing an $\varepsilon$-neighborhood around the desired state and extending the analysis presented here to ``practical optimality." 

As opposed to the previous example, we note that the average minimal stabilizing horizon for this numerical example is higher for smaller $\lambda$. This is due to the fact that a smaller $\lambda$ results in a higher control effort and the desired output is thus reached sooner for smaller values of $\lambda$ than for larger values of $\lambda$. Since the stabilizing horizon is larger in this case as just mentioned, we also obtain a larger horizon on average.

\begin{table}[t]
\begin{center}
\pgfplotstableread{tables/data_welding_DU_uc_lambda6_4_K75_MS5_fine_minh_lambda2e-06_I30_omega0.0_simple.dat}\dataA
\pgfplotstableread{tables/data_welding_DU_uc_lambda6_4_K75_MS5_fine_minh_lambda5e-06_I30_omega0.0_simple.dat}\dataB
\pgfplotstableread{tables/data_welding_DU_uc_lambda6_4_K75_MS5_fine_minh_lambda8e-06_I30_omega0.0_simple.dat}\dataC
\pgfplotstableread{tables/data_welding_DU_uc_lambda6_4_K75_MS5_fine_minh_lambda5e-05_I30_omega0.0_simple.dat}\dataD

\pgfplotstablecreatecol[copy column from table={\dataB}{[index] 3}] {56max1} {\dataA}
\pgfplotstablecreatecol[copy column from table={\dataB}{[index] 6}] {56max2} {\dataA}

\pgfplotstablecreatecol[copy column from table={\dataC}{[index] 3}] {86max1} {\dataA}
\pgfplotstablecreatecol[copy column from table={\dataC}{[index] 6}] {86max2} {\dataA}

\pgfplotstablecreatecol[copy column from table={\dataD}{[index] 3}] {55max1} {\dataA}
\pgfplotstablecreatecol[copy column from table={\dataD}{[index] 6}] {55max2} {\dataA}

	\pgfplotstabletypeset[font=\footnotesize, 
		every row no 4/.style={before row=\hline \hline},
		every head row/.style={
			before row={
				 $\lambda$  & \multicolumn{2}{c}{$2\cdot10^{-6}$ } & \multicolumn{2}{c}{$5\cdot10^{-6}$ }  & \multicolumn{2}{c}{$8\cdot10^{-6}$ } & \multicolumn{2}{c}{$5\cdot10^{-5}$ }\\
				\hline & &  & &  & &  & &  \\[-7pt]
				$\mu_1$ & 0.75& 1.75  & 0.75& 1.75 & 0.75& 1.75  & 0.75& 1.75\\
				\hline & &   & &  & &  & &   \\[-7pt]
				$N$ & \multicolumn{2}{c|}{$\bar{K}$} & \multicolumn{2}{c|}{$\bar{K}$} & \multicolumn{2}{c|}{$\bar{K}$} & \multicolumn{2}{c}{$\bar{K}$} \\
				\hline & &  & &  & &  & &    \\[-7pt]
			},
			output empty row
		},
		columns={N, K-max1, K-max2, 56max1, 56max2, 86max1, 86max2, 55max1, 55max2},
        columns/N/.append style={string replace={-1}{}},
		columns/K-max2/.append style={column type={c|}},
		columns/56max2/.append style={column type={c|}},
		columns/86max2/.append style={column type={c|}},
        empty cells with={$\kappa$}
        ]{\dataA}
\end{center} 
\caption{Unconstrained case: maximal minimal stabilizing horizon $\bar{K}(\le K_{\max})$ over $[0,3]$.}
\label{tab:WP_MinH_uc_Kmax}
\end{table}

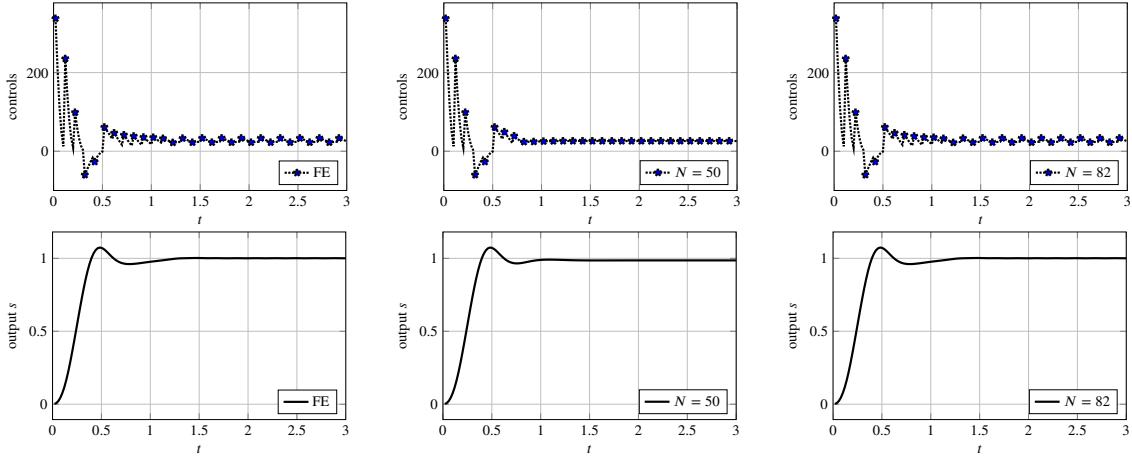
\begin{figure}
\centering
\begin{minipage}[b]{0.3 \textwidth}
\centering
	\begin{tikzpicture}[scale=0.55]
	\begin{axis}[xlabel=$t$, ylabel=controls, xmin=0, xmax=3, 
	width=1.7\textwidth,
    height=1.2\textwidth,
	xmajorgrids,
	ymajorgrids,
	yminorgrids,
	yminorticks=true, 
	axis background/.style={fill=white},
		legend pos=south east]
		\addplot+[densely dotted, color=black, line width=1.5pt, mark=*, mark repeat=5, mark phase = 1]
			table[x=t,y=u-FE] {figures/plotu_welding_DU_uc_lambda6_4_K75_MS5_fine_minh_lambda2e-06_kappa1.75_I30_omega0.2_simple.dat};
			\addlegendentry{FE}
	\end{axis}
	\end{tikzpicture}
\end{minipage}
\begin{minipage}[b]{0.3 \textwidth}
\centering
	\begin{tikzpicture}[scale=0.55]
	\begin{axis}[xlabel=$t$, ylabel=controls, xmin=0, xmax=3, 
	width=1.7\textwidth,
    height=1.2\textwidth,
	xmajorgrids,
	ymajorgrids,
	yminorgrids,
	yminorticks=true, 
	axis background/.style={fill=white},
		legend pos= south east]
		\addplot+[densely dotted, color=black, line width=1.5pt, mark=*, mark repeat=5, mark phase = 1]
			table[x=t,y=uN-2] {figures/plotu_welding_DU_uc_lambda6_4_K75_MS5_fine_minh_lambda2e-06_kappa1.75_I30_omega0.2_simple.dat};
			\addlegendentry{$N=50$}
	\end{axis}
	\end{tikzpicture}
\end{minipage}
\begin{minipage}[b]{0.3 \textwidth}
\centering
	\begin{tikzpicture}[scale=0.55]
	\begin{axis}[xlabel=$t$, ylabel=controls, xmin=0, xmax=3, 
	width=1.7\textwidth,
    height=1.2\textwidth,
	xmajorgrids,
	ymajorgrids,
	yminorgrids,
	yminorticks=true, 
	axis background/.style={fill=white},
		legend pos= south east]
		\addplot+[densely dotted, color=black, line width=1.5pt, mark=*, mark repeat=5, mark phase = 1]
			table[x=t,y=uN-4] {figures/plotu_welding_DU_uc_lambda6_4_K75_MS5_fine_minh_lambda2e-06_kappa1.75_I30_omega0.2_simple.dat};
			\addlegendentry{$N=82$}
	\end{axis}
	\end{tikzpicture}
\end{minipage}
%
\begin{minipage}[b]{0.3 \textwidth}
\centering
	\begin{tikzpicture}[scale=0.55]
	\begin{axis}[xlabel=$t$, ylabel=output $s$, xmin=0, xmax=3, 
	width=1.7\textwidth,
    height=1.2\textwidth,
	xmajorgrids,
	ymajorgrids,
	yminorgrids,
	yminorticks=true, 
	axis background/.style={fill=white},
		legend pos=south east]
		\addplot+[color=black, line width=1.5pt, mark=none]
			table[x=t,y=s-FE] {figures/plotu_welding_DU_uc_lambda6_4_K75_MS5_fine_minh_lambda2e-06_kappa1.75_I30_omega0.2_simple.dat};
			\addlegendentry{FE}
	\end{axis}
	\end{tikzpicture}
\end{minipage}
\begin{minipage}[b]{0.3 \textwidth}
\centering
	\begin{tikzpicture}[scale=0.55]
	\begin{axis}[xlabel=$t$, ylabel=output $s$, xmin=0, xmax=3, 
	width=1.7\textwidth,
    height=1.2\textwidth,
	xmajorgrids,
	ymajorgrids,
	yminorgrids,
	yminorticks=true, 
	axis background/.style={fill=white},
		legend pos= south east]
		\addplot+[color=black, line width=1.5pt,mark=none]
			table[x=t,y=s-2] {figures/plotu_welding_DU_uc_lambda6_4_K75_MS5_fine_minh_lambda2e-06_kappa1.75_I30_omega0.2_simple.dat};
			\addlegendentry{$N=50$}
	\end{axis}
	\end{tikzpicture}
\end{minipage}
\begin{minipage}[b]{0.3 \textwidth}
\centering
	\begin{tikzpicture}[scale=0.55]
	\begin{axis}[xlabel=$t$, ylabel=output $s$, xmin=0, xmax=3, 
	width=1.7\textwidth,
    height=1.2\textwidth,
	xmajorgrids,
	ymajorgrids,
	yminorgrids,
	yminorticks=true, 
	axis background/.style={fill=white},
		legend pos= south east]
		\addplot+[color=black, line width=1.5pt,mark=none]
			table[x=t,y=s-4] {figures/plotu_welding_DU_uc_lambda6_4_K75_MS5_fine_minh_lambda2e-06_kappa1.75_I30_omega0.2_simple.dat};
			\addlegendentry{$N=82$}
	\end{axis}
	\end{tikzpicture}
\end{minipage}
\caption{Unconstrained case with $\omega_{\min}=0.2$: controls and outputs for $\mu_1=1.75$, $\lambda= 2\cdot10^{-6}$, FE-MPC on the left, RB-MPC with $N = 50$ (middle) and $N = 82$ (right).}
\label{fig:WP_MinH_uc_omega_mu1_1.75}
\end{figure}

So far we accepted the RB-MPC feedback control $\kappa_N$ in Algorithm~\ref{alg:Adaptive_RB_MPC} as soon as the suboptimality degree $\tilde{\omega}_{N,K}(\mu)$ became positive. However, we can increase the stability margin of the feedback law by incorporating a threshold $\omega_{\min}>0$ in the algorithm, i.e., we modify the stability condition and accept $K$ as soon as $\tilde{\omega}_{N,K}(\mu)>\omega_{\min}>0$ holds. Although this approach will increase the minimal stabilizing horizon, it may improve the smoothness of the control and the control performance. In Table~\ref{tab:WP_MinH_uc_omega} and Figure~\ref{fig:WP_MinH_uc_omega_mu1_1.75} we show the results corresponding to Table~\ref{tab:WP_MinH_uc} and Figure~\ref{fig:WP_MinH_uc_mu1_1.75} for a threshold $\omega_{\min}=0.2$. As expected, we observe an increase in the minimal stabilizing horizon, but at the same time a smoother control and hence output as a function of time. This effect will be more pronounced in the control constrained case discussed next.

\begin{table}[t]
\begin{center}
\pgfplotstableread{tables/data_welding_DU_uc_lambda6_4_K75_MS5_fine_minh_lambda2e-06_I30_omega0.2_simple.dat}\dataA
\pgfplotstableread{tables/data_welding_DU_uc_lambda6_4_K75_MS5_fine_minh_lambda5e-06_I30_omega0.2_simple.dat}\dataB
\pgfplotstableread{tables/data_welding_DU_uc_lambda6_4_K75_MS5_fine_minh_lambda8e-06_I30_omega0.2_simple.dat}\dataC
\pgfplotstableread{tables/data_welding_DU_uc_lambda6_4_K75_MS5_fine_minh_lambda5e-05_I30_omega0.2_simple.dat}\dataD

\pgfplotstablecreatecol[copy column from table={\dataB}{[index] 1}] {56ave1} {\dataA}
\pgfplotstablecreatecol[copy column from table={\dataB}{[index] 4}] {56ave2} {\dataA}

\pgfplotstablecreatecol[copy column from table={\dataC}{[index] 1}] {86ave1} {\dataA}
\pgfplotstablecreatecol[copy column from table={\dataC}{[index] 4}] {86ave2} {\dataA}

\pgfplotstablecreatecol[copy column from table={\dataD}{[index] 1}] {55ave1} {\dataA}
\pgfplotstablecreatecol[copy column from table={\dataD}{[index] 4}] {55ave2} {\dataA}

\pgfplotstabletypeset[font=\footnotesize, 
		every row no 4/.style={before row=\hline \hline},
		every head row/.style={
			before row={
				 $\lambda$  & \multicolumn{2}{c}{$2\cdot10^{-6}$ } & \multicolumn{2}{c}{$5\cdot10^{-6}$ }  & \multicolumn{2}{c}{$8\cdot10^{-6}$ } & \multicolumn{2}{c}{$5\cdot10^{-5}$ }\\
				\hline & &  & &  & &  & &  \\[-7pt]
				$\mu_1$ & 0.75& 1.75  & 0.75& 1.75 & 0.75& 1.75  & 0.75& 1.75\\
				\hline & &   & &  & &  & &   \\[-7pt]
				$N$ & \multicolumn{2}{c|}{$K_{\textrm{ave}}$} & \multicolumn{2}{c|}{$K_{\textrm{ave}}$} & \multicolumn{2}{c|}{$K_{\textrm{ave}}$} & \multicolumn{2}{c}{$K_{\textrm{ave}}$} \\
				\hline & &  & &  & &  & &    \\[-7pt]
			},
			output empty row
		},
		columns={N, K-ave1, K-ave2, 56ave1, 56ave2, 86ave1, 86ave2, 55ave1, 55ave2},
        columns/N/.append style={string replace={-1}{}},
		columns/K-ave2/.append style={column type={c|}},
		columns/56ave2/.append style={column type={c|}},
		columns/86ave2/.append style={column type={c|}},
        empty cells with={$\kappa$}
        ]{\dataA}
\end{center} 
\caption{Unconstrained case: average minimal stabilizing horizons over $[0,3]$ with threshold $\omega_{\min}=0.2$.}
\label{tab:WP_MinH_uc_omega}
\end{table}

\medskip

\paragraph{\bf Constrained Case.}

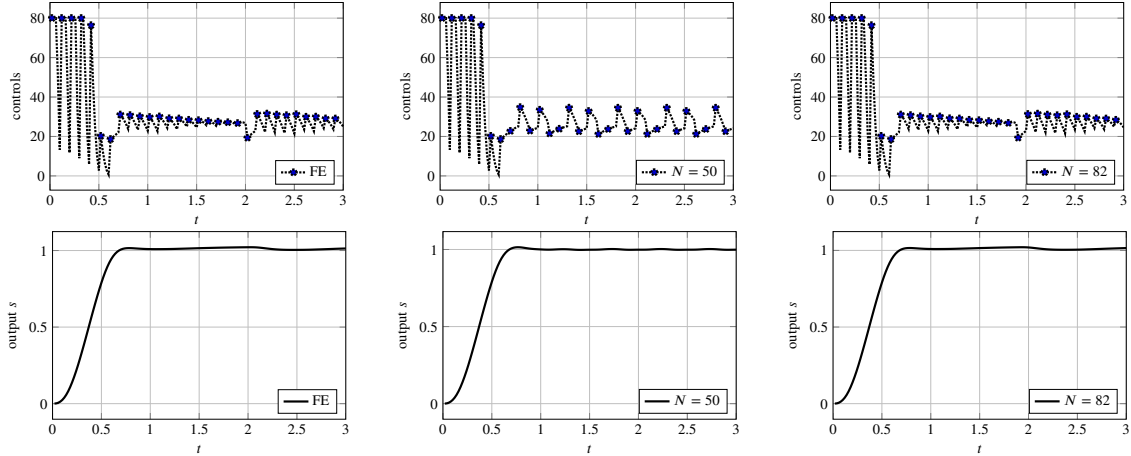
\begin{figure}
\centering
\begin{minipage}[b]{0.3 \textwidth}
\centering
	\begin{tikzpicture}[scale=0.55]
	\begin{axis}[xlabel=$t$, ylabel=controls, xmin=0, xmax=3, 
	width=1.7\textwidth,
    height=1.2\textwidth,
	xmajorgrids,
	ymajorgrids,
	yminorgrids,
	yminorticks=true, 
	axis background/.style={fill=white},
		legend pos=south east]
		\addplot+[densely dotted, color=black, line width=1.5pt, mark=*, mark repeat=5, mark phase = 1]
			table[x=t,y=u-FE] {figures/plotu_welding_DU_uc_lambda6_4_K75_MS5_fine_minh_lambda2e-06_kappa1.75_I30_ua0_ub80_omega0.0.dat};
			\addlegendentry{FE}
	\end{axis}
	\end{tikzpicture}
\end{minipage}
\begin{minipage}[b]{0.3 \textwidth}
\centering
	\begin{tikzpicture}[scale=0.55]
	\begin{axis}[xlabel=$t$, ylabel=controls, xmin=0, xmax=3, 
	width=1.7\textwidth,
    height=1.2\textwidth,
	xmajorgrids,
	ymajorgrids,
	yminorgrids,
	yminorticks=true, 
	axis background/.style={fill=white},
		legend pos= south east]
		\addplot+[densely dotted, color=black, line width=1.5pt, mark=*, mark repeat=5, mark phase = 1]
			table[x=t,y=uN-2] {figures/plotu_welding_DU_uc_lambda6_4_K75_MS5_fine_minh_lambda2e-06_kappa1.75_I30_ua0_ub80_omega0.0.dat};
			\addlegendentry{$N=50$}
	\end{axis}
	\end{tikzpicture}
\end{minipage}
\begin{minipage}[b]{0.3 \textwidth}
\centering
	\begin{tikzpicture}[scale=0.55]
	\begin{axis}[xlabel=$t$, ylabel=controls, xmin=0, xmax=3, 
	width=1.7\textwidth,
    height=1.2\textwidth,
	xmajorgrids,
	ymajorgrids,
	yminorgrids,
	yminorticks=true, 
	axis background/.style={fill=white},
		legend pos= south east]
		\addplot+[densely dotted, color=black, line width=1.5pt, mark=*, mark repeat=5, mark phase = 1]
			table[x=t,y=uN-4] {figures/plotu_welding_DU_uc_lambda6_4_K75_MS5_fine_minh_lambda2e-06_kappa1.75_I30_ua0_ub80_omega0.0.dat};
			\addlegendentry{$N=82$}
	\end{axis}
	\end{tikzpicture}
\end{minipage}
%
\begin{minipage}[b]{0.3 \textwidth}
\centering
	\begin{tikzpicture}[scale=0.55]
	\begin{axis}[xlabel=$t$, ylabel=output $s$, xmin=0, xmax=3, 
	width=1.7\textwidth,
    height=1.2\textwidth,
	xmajorgrids,
	ymajorgrids,
	yminorgrids,
	yminorticks=true, 
	axis background/.style={fill=white},
		legend pos=south east]
		\addplot+[color=black, line width=1.5pt, mark=none]
			table[x=t,y=s-FE] {figures/plotu_welding_DU_uc_lambda6_4_K75_MS5_fine_minh_lambda2e-06_kappa1.75_I30_ua0_ub80_omega0.0.dat};
			\addlegendentry{FE}
	\end{axis}
	\end{tikzpicture}
\end{minipage}
\begin{minipage}[b]{0.3 \textwidth}
\centering
	\begin{tikzpicture}[scale=0.55]
	\begin{axis}[xlabel=$t$, ylabel=output $s$, xmin=0, xmax=3, 
	width=1.7\textwidth,
    height=1.2\textwidth,
	xmajorgrids,
	ymajorgrids,
	yminorgrids,
	yminorticks=true, 
	axis background/.style={fill=white},
		legend pos= south east]
		\addplot+[color=black, line width=1.5pt,mark=none]
			table[x=t,y=s-2] {figures/plotu_welding_DU_uc_lambda6_4_K75_MS5_fine_minh_lambda2e-06_kappa1.75_I30_ua0_ub80_omega0.0.dat};
			\addlegendentry{$N=50$}
	\end{axis}
	\end{tikzpicture}
\end{minipage}
\begin{minipage}[b]{0.3 \textwidth}
\centering
	\begin{tikzpicture}[scale=0.55]
	\begin{axis}[xlabel=$t$, ylabel=output $s$, xmin=0, xmax=3, 
	width=1.7\textwidth,
    height=1.2\textwidth,
	xmajorgrids,
	ymajorgrids,
	yminorgrids,
	yminorticks=true, 
	axis background/.style={fill=white},
		legend pos= south east]
		\addplot+[color=black, line width=1.5pt,mark=none]
			table[x=t,y=s-4] {figures/plotu_welding_DU_uc_lambda6_4_K75_MS5_fine_minh_lambda2e-06_kappa1.75_I30_ua0_ub80_omega0.0.dat};
			\addlegendentry{$N=82$}
	\end{axis}
	\end{tikzpicture}
\end{minipage}
\caption{Constrained case: controls and outputs for $\mu_1=1.75$, $\lambda= 2\cdot10^{-6}$, FE-MPC on the left, RB-MPC with $N = 50$ (middle) and $N = 82$ (right).}
\label{fig:WP_MinH_con_mu1_1.75}
\end{figure}

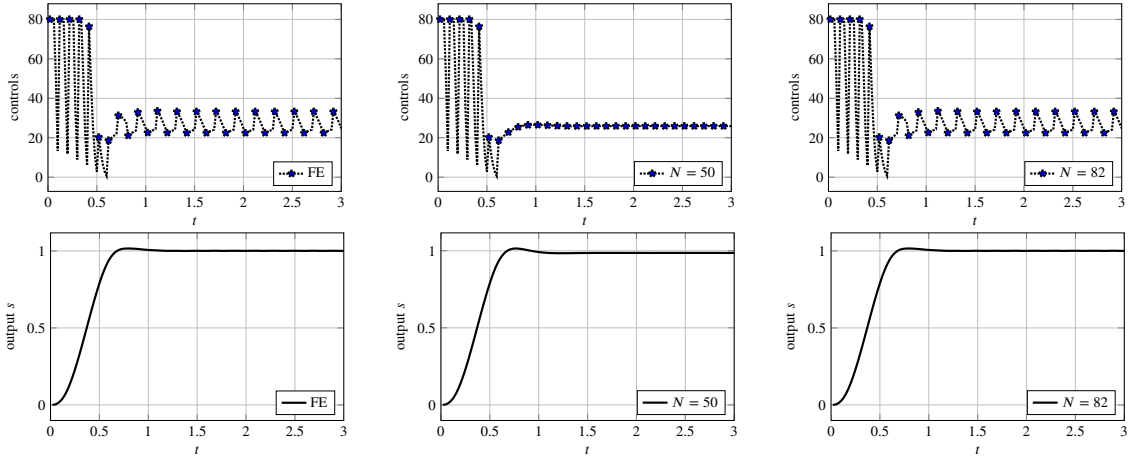
\begin{figure}
\centering
\begin{minipage}[b]{0.3 \textwidth}
\centering
	\begin{tikzpicture}[scale=0.55]
	\begin{axis}[xlabel=$t$, ylabel=controls, xmin=0, xmax=3, 
	width=1.7\textwidth,
    height=1.2\textwidth,
	xmajorgrids,
	ymajorgrids,
	yminorgrids,
	yminorticks=true, 
	axis background/.style={fill=white},
		legend pos=south east]
		\addplot+[densely dotted, color=black, line width=1.5pt, mark=*, mark repeat=5, mark phase = 1]
			table[x=t,y=u-FE] {figures/plotu_welding_DU_uc_lambda6_4_K75_MS5_fine_minh_lambda2e-06_kappa1.75_I30_ua0_ub80_omega0.2.dat};
			\addlegendentry{FE}
	\end{axis}
	\end{tikzpicture}
\end{minipage}
\begin{minipage}[b]{0.3 \textwidth}
\centering
	\begin{tikzpicture}[scale=0.55]
	\begin{axis}[xlabel=$t$, ylabel=controls, xmin=0, xmax=3, 
	width=1.7\textwidth,
    height=1.2\textwidth,
	xmajorgrids,
	ymajorgrids,
	yminorgrids,
	yminorticks=true, 
	axis background/.style={fill=white},
		legend pos= south east]
		\addplot+[densely dotted, color=black, line width=1.5pt, mark=*, mark repeat=5, mark phase = 1]
			table[x=t,y=uN-2] {figures/plotu_welding_DU_uc_lambda6_4_K75_MS5_fine_minh_lambda2e-06_kappa1.75_I30_ua0_ub80_omega0.2.dat};
			\addlegendentry{$N=50$}
	\end{axis}
	\end{tikzpicture}
\end{minipage}
\begin{minipage}[b]{0.3 \textwidth}
\centering
	\begin{tikzpicture}[scale=0.55]
	\begin{axis}[xlabel=$t$, ylabel=controls, xmin=0, xmax=3, 
	width=1.7\textwidth,
    height=1.2\textwidth,
	xmajorgrids,
	ymajorgrids,
	yminorgrids,
	yminorticks=true, 
	axis background/.style={fill=white},
		legend pos= south east]
		\addplot+[densely dotted, color=black, line width=1.5pt, mark=*, mark repeat=5, mark phase = 1]
			table[x=t,y=uN-4] {figures/plotu_welding_DU_uc_lambda6_4_K75_MS5_fine_minh_lambda2e-06_kappa1.75_I30_ua0_ub80_omega0.2.dat};
			\addlegendentry{$N=82$}
	\end{axis}
	\end{tikzpicture}
\end{minipage}
%
\begin{minipage}[b]{0.3 \textwidth}
\centering
	\begin{tikzpicture}[scale=0.55]
	\begin{axis}[xlabel=$t$, ylabel=output $s$, xmin=0, xmax=3, 
	width=1.7\textwidth,
    height=1.2\textwidth,
	xmajorgrids,
	ymajorgrids,
	yminorgrids,
	yminorticks=true, 
	axis background/.style={fill=white},
		legend pos=south east]
		\addplot+[color=black, line width=1.5pt, mark=none]
			table[x=t,y=s-FE] {figures/plotu_welding_DU_uc_lambda6_4_K75_MS5_fine_minh_lambda2e-06_kappa1.75_I30_ua0_ub80_omega0.2.dat};
			\addlegendentry{FE}
	\end{axis}
	\end{tikzpicture}
\end{minipage}
\begin{minipage}[b]{0.3 \textwidth}
\centering
	\begin{tikzpicture}[scale=0.55]
	\begin{axis}[xlabel=$t$, ylabel=output $s$, xmin=0, xmax=3, 
	width=1.7\textwidth,
    height=1.2\textwidth,
	xmajorgrids,
	ymajorgrids,
	yminorgrids,
	yminorticks=true, 
	axis background/.style={fill=white},
		legend pos= south east]
		\addplot+[color=black, line width=1.5pt,mark=none]
			table[x=t,y=s-2] {figures/plotu_welding_DU_uc_lambda6_4_K75_MS5_fine_minh_lambda2e-06_kappa1.75_I30_ua0_ub80_omega0.2.dat};
			\addlegendentry{$N=50$}
	\end{axis}
	\end{tikzpicture}
\end{minipage}
\begin{minipage}[b]{0.3 \textwidth}
\centering
	\begin{tikzpicture}[scale=0.55]
	\begin{axis}[xlabel=$t$, ylabel=output $s$, xmin=0, xmax=3, 
	width=1.7\textwidth,
    height=1.2\textwidth,
	xmajorgrids,
	ymajorgrids,
	yminorgrids,
	yminorticks=true, 
	axis background/.style={fill=white},
		legend pos= south east]
		\addplot+[color=black, line width=1.5pt,mark=none]
			table[x=t,y=s-4] {figures/plotu_welding_DU_uc_lambda6_4_K75_MS5_fine_minh_lambda2e-06_kappa1.75_I30_ua0_ub80_omega0.2.dat};
			\addlegendentry{$N=82$}
	\end{axis}
	\end{tikzpicture}
\end{minipage}
\caption{Constrained case with $\omega_{\min}=0.2$: controls and outputs for $\mu_1=1.75$, $\lambda= 2\cdot10^{-6}$, FE-MPC on the left, RB-MPC with $N = 50$ (middle) and $N = 82$ (right).}
\label{fig:WP_MinH_con_omega_mu1_1.75}
\end{figure}

We finally consider the control constrained case with lower bound $u_a\equiv 0$ and upper bound $u_b\equiv 80$. We present analogous results to the unconstrained case: the minimal stabilizing horizons for the RB-MPC and FE-MPC controller as well as the controls and outpus as a function of time  for $\omega_{\min} = 0$ are shown in Table~\ref{tab:WP_MinH_con} and Figure~\ref{fig:WP_MinH_con_mu1_1.75}. The corresponding results with a suboptimality threshold $\omega_{\min}=0.2$ are shown in Table~\ref{tab:WP_MinH_con_omega} and Figure~\ref{fig:WP_MinH_con_omega_mu1_1.75}.

We observe that the results are comparable to the unconstrained case. The minimal stabilizing horizons converge towards the reference values for all parameters and they are higher for  $\omega_{\min}=0.2$ compared to the setting with no threshold. The smoothing effect of $\omega_{\min}$ is slightly more evident in the constrained case especially if the output is close to the desired output. Overall, the performance of the RB-MPC approach is very close to the original MPC approach. We would like to note, however, that the feasibility of the stability condition $\tilde{\omega}_{N,K}>0$ for the constrained case, see Proposition \ref{prop:Alpha_C}, strongly depends on the ratio of the cost for the constrained and unconstrained problem. There are certainly examples conceivable in which the cost of the former will significantly exceed the cost of the latter and our approach would thus fail.


\begin{table}[t]
\begin{center}
\pgfplotstableread{tables/data_welding_DU_uc_lambda6_4_K75_MS5_fine_minh_lambda2e-06_I30_ua0_ub80_omega0.0.dat}\dataA
\pgfplotstableread{tables/data_welding_DU_uc_lambda6_4_K75_MS5_fine_minh_lambda5e-06_I30_ua0_ub80_omega0.0.dat}\dataB
\pgfplotstableread{tables/data_welding_DU_uc_lambda6_4_K75_MS5_fine_minh_lambda8e-06_I30_ua0_ub80_omega0.0.dat}\dataC
\pgfplotstableread{tables/data_welding_DU_uc_lambda6_4_K75_MS5_fine_minh_lambda5e-05_I30_ua0_ub80_omega0.0.dat}\dataD

\pgfplotstablecreatecol[copy column from table={\dataB}{[index] 1}] {56ave1} {\dataA}
\pgfplotstablecreatecol[copy column from table={\dataB}{[index] 4}] {56ave2} {\dataA}

\pgfplotstablecreatecol[copy column from table={\dataC}{[index] 1}] {86ave1} {\dataA}
\pgfplotstablecreatecol[copy column from table={\dataC}{[index] 4}] {86ave2} {\dataA}

\pgfplotstablecreatecol[copy column from table={\dataD}{[index] 1}] {55ave1} {\dataA}
\pgfplotstablecreatecol[copy column from table={\dataD}{[index] 4}] {55ave2} {\dataA}

	\pgfplotstabletypeset[font=\footnotesize, 
		every row no 4/.style={before row=\hline \hline},
		every head row/.style={
			before row={
				 $\lambda$  & \multicolumn{2}{c}{$2\cdot10^{-6}$ } & \multicolumn{2}{c}{$5\cdot10^{-6}$ }  & \multicolumn{2}{c}{$8\cdot10^{-6}$ } & \multicolumn{2}{c}{$5\cdot10^{-5}$ }\\
				\hline & &  & &  & &  & &  \\[-7pt]
				$\mu_1$ & 0.75& 1.75  & 0.75& 1.75 & 0.75& 1.75  & 0.75& 1.75\\
				\hline & &   & &  & &  & &   \\[-7pt]
				$N$ & \multicolumn{2}{c|}{$K_{\textrm{ave}}$} & \multicolumn{2}{c|}{$K_{\textrm{ave}}$} & \multicolumn{2}{c|}{$K_{\textrm{ave}}$} & \multicolumn{2}{c}{$K_{\textrm{ave}}$} \\
				\hline & &  & &  & &  & &    \\[-7pt]
			},
			output empty row
		},
		columns={N, K-ave1, K-ave2, 56ave1, 56ave2, 86ave1, 86ave2, 55ave1, 55ave2},
        columns/N/.append style={string replace={-1}{}},
		columns/K-ave2/.append style={column type={c|}},
		columns/56ave2/.append style={column type={c|}},
		columns/86ave2/.append style={column type={c|}},
        empty cells with={$\kappa$}
        ]{\dataA}
\end{center} 
\caption{Constrained case: average minimal stabilizing horizons over $[0,3]$.}
\label{tab:WP_MinH_con}
\end{table}

\begin{table}[t]
\begin{center}
\pgfplotstableread{tables/data_welding_DU_uc_lambda6_4_K75_MS5_fine_minh_lambda2e-06_I30_ua0_ub80_omega0.2.dat}\dataA
\pgfplotstableread{tables/data_welding_DU_uc_lambda6_4_K75_MS5_fine_minh_lambda5e-06_I30_ua0_ub80_omega0.2.dat}\dataB
\pgfplotstableread{tables/data_welding_DU_uc_lambda6_4_K75_MS5_fine_minh_lambda8e-06_I30_ua0_ub80_omega0.2.dat}\dataC
\pgfplotstableread{tables/data_welding_DU_uc_lambda6_4_K75_MS5_fine_minh_lambda5e-05_I30_ua0_ub80_omega0.2.dat}\dataD

\pgfplotstablecreatecol[copy column from table={\dataB}{[index] 1}] {56ave1} {\dataA}
\pgfplotstablecreatecol[copy column from table={\dataB}{[index] 4}] {56ave2} {\dataA}

\pgfplotstablecreatecol[copy column from table={\dataC}{[index] 1}] {86ave1} {\dataA}
\pgfplotstablecreatecol[copy column from table={\dataC}{[index] 4}] {86ave2} {\dataA}

\pgfplotstablecreatecol[copy column from table={\dataD}{[index] 1}] {55ave1} {\dataA}
\pgfplotstablecreatecol[copy column from table={\dataD}{[index] 4}] {55ave2} {\dataA}

	\pgfplotstabletypeset[font=\footnotesize, 
		every row no 4/.style={before row=\hline \hline},
		every head row/.style={
			before row={
				 $\lambda$  & \multicolumn{2}{c}{$2\cdot10^{-6}$ } & \multicolumn{2}{c}{$5\cdot10^{-6}$ }  & \multicolumn{2}{c}{$8\cdot10^{-6}$ } & \multicolumn{2}{c}{$5\cdot10^{-5}$ }\\
				\hline & &  & &  & &  & &  \\[-7pt]
				$\mu_1$ & 0.75& 1.75  & 0.75& 1.75 & 0.75& 1.75  & 0.75& 1.75\\
				\hline & &   & &  & &  & &   \\[-7pt]
				$N$ & \multicolumn{2}{c|}{$K_{\textrm{ave}}$} & \multicolumn{2}{c|}{$K_{\textrm{ave}}$} & \multicolumn{2}{c|}{$K_{\textrm{ave}}$} & \multicolumn{2}{c}{$K_{\textrm{ave}}$} \\
				\hline & &  & &  & &  & &    \\[-7pt]
			},
			output empty row
		},
		columns={N, K-ave1, K-ave2, 56ave1, 56ave2, 86ave1, 86ave2, 55ave1, 55ave2},
        columns/N/.append style={string replace={-1}{}},
		columns/K-ave2/.append style={column type={c|}},
		columns/56ave2/.append style={column type={c|}},
		columns/86ave2/.append style={column type={c|}},
        empty cells with={$\kappa$}
        ]{\dataA}
\end{center} 
\caption{Constrained case: average minimal stabilizing horizons over $[0,3]$ with threshold $\omega_{\min}=0.2$.}
\label{tab:WP_MinH_con_omega}
\end{table}

\subsubsection{Time}
\label{sec:NumExp_WP_Time}

We finally report on the online computational times to evaluate the FE-MPC and RB-MPC feedback laws and associated suboptimality degrees. In the FE case, we include the time to compute the feedback law $\kappa$ in~\eqref{eq:uast} which corresponds to solving the optimality systems~\eqref{eq:FE_OCPSystem} and the time to evaluate the suboptimality degree $\omega(\mu)$ from \eqref{eq:FE_RelaxedDPP_Lyapunov}. In the RB case, we include the time to compute the RB-MPC feedback law~\eqref{eq:kappa_N} and the time to evaluate the suboptimality degree $\tilde{\omega}_{N,K}(\mu)$ from~\eqref{eq:omegatilde_uc} (resp.\ \eqref{eq:omegatilde_c} in the constrained case). Note that evaluating the suboptimality degree in the constrained case requires an additional solution of the unconstrained RB optimal control problem due to the bound in~\eqref{eq:Alpha_C_add_est}. The average computational times in seconds are presented for various values of the RB dimension $N$ and the prediction horizon $K$ in Table~\ref{tab:WP_time_UC} and Table~\ref{tab:WP_time_Con} for the unconstrained and constrained case, respectively.  We used MATLAB for all computations on a computer with a $3.30$ GHz Intel Core i5-6600 processor, where only three of the four cores were available to the system. The average is taken over the parameter test set $\Xi_{test}$ of size $60$ introduced in Section \ref{sec:NumExp_WP_Effectivity}.

As expected, the RB times increase with $N$ and $K$ and the FE time also increases with $K$. 
For $N = 66$ and $K = 10$ (say), the speed-up gained by the RB approach is approximately a factor of 100. We also observe that the speed-up remains almost the same for different prediction lengths $K$, the influence of $N$ on the speed-up is obviously more pronounced. We note that the speed-up of the RB vs. FE-MPC controller over the whole time interval would approximately be the same order of magnitude, assuming similar prediction lengths in the RB and FE case (which is guaranteed as long as $N$ is large enough as seen in the last section). To summarize, the RB-MPC approach provides a significant speed-up compared to the classical (FE-)MPC approach while guaranteeing stability of the closed-loop system.

\begin{table}[t]
\begin{center}
\pgfplotstableread{tables/data_welding_DU_uc_lambda6_4_K75_MS5_fine_time.txt}\dataA

\pgfplotstabletypeset[font=\footnotesize, 
		every row no 4/.style={before row=\hline \hline},
		every head row/.style={
			before row={
				\hline & &  & &  &   \\[-7pt]
				\backslashbox{$N$}{$K$} & $5$ & $10$& $20$ & $40$ & $75$\\
				\hline & &  & &  &   \\[-7pt]
			},
			output empty row
		},
		columns={N, K5, K10, K20, K40, K75},
        columns/N/.append style={column type={|c|}, string replace={-1}{}},
        columns/K5/.style={sci,	sci zerofill},
        columns/K10/.style={sci, sci zerofill},
        columns/K20/.style={sci, sci zerofill},
        columns/K40/.style={sci, sci zerofill},
        columns/K75/.style={column type={c|},sci, sci zerofill},
        empty cells with={FE}
        ]{\dataA}
\end{center} 
\caption{Unconstrained case: average online computational time in sec.}
\label{tab:WP_time_UC}
\end{table}

\begin{table}[t]
\begin{center}
\pgfplotstableread{tables/data_welding_DU_uc_lambda6_4_K75_MS5_fine_time_ua0_ub80.txt}\dataA

\pgfplotstabletypeset[font=\footnotesize, 
		every row no 4/.style={before row=\hline \hline},
		every head row/.style={
			before row={
			    \hline & &  & &  &   \\[-7pt]
				\backslashbox{$N$}{$K$} & $5$ & $10$& $20$ & $40$ & $75$\\
				\hline & &  & &  &   \\[-7pt]
			},
			output empty row
		},
		columns={N, K5, K10, K20, K40, K75},
        columns/N/.append style={column type={|c|}, string replace={-1}{}},
        columns/K5/.style={sci,	sci zerofill},
        columns/K10/.style={sci, sci zerofill},
        columns/K20/.style={sci, sci zerofill},
        columns/K40/.style={sci, sci zerofill},
        columns/K75/.style={column type={c|},sci, sci zerofill},
        empty cells with={FE}
        ]{\dataA}
\end{center} 
\caption{Constrained case: average online computational time in sec.}
\label{tab:WP_time_Con}
\end{table}

\section{Conclusion}

The goal of this work is to improve the applicability of MPC for high dimensional systems arising e.g.\ through the discretization of partial differential equations. To this end, we proposed an RB-MPC controller for parametrized parabolic partial differential equations. The key ingredients are reduced basis approximations and {\it a posteriori} error bounds which can be evaluated efficiently in the online stage. This not only allows to compute the feedback control efficiently but also to assess the stability of the closed-loop system during the online stage. We also proposed an adaptive algorithm for choosing the minimal stabilizing prediction horizon in the MPC loop. Furthermore, we extended the theory to weakly coercive problems and numerically investigated the effect of a suboptimality threshold combined with a multi-step approach.

Despite the pretty conservative error bounds, the RB-MPC approach worked very well for the two numerical examples considered in this paper. The minimum stabilizing prediction horizon of the RB-MPC approach is very close to the minimum stabilizing horizon of the reference FE-MPC approach as long as the RB approximation is sufficiently accurate. At the same time, the RB-MPC approach delivers a substantial speed-up compared to the FE-MPC approach thus improving the applicability of MPC for high-dimensional systems.

We are currently extending the methodolgy presented here to incorporate an $\varepsilon$-neighborhood around the desired state following the ideas on ``practical optimality.'' Another obvious next step is the extension of the methodology presented here for state feedback MPC to output feedback MPC. In actual practice -- and especially when considering PDEs -- the entire state is often unknown and only a (small) number of measurements of the state are available. Combining the presented approach with a suitable observer is one topic for future work. Although our approach cannot be applied or extended directly to nonlinear problems, another topic is to investigate how it may provide guidelines even for MPC of nonlinear problems. 

\begin{appendix}
\section{Proofs}

\subsection{Proof of Lemma \ref{lem:Hilfslem_OptStateErr}}
\begin{proof}
Choosing $\phi= e_y^k$ in the state error residual equation 
\begin{align} \label{eq:StateResEq}
  m(e_y^{k}-e_y^{k-1},\phi) + \tau a(e_y^{k},\phi;\mu) - \tau b(e_u^{k},\phi) &= \tau  r_y^k (\phi;\mu)  , \quad \forall \phi\in Y, \ \forall k\in \K,
  \end{align}
 with $  e_y^{0} = y_0-y_{N,0} $ we obtain
\begin{align}\label{eq:Delta_u_6ErrRes}
 m(e_y^k  , e_y^k) +\tau a(e_y^k, e_y^k ) =  m(e_y^{k-1} , e_y^k) + \tau \sum_{i=1}^m b_i(e_y^k) ( e_u^k )_i +  \tau  r_y^k(e_y^k;\mu) .
\end{align}
From Cauchy-Schwarz and Young's inequality ($2ab \le a^2 + b^2$) it follows that
\begin{align}\label{eq:Delta_u_Y1}
2 m(e_y^{k-1} , e_y^k) \le m(e_y^{k} , e_y^k) + m(e_y^{k-1} , e_y^{k-1}),
\end{align}
and from boundedness of the linear functionals $b_i$, Young's inequality, and Cauchy-Schwarz, we get
\begin{equation}
\begin{aligned}
2 \tau  \sum_{i=1}^m b_i(e_y^k) ( e_u^k )_i &\le 2 \tau  \sum_{i=1}^m \|b_i\|_{Y'} \|e_y^k\|_Y ( e_u^k )_i \frac{(2\alphaLB)^{1/2}}{(2\alphaLB)^{1/2}} 
\le \dfrac{2 \tau}{\alphaLB}  \left( \sum_{i=1}^m \|b_i\|_{Y'}  ( e_u^k )_i \right)^2 + \dfrac{\tau \alphaLB}{2}\|e_y^k\|_Y^2\\
&\le \dfrac{2 \tau}{\alphaLB}  \left[ \left( \sum_{i=1}^m \|b_i\|_{Y'}^2 \right)^{\frac{1}{2}}  \left( \sum_{i=1}^m ( e_u^k )_i^2 \right)^{\frac{1}{2}} \right]^2 + \dfrac{\tau \alphaLB}{2}\|e_y^k\|_Y^2 \\
&= \dfrac{2\tau}{\alphaLB} \left( \sum_{i=1}^m \|b_i\|_{Y'}^2 \right) \| e_u^k \|_{R^m}^2  + \dfrac{\tau \alphaLB}{2}\|e_y^k\|_Y^2 .
\end{aligned}
\end{equation}
The definition of the dual norm and Young's inequality also leads to
\begin{align}\label{eq:Delta_u_Y3}
2\tau r_y^k(e_y^k;\mu) \le 2\tau \|r_y^k(\cdot;\mu) \|_{Y'} \|e_y^k\|_Y 
\le  \dfrac{2 \tau}{\alphaLB}  \|r_y^k(\cdot;\mu) \|_{Y'}^2 + \dfrac{\tau \alphaLB}{2} \|e_y^k\|_Y^2.
\end{align}
It thus follows from \eqref{eq:Delta_u_6ErrRes} - \eqref{eq:Delta_u_Y3}, invoking \eqref{eq:truth_ellipticity} and \eqref{eq:Alpha_LB}, and summing from $k=1$ to $K$ that
\begin{align}
m(e_y^K  , e_y^K) +\tau  \sum_{k=1}^K a(e_y^k, e_y^k ) 
		\le  \dfrac{2 \tau}{\alphaLB}   \sum_{k=1}^K \|r_y^k(\cdot;\mu) \|_{Y'}^2 +  \dfrac{2}{\alphaLB} \left( \sum_{i=1}^m \|b_i\|_{Y'}^2 \right)  \tau \sum_{k=1}^K \| e_u^k \|_{R^m}^2  + m(e_y^0,e_y^0).
\end{align}
From the previous equation we can also derive that
\begin{equation}
\begin{aligned}
\tau  \sum_{k=1}^K \| e_y^k \|_Y^2 &\le \dfrac{\tau}{\alphaLB} \sum_{k=1}^K a(e_y^k,e_y^k) 
\le \dfrac{1}{\alphaLB} \left( m(e_y^k,e_y^k)  + \sum_{k=1}^K \tau a(e_y^k,e_y^k) \right) \\
& \le \dfrac{2 \tau}{\alphaLB^2}   \sum_{k=1}^K \|r_y^k(\cdot;\mu) \|_{Y'}^2 +  \dfrac{2}{\alphaLB^2} \left( \sum_{i=1}^m \|b_i\|_{Y'}^2 \right)  \tau \sum_{k=1}^K \| e_u^k \|_{R^m}^2 + \dfrac{\norm{e^0_y}_{L^2(\Omega)}^2}{\alphaLB}  .
\end{aligned}
\end{equation}
\end{proof}

\subsection{Proof of Lemma \ref{lem:Hilfslem_OptAdjErr}}
\begin{proof}
We choose $\varphi= e_p^{k}$ in the adjoint error residual equation
\begin{align} \label{eq:AdjResEq}
m(\varphi, e_p^{k}-e_p^{k+1}) + \tau a(\varphi, e_p^{k};\mu) + \tau  \sigma_1( e_y^{k},   \varphi)_D&=  \tau  r_p^k (\varphi;\mu) , \quad \forall \varphi\in Y, \ \forall k\in \K \backslash \{K\},
\end{align}
 and $\varphi= e_p^{K}$ in the finial condition 
 \begin{align} \label{eq:AdjResEqFinal}
m(\varphi, e_p^{K}) + \tau a(\varphi, e_p^{K};\mu) + \tau  \sigma_1( e_y^{K},   \varphi)_D  +\sigma_2( e_y^{K},   \varphi)_D &= \tau r_p^K (\varphi;\mu)  , \quad \forall \varphi\in Y,
\end{align} 
to get
\begin{align}\label{eq:Delta_u_7AdjRes}
m(e_p^{k}, e_p^{k}) + \tau a(e_p^{k}, e_p^{k};\mu) &= m(e_p^{k}, e_p^{k+1}) + \tau  r_p^k (e_p^{k};\mu) - \tau  \sigma_1( e_y^{k},   e_p^{k})_D , 
\end{align}
for all $k\in \K \backslash \{K\}  $ and
\begin{align}\label{eq:Delta_u_7AdjResF}
m(e_p^{K}, e_p^{K}) + \tau a(e_p^{K}, e_p^{K};\mu)  &= \tau r_p^K (e_p^{K};\mu) - \tau  \sigma_1( e_y^{K},   e_p^{K})_D  -\sigma_2( e_y^{K},   e_p^{K})_D.
\end{align}
From Cauchy-Schwarz and Young's Inequality it follows that
\begin{align}\label{eq:Delta_u_P1}
2 m(e_p^{k-1} , e_p^k) \le m(e_p^{k} , e_p^k) + m(e_p^{k-1} , e_p^{k-1}),
\end{align}
and from the definition of the dual norm and Young's Inequality 
\begin{align}\label{eq:Delta_u_P2}
2\tau r_p^k(e_p^k;\mu) \le 2\tau \|r_p^k(\cdot;\mu) \|_{Y'} \|e_p^k\|_Y 
\le  \dfrac{2 \tau}{\alphaLB}  \|r_p^k(\cdot;\mu) \|_{Y'}^2 + \dfrac{\tau \alphaLB}{2} \|e_p^k\|_Y^2  .
\end{align}
We again use Cauchy-Schwarz, Young's Inequality, and the definition of the constant $C_D$ defined in \eqref{eq:ConstantC_D} to obtain
\begin{equation}
\begin{aligned}\label{eq:Delta_u_P3}
2 \tau \sigma_1 (e_y^k, e_p^k)_D &\le 2 \tau \sigma_1 (e_y^k, e_y^k)_D^{\frac{1}{2}} (e_p^k, e_p^k)_D^{\frac{1}{2}}
		\le 2 \tau \sigma_1 (e_y^k, e_y^k)_D^{\frac{1}{2}} C_D \norm{e_p^k}_Y \\
		 &\le \frac{2C_D^2 \sigma_1^2 \tau}{\alphaLB}(e_y^k, e_y^k)_D + \frac{\tau \alphaLB}{2} \norm{e_p^k}_Y ^2  .
\end{aligned}
\end{equation}
Similarly, we get
\begin{align}\label{eq:Delta_u_P4}
2 \sigma_2 (e_y^K, e_p^K)_D &\le 2  \sigma_2 (e_y^K, e_y^K)_D^{\frac{1}{2}} (e_p^K, e_p^K)_D^{\frac{1}{2}}
		\le   \sigma_2^2 (e_y^K, e_y^K)_D  + (e_p^K, e_p^K)_D.
\end{align}
It then follows from~\eqref{eq:Delta_u_P2} - \eqref{eq:Delta_u_P4}, summing \eqref{eq:Delta_u_7AdjRes} from $k=1$ to $K-1$, adding \eqref{eq:Delta_u_7AdjResF}, 
invoking \eqref{eq:truth_ellipticity} and \eqref{eq:Alpha_LB} that
\begin{align}
m(e_p^{1}, e_p^{1}) + \tau \sum_{k=1}^K a(e_p^{k}, e_p^{k};\mu) 
\le  \dfrac{2 \tau}{\alphaLB} \sum_{k=1}^K \|r_p^k(\cdot;\mu) \|_{Y'}^2 + \frac{2C_D^2 \sigma_1^2 \tau}{\alphaLB}\sum_{k=1}^K (e_y^k, e_y^k)_D  +  \sigma_2^2 (e_y^K, e_y^K)_D,
\end{align}
where we also used the fact that $(\cdot,\cdot)_D\le m(\cdot,\cdot)$. From the previous equation we directly obtain
\begin{align}
\tau  \sum_{k=1}^K \|e_p^k \|_{Y}^2
		\le  \dfrac{2 \tau}{\alphaLB^2}   \sum_{k=1}^K \|r_p^k(\cdot;\mu) \|_{Y'}^2 +  \dfrac{2C_D^2\sigma_1^2 \tau}{\alphaLB^2} \sum_{k=1}^K (e_y^k, e_y^k)_D + \dfrac{\sigma_2^2}{\alphaLB} (e_y^K, e_y^K)_D. 
\end{align}
\end{proof}

\subsection{Proof of Proposition \ref{prop:OptControlErrBound}}
\begin{proof}
To begin, we recall the control residual equation
\begin{align}\label{eq:ControlResEq}
\left(\lambda e_u^{k} -(\cB^Te_p)^k, \psi \right)_{\xR^m} = r_u^k (\psi;\mu), \quad \forall \psi\in U, \ \forall k\in \K.
\end{align}
Summing up \eqref{eq:ControlResEq} multiplied by $\tau$ from $k=1$ to $K$ with $\psi = e_u^{k}$, adding the sum of the adjoint \eqref{eq:AdjResEq} from $k=1$ to $K-1$ with $\varphi = e_y^{k}$ and the final condition \eqref{eq:AdjResEqFinal} with $\varphi = e_y^{K}$, and subtracting the sum of the state \eqref{eq:StateResEq} with $\phi= e_p^{k}$ from $k=1$ to $K$, we obtain 
\begin{multline}
\label{eq:ProofDeltaALT_Step5}
\lambda \tau \sum_{k=1}^K (e_u^{k}, e_u^{k}  )_{\xR^m} - \tau \sum_{k=1}^K ((\cB^Te_p)^k, e_u^{k} )_{\xR^m} 
 + m(e_y^{K}, e_p^{K})  + \tau a(e_y^{K}, e_p^{K})  + \tau \sigma_1 (e_y^{K}, e_y^{K})_D + \sigma_2 (e_y^{K}, e_y^{K})_D  \\
 + \sum_{k=1}^{K-1}  m(e_y^{k}, e_p^{k} -e_p^{k+1} ) + \cby{\tau}\sum_{k=1}^{K-1}  a(e_y^{k}, e_p^{k} ) + \sum_{k=1}^{K-1}  \tau \sigma_1 (e_y^{k}, e_y^{k})_D 
- \sum_{k=1}^K m(e_y^{k} -e_y^{k-1} , e_p^{k}) - \cby{\tau} \sum_{k=1}^{K}  a(e_y^{k}, e_p^{k} ) 
+ \tau \sum_{k=1}^K ((\cB^Te_p)^k, e_u^{k} )_{\xR^m} \\
\le 0 + \tau \sum_{k=1}^K r_p^k(e_y^{k}) - \tau \sum_{k=1}^K r_y^k(e_p^{k}) .
\end{multline}
Here, we also used the fact that 
\begin{align}
\left(\lambda e_u^{k} -(\cB^Te_p)^k, e_u^{k} \right)_{\xR^m}\le  0,
\end{align}
which directly follows from \eqref{eq:FE_OCPOptimalityInequality} and \eqref{eq:RB_OCPOptimalityInequality}  using the fact that both controls $u^\ast$ and $ u_N^\ast$ are admissible, i.e, $u^\ast, u_N^\ast\in U_{ad}$. 
We observe that \eqref{eq:ProofDeltaALT_Step5} is equivalent to 
\begin{align}
&\lambda \tau \sum_{k=1}^K (e_u^k, e_u^k  )_{\xR^m} +\sum_{k=1}^{K}  \tau \sigma_1 (e_y^k, e_y^k)_D  + \sigma_2 (e_y^K, e_y^K)_D  
\le \tau \sum_{k=1}^K r_p^k(e_y^k) - \tau \sum_{k=1}^K r_y^k(e_p^k) \cby{+} m(e_y^0, e_p^1).
\end{align}
Next, we apply the Cauchy-Schwarz inequality and obtain
\begin{equation}
\begin{aligned}\label{eq:Step5}
&\lambda \tau \sum_{k=1}^K \|e_u^k\|_U^2 +\sum_{k=1}^{K}  \tau \sigma_1 |e_y^k|_D^2  + \sigma_2 |e_y^K|_D^2 \\  
&\le \left(\tau \sum_{k=1}^K \|r_p^k(\cdot)\|_{Y'}^2 \right)^\frac{1}{2} \underbrace{\left(\tau \sum_{k=1}^K \|e_y^k\|_{Y}^2 \right)^\frac{1}{2} }_{(I)}
+  \left(\tau \sum_{k=1}^K \|r_y^k(\cdot)\|_{Y'}^2 \right)^\frac{1}{2} \underbrace{\left(\tau \sum_{k=1}^K \|e_p^k\|_{Y}^2 \right)^\frac{1}{2}}_{(II)} 
+\underbrace{\norm{e_y^0}_{L^2(\Omega)}\norm{e_p^1}_{L^2(\Omega)}}_{(III)}  .
\end{aligned}
\end{equation}
We now apply Lemma \ref{lem:Hilfslem_OptStateErr} and Lemma \ref{lem:Hilfslem_OptAdjErr} to $(I)$ and $(II)$, respectively, and exploit that $a^2 + b^2 + c^2\le (a+b+c)^2$ for $a,b,c \ge 0$. This gives us  
\begin{align}\label{eq:StateOptErrBoundsqrt_ProofDeltaALT}
\sqrt{\tau  \sum_{k=1}^K \| e_y^k \|_Y^2  }  \le  
	\sqrt{\dfrac{2 \tau}{\alphaLB^2}   \sum_{k=1}^K \|r_y^k(\cdot) \|_{Y'}^2 }
		+ \sqrt{ \dfrac{2}{\alphaLB^2} \left( \sum_{i=1}^m \|b_i\|_{Y'}^2 \right) \tau \sum_{k=1}^K \| e_u^k \|_U^2}
+ \sqrt{\dfrac{\norm{e^0_y}_{L^2(\Omega)}^2}{\alphaLB}  }.
\end{align}
and
\begin{align}\label{eq:AdjOptErrBoundsqrt_ProofDeltaALT}
\sqrt{\tau  \sum_{k=1}^K \|e_p^k \|_{Y}^2}
		\le \sqrt{ \dfrac{2 \tau}{\alphaLB^2}   \sum_{k=1}^K \|r_p^k(\cdot) \|_{Y'}^2 }+  \sqrt{\dfrac{2C_D^2\sigma_1^2 \tau}{\alphaLB^2} \sum_{k=1}^K (e_y^k, e_y^k)_D }+ \sqrt{\dfrac{\sigma_2^2}{\alphaLB} (e_y^K, e_y^K)_D }.
\end{align}
For the last term $(III)$ we apply Lemma \ref{lem:Hilfslem_OptAdjErr} again, since $\norm{e_p^1}_{L^2(\Omega)} = m(e_p^1, e_p^1) \le \interleave e_p^1\interleave_p$. 
We define the following shortcuts to simplify the representation
\begin{align}
 E_u = \left( \tau \sum_{k=1}^K  \| e_u^{k}\|^2_{\xR^m}\right)^{\frac{1}{2}},
 \quad R_0  = \norm{e^0_y}_{L^2(\Omega)}, 
\quad R_\bullet = \left( \tau \sum_{k=1}^K  \| r_\bullet^{k} (\cdot;\mu) \|^2_{Y'}\right)^{\frac{1}{2}}, \ \bullet \in \{y,p\}.
\end{align}
Plugging \eqref{eq:StateOptErrBoundsqrt_ProofDeltaALT} and \eqref{eq:AdjOptErrBoundsqrt_ProofDeltaALT} into \eqref{eq:Step5},  using \eqref{eq:AdjOptErrBound_ProofDeltaALT_Enorm} in term $(III)$, and $(v,v)_D = |v|_D^2$ we obtain 
\begin{multline}\label{eq:Step8}
\lambda E_u^2 + \sigma_1 \tau \sum_{k=1}^K |e_y^k|^2_D + \sigma_2 |e_y^K|^2_D
\le  
R_p \left(\dfrac{\sqrt{2}}{\alphaLB}   R_y + \dfrac{\sqrt{2}}{\alphaLB} \left( \sum_{i=1}^m \|b_i\|_{Y'}^2 \right)^{\frac{1}{2}} E_u  
+ \dfrac{R_0}{\sqrt{\alphaLB}  }\right) \\
+ R_y \left( \dfrac{\sqrt{ 2}}{\alphaLB}   R_p +  \sqrt{\dfrac{2C_D^2\sigma_1}{\alphaLB^2}} \sqrt{ \sigma_1 \tau \sum_{k=1}^K |e_y^k|_D^2}+ \sqrt{\dfrac{\sigma_2}{\alphaLB}} \sqrt{ \sigma_2|e_y^K|_D^2 } \right)  \\
+ R_0 \left( \sqrt{\dfrac{ 2}{\alphaLB}}   R_p +  \sqrt{\dfrac{2C_D^2\sigma_1}{\alphaLB}} \sqrt{ \sigma_1 \tau \sum_{k=1}^K |e_y^k|_D^2}+ \sqrt{\sigma_2} \sqrt{ \sigma_2|e_y^K|_D^2 } \right)   .
\end{multline}
We use Young's inequality twice in the second term on the right hand side as follows
\begin{align}
R_y \sqrt{\dfrac{2C_D^2\sigma_1}{\alphaLB^2}} \sqrt{ \sigma_1 \tau \sum_{k=1}^K |e_y^k|_D^2}  
\le \dfrac{C_D^2\sigma_1}{\alphaLB^2} R_y^2 + \dfrac{\sigma_1}{2} \tau \sum_{k=1}^K |e_y^k|_D^2
\end{align}
and 
\begin{align}
R_y\sqrt{\dfrac{\sigma_2}{\alphaLB}} \sqrt{ \sigma_2|e_y^K|_D^2 } 
 \le \dfrac{\sigma_2}{2\alphaLB} R_y^2 + \dfrac{\sigma_2}{2}|e_y^K|_D^2.
\end{align}
In addition, we invoke Young's inequality also twice in the third term to get
\begin{align}
R_0 \sqrt{\dfrac{2C_D^2\sigma_1}{\alphaLB}} \sqrt{ \sigma_1 \tau \sum_{k=1}^K |e_y^k|_D^2}  
\le \dfrac{C_D^2\sigma_1}{\alphaLB} R_0^2 + \dfrac{\sigma_1}{2} \tau \sum_{k=1}^K |e_y^k|_D^2
\end{align}
and 
\begin{align}
\sqrt{\sigma_2}R_0 \sqrt{ \sigma_2|e_y^K|_D^2 } 
 \le \dfrac{\sigma_2}{2} R_0^2 + \dfrac{\sigma_2}{2}|e_y^K|_D^2
\end{align}
such that the second and third term in \eqref{eq:Step8} cancel.
We now have 
\begin{equation}
\begin{aligned}
\lambda E_u^2 \le& 
\dfrac{\sqrt{2}}{\alphaLB} R_p   R_y + \dfrac{\sqrt{2}}{\alphaLB} \left( \sum_{i=1}^m \|b_i\|_{Y'}^2 \right)^{\frac{1}{2}} R_p  E_u 
+ R_pR_0\dfrac{1}{\sqrt{\alphaLB}  }
+\dfrac{\sqrt{2}}{\alphaLB} R_p   R_y  + \dfrac{C_D^2\sigma_1}{\alphaLB^2} R_y^2 + \dfrac{\sigma_2}{2\alphaLB} R_y^2 \\
&+ \sqrt{\dfrac{ 2}{\alphaLB}} R_0 R_p  
			+\dfrac{C_D^2\sigma_1}{\alphaLB} R_0^2 
			+ \dfrac{\sigma_2}{2} R_0^2\\
=& 		 \dfrac{\sqrt{2}}{\alphaLB} \left( \sum_{i=1}^m \|b_i\|_{Y'}^2 \right)^{\frac{1}{2}} R_p  E_u
+ \left( \dfrac{2\sqrt{2}}{\alphaLB}R_y 
+ \dfrac{1+\sqrt{2}}{\sqrt{\alphaLB}} R_0 \right)R_p    
 + \left( \dfrac{C_D^2\sigma_1}{\alphaLB} + \dfrac{\sigma_2}{2}\right) R_0^2  \\
&+ \left( \dfrac{C_D^2\sigma_1}{\alphaLB^2} + \dfrac{\sigma_2}{2\alphaLB}\right) R_y^2  .
\end{aligned}
\end{equation}
We thus obtain the quadratic inequality for $E_u$ given by
\begin{align}
 E_u^2 -2c_1(\mu) E_u -c_2(\mu)\le 0, 
\end{align}
where 
\begin{align}
c_1(\mu) &=  \dfrac{1}{\sqrt{2}\alphaLB\lambda} \left( \sum_{i=1}^m \|b_i\|_{Y'}^2 \right)^{\frac{1}{2}} R_p  ,\\
c_2(\mu) &= \dfrac{1}{\lambda}\left[ \left(
 \dfrac{2\sqrt{2}}{\alphaLB}R_y 
+ \dfrac{1+\sqrt{2}}{\sqrt{\alphaLB}} R_0 \right)R_p    
+   \left( \dfrac{C_D^2\sigma_1}{\alphaLB} + \dfrac{\sigma_2}{2}\right) R_0^2  
+ \left( \dfrac{C_D^2\sigma_1}{\alphaLB^2} + \dfrac{\sigma_2}{2\alphaLB}\right) R_y^2 
\right].
\end{align}
The inequality is satisfied if and only if 
\begin{align}
\Delta^-_N \le E_u \le \Delta^+_N ,\ \textrm{where } 
\Delta^{\pm}_N = c_1(\mu) \pm \sqrt{c_1(\mu)^2 +c_2(\mu)} .
\end{align} 
\end{proof}

\end{appendix}

\bibliographystyle{abbrv}
\bibliography{literature}

\end{document}